\documentclass[10pt]{article}

\usepackage[utf8]{inputenc}
\usepackage[dvipsnames]{xcolor}
\usepackage{cancel}
\usepackage{leftidx}
\usepackage{amsmath}
\usepackage{amscd}
\usepackage{amssymb}
\usepackage{rotating}
\usepackage{amsfonts}
\usepackage{mathrsfs}
\usepackage{amsthm}
\usepackage{verbatim}
\usepackage{nicematrix}
\usepackage{lscape}

\usepackage{cases}
\usepackage[top=1in, bottom=1in, left=1.25in, right=1.25in]{geometry}
\usepackage[colorlinks,linkcolor=Cyan, anchorcolor=Cyan,urlcolor=blue,  citecolor=blue]{hyperref}
\usepackage{cleveref} 
\usepackage[hyperpageref]{backref}

\usepackage{titlesec}
       \titleformat{\chapter}[display]
             {\normalfont\Large\bfseries}{\thechapter}{11pt}{\Large}
       \titleformat{\section}
             {\normalfont\large\bfseries}{\thesection}{11pt}{\large}
       \titlespacing*{\chapter}{0pt}{0pt}{15pt} 
       \titlespacing*{\section}{0pt}{3.5ex plus 1ex minus .2ex}{2.3ex plus .2ex}

    \usepackage[titletoc]{appendix}

\input xy
\xyoption{all}

\allowdisplaybreaks[4]

\newcommand{\rvline}{\hspace*{-\arraycolsep}\vline\hspace*{-\arraycolsep}}

\newcommand{\nn}{\nonumber}
\newcommand{\sfh}{\mathsf{h}}
\newcommand{\tsfh}{\tilde{\mathsf{h}}}

\newcommand{\Mbar}{\overline{\mathcal{M}}}
\newcommand{\Spec}{\mathrm{Spec}}
\newcommand{\sqp}{\diamond}
\newcommand{\bqp}{\star}

\newtheorem{theorem}{Theorem}[section]
\newtheorem{proposition}[theorem]{Proposition}
\newtheorem{lemma}[theorem]{Lemma}
\newtheorem{corollary}[theorem]{Corollary}
\newtheorem{conjecture}[theorem]{Conjecture}
\newtheorem{theorem/definition}[theorem]{Theorem/Definition}

\theoremstyle{definition}
 \newtheorem{example}[theorem]{Example}
\newtheorem{definition}[theorem]{Definition}

\newtheorem{problem}[theorem]{Problem}
\newtheorem{question}[theorem]{Question}
\newtheorem{algorithm}[theorem]{Algorithm}
\newtheorem{remark}[theorem]{Remark}
\newtheorem{notation}[theorem]{Notation}

\begin{document}

\title
{\large{\textbf{BIG QUANTUM COHOMOLOGY OF EVEN DIMENSIONAL INTERSECTIONS OF TWO QUADRICS}}}
\author{\normalsize Xiaowen Hu\\
\date{}
\\
}
\maketitle

\begin{abstract}
For  even dimensional smooth complete intersections, of dimension at least 4, of two quadric hypersurfaces in a projective space, we study the genus zero Gromov-Witten invariants by the monodromy group of its whole family. We compute the invariants of length 4 and show that, besides a special invariant, all genus zero Gromov-Witten invariants can be reconstructed from the invariants of length 4.
In dimension 4, we compute the special invariant by solving a curve counting problem.  We show that the generating function of genus zero Gromov-Witten invariants has a positive radius of  convergence. We show that, although the small quantum cohomology is not semisimple, the associated Frobenius manifold is generically tame semisimple.
\end{abstract}

\tableofcontents

\section{Introduction}
Quantum cohomology, or more general, higher genus Gromov-Witten invariants are defined for symplectic and algebraic manifolds in e.g. \cite{Ruan96}, \cite{Beh97}, \cite{LT98}. As a kind of basic varieties, quantum cohomology of the smooth complete intersections in projective spaces has been studied intensively. A kind of (numerical) mirror symmetry for Fano and Calabi-Yau complete intersections has been established (\cite{Giv96}), which expresses the small J-function in terms of hypergeometric series. From small J-functions one can reconstruct all genus zero Gromov-Witten invariants with only \emph{ambient cohomology classes} as insertions, i.e. the classes obtained by intersecting linear spaces with the complete intersections.  The invariants with \emph{primitive cohomology classes} as insertions are not well-understood as those with ambient insertions; for 3 point invariants see \cite{Bea95}. But they are necessary for understanding, for example, Dubrovin's conjecture on the relation between quantum cohomology and derived categories.   

From now on, we call genus zero primary (i.e. without $\psi$-classes) Gromov-Witten invariants \emph{correlators}, for brevity. 
In \cite{Hu15} we studied the big quantum cohomology of Fano complete intersections, especially the correlators with primitive insertions. One of our tools is the monodromy group of the whole family of complete intersections of a given multidegree. The Zariski closure $\overline{G}$ of such monodromy groups are either orthogonal groups, symplectic groups, or finite groups. If $\overline{G}$ is  a finite group, we call the complete intersections \emph{exceptional}, and otherwise \emph{non-exceptional}. There are only three classes of exceptional complete intersections: the cubic surfaces, the quadric hypersurfaces, and the even dimensional complete intersections of two quadric hypersurfaces (which we call  even $(2,2)$-complete intersections for brevity). The first two classes have been well-studied, in the sense that all correlators can be effectively computed from the WDVV equation, and the corresponding Frobenius manifold is semisimple at the origin. The study of the remaining exceptional complete intersections, i.e. the even $(2,2)$-complete intersections, is the topic of this paper.

To explain our approach and its limitations,  let us recall more details in \cite{Hu15}. The basic idea is still to use WDVV equations
\begin{equation}\label{eq-WDVV-intro}
\sum_e\sum_f	(\partial_{t^a}\partial_{t^b}\partial_{t^e}F)g^{ef}(\partial_{t^f}\partial_{t^c}\partial_{t^d}F)=\sum_e\sum_f(\partial_{t^a}\partial_{t^c}\partial_{t^e}F)g^{ef}(\partial_{t^f}\partial_{t^b}\partial_{t^d}F).
\end{equation}
 The most direct way to use WDVV to get recursions is to use the leading terms. Namely, by selecting a monomial $t^I$, where $I$ is a multi-index, and extracting the coefficients of $t^I$ on both sides, we get an equation of the form
 \begin{eqnarray*}
 	&&\mathrm{Coeff}_{t^I}(\partial_{t^a}\partial_{t^b}\partial_{t^e}F)g^{ef}(\partial_{t^f}\partial_{t^c}\partial_{t^d}F)(0)
 	+	(\partial_{t^a}\partial_{t^b}\partial_{t^e}F)(0)g^{ef}\mathrm{Coeff}_{t^I}(\partial_{t^f}\partial_{t^c}\partial_{t^d}F)\nn\\
 	&&-\mathrm{Coeff}_{t^I}(\partial_{t^a}\partial_{t^c}\partial_{t^e}F)g^{ef}(\partial_{t^f}\partial_{t^b}\partial_{t^d}F)(0)
 	-(\partial_{t^a}\partial_{t^c}\partial_{t^e}F)(0)g^{ef}\mathrm{Coeff}_{t^I}(\partial_{t^f}\partial_{t^b}\partial_{t^d}F)\nn\\
 	&=& \mbox{combinations of coefficients of lower order terms}.
 \end{eqnarray*}
 Here we have adopted Einstein's summation convention, i.e. omitting the summation notations of the repeated indices $e$ and $f$. We call the resulted recursions \emph{essentially linear recursions}. Such recursions are effective only when the correlators of length 3 have good properties. For non-exceptional complete intersections and the even $(2,2)$-complete intersections, the correlators with primitive insertions cannot provide recursion to compute  correlators of arbitrary lengths, essentially due to the monodromy reason. For non-exceptional ones, we remedy this by considering the symmetric reduction of the WDVV equations. From the classical theory of polynomial invariants of orthogonal groups or symplectic groups, the dependence of the generating function $F$ on the variables dual to the primitive classes are encoded in only one variable $s$. When the dimension $n$ of the complete intersection $X$ is even, $s$ is defined as
 \begin{equation*}
 	s=\frac{1}{2}\sum (t^i)^2,
 \end{equation*}
 where $t^i$ are dual variables of an orthonormal basis of $H^n_{\mathrm{prim}}(X)$. The main novelty of \cite{Hu15} is that we found that the coefficient of $s^k$ in $F$ can be recursively computed by solving quadratic equations of one variable with double roots. For a precise conjectural algorithm and partial results on this we refer the reader to \cite[Section 1]{Hu15}. There is only one case that is excluded in this algorithm: the cubic hypersurfaces. For cubic hypersurfaces we showed in \cite{Hu15} that the generating function $F$ can be computed by an essentially linear recursion with respect to the variable $s$, using the symmetric-reduced WDVV equations. Note that $s$ is quadratic in $t^i$'s. So this means that we can compute $F$ by recursion on the \emph{sub-leading terms} of (\ref{eq-WDVV-intro}) and the monodromy symmetries. As (even dimensional) cubic hypersurfaces and the even $(2,2)$-complete intersections are the only even dimensional complete intersections of Fano index $n-1$ (recall $n=\dim X$), it is reasonable to expect an analogy between them.

The monodromy group of the family of even $n$-dimensional $(2,2)$-complete intersections, or more precisely, it image in $\mathrm{Aut}(H^n_{\mathrm{prim}}(X))$, is the Weyl group $D_{n+3}$. Moreover the lattice $H^n_{\mathrm{prim}}(X)\cap H^n(X;\mathbb{Z})$ with the Poincaré pairing is isomorphic to a $D_{n+3}$-lattice with $(-1)^{\frac{n}{2}}$ times its standard inner product. We can choose a basis $\alpha_1,\dots,\alpha_{n+3}$ of the lattice $H^n_{\mathrm{prim}}(X)\cap H^n(X;\mathbb{Z})$ which maps to the $D_{n+3}$ roots under such an isomorphism, and then define a specific orthonormal basis $\epsilon_1,\dots,\epsilon_{n+3}$ of $H^n_{\mathrm{prim}}(X)$ (see (\ref{eq-roots-D}) and (\ref{eq-normalizedOrthonormalBasis})). Let $\sfh_i$ be the $i$-th cup product of the hyperplane class $\sfh$. Let $t^0,\dots,t^{n},t^{n+1},\dots,t^{2n+3}$ be the basis dual to $1,\sfh,\dots,\sfh_n,\epsilon_1,\dots,\epsilon_{n+3}$, then the generating function of correlators of $X$ is a function of $t^0,\dots,t^n$ and $s_1,\dots,s_{n+3}$, where
\begin{equation}\label{eq-invariantsOf-typeD-intro}
	s_{i}=\begin{cases}
	\vspace{0.2cm}
	\frac{1}{(2i)!}\sum_{j=n+1}^{2n+3}(t^j)^{2i},& \mbox{for}\ 1\leq i\leq n+2,\\
	\prod_{j=n+1}^{2n+3}t^j,& i=n+3.
	\end{cases}
\end{equation}
Instead of doing the  reduction of WDVV with the $D_{n+3}$ symmetry, we are going to find recursions based on the sub-leading terms of the WDVV equations, as mentioned above. For this we need first compute all the correlators of length 4. The correlators of length 4 with at most 2 primitive insertions are computed in \cite{Hu15}. For the correlators of length 4 with 4 primitive insertions, we have a uniform partial result in \cite[Section 9]{Hu15} which holds for all complete intersections of Fano index $n-1$, which is obtained by an application of Zinger's reduced genus 1 Gromov-Witten invariants \cite{Zin08}. Combining this result with the above monodromy reason and some integrality, we obtain:

\begin{theorem}\label{thm-4points-fanoIndex-even(2,2)-intro}(= Theorem \ref{thm-4points-fanoIndex-even(2,2)})
Let   $X$ be an even dimensional complete intersection of two quadrics in $\mathbb{P}^{n+2}$,  with $n\geq 4$. Then
\begin{equation}\label{eq-4points-fanoIndex-even(2,2)-intro}
	\frac{\partial^2 F}{(\partial s_1)^2}(0)=1,\ \frac{\partial F}{\partial s_2}(0)=-2.
 \end{equation}
Equivalently, for $1\leq a,b\leq n+3$,
\begin{equation}\label{eq-4points-fanoIndex-even(2,2)-ab-intro}
	\langle \epsilon_a,\epsilon_a,\epsilon_b,\epsilon_b\rangle_{0,1}=1.
\end{equation}
\end{theorem}
Using Theorem \ref{thm-4points-fanoIndex-even(2,2)-intro}, a careful study of  the sub-leading terms of the WDVV equations leads to the following.
\begin{theorem}[Reconstruction]\label{thm-reconstruction-even(2,2)-intro}(= Theorem \ref{thm-reconstruction-even(2,2)})
Let   $X$ be an even dimensional complete intersection of two quadrics in $\mathbb{P}^{n+2}$,  with $n\geq 4$.
With the knowledge of the 4-point invariants, all the invariants can be reconstructed from the WDVV, the deformation invariance, and the \emph{special correlator}
\begin{equation}\label{eq-intro-length(n+3)Invariant-even(2,2)}
	\langle \epsilon_{1},\dots,\epsilon_{n+3}\rangle_{0,n+3,\frac{n}{2}}.
\end{equation}
\end{theorem}

We have made more intricate attempts on the WDVV equation, e.g. some equations that  a priori may give quadratic equations of the special correlator (\ref{eq-intro-length(n+3)Invariant-even(2,2)}). But  such equations in examples turn out to be trivial. This gives reason to suspect that WDVV does not give any new information on the special correlator; see Conjecture \ref{conj-specialCorrelator-free} for a precise statement. We can show a weaker result: 
\begin{lemma}\label{lem-specialCorrelator-freedomOfSign}
For even $(2,2)$-complete intersections, the WDVV equations and the knowledge of correlators of length 4 can at most determine the special correlator with freedom of signs, unless it vanishes.
\end{lemma}
The proof is immediate  : we change the signs of the basis $\alpha_1,\dots,\alpha_{n+3}$. The resulted basis is still allowable as above. Then each correlator will change by a sign $(-1)^k$, where $k$ is the number of primitive insertions. 

We refer the reader to Remark \ref{rem:choiceOfBasis}, and Section \ref{sec:conjecturesOnSpecialCorrelator}, on the choice of the basis. In Section \ref{sec:explictD-Lattice} we give an explicit construction of the basis, which we take as our standard choice.

Nevertheless, Theorem \ref{thm-4points-fanoIndex-even(2,2)-intro} and \ref{thm-reconstruction-even(2,2)-intro} are enough to imply the following.

\begin{theorem}[Analyticity and Semisimplicity]\label{thm-intro-convergence-semisimplicity}
There exists an open (in the classical topology) neighborhood of the origin of $\mathbb{C}^{2n+4}$, on which the generating function $F(t^0,\dots,t^{2n+3})$ is analytic and defines a (generically) tame semisimple Frobenius manifold.
\end{theorem}
The analyticity follows from a bound of the correlators, which follows from an induction based on an algorithm given by the proof of Theorem \ref{thm-reconstruction-even(2,2)-intro}. Let $\widetilde{E}$ be the matrix of the big quantum multiplication by the Euler field $E$.   We show the semisimplicity by showing that 
\begin{equation}\label{eq-statement-cutoffEulerField-distinctEigenValues}
 	\mbox{the cutoff of $\widetilde{E}$ at order 2 has pairwise distinct eigenvalues.}
\end{equation}
I would like to provide a prophetic view of why this is possible (see also \cite[Remark 3.2]{Hu15}). If the generating function $F$ has continuous symmetries such as the orthogonal or symplectic groups in the cases of non-exceptional complete intersections, one can construct a family of (normalized)  Euler  fields for the associated Frobenius manifold $\mathcal{M}$. But  semisimple Frobenius manifold has a unique normalized Euler field \cite[Theorem I.3.6]{Man99}. So such $\mathcal{M}$ cannot be generically semisimple. For the same reason, if one wants to show that the Frobenius manifold associated with the generating function $F$ of quantum cohomology of an even $(2,2)$-complete intersection $X$ is semisimple, one needs to use the expansion of $F$ to a certain order  such that it has no continuous symmetry. This explains why the small quantum cohomology of $X$ is not semisimple: the cutoff of $F$ at order 3 is a function of $t^0,\dots,t^n$ and 
\[
s_1=\sum_{i=n+1}^{2n+3}(t^i)^2.
\]
So it has symmetries from the orthogonal group $\mathrm{O}\big(H^n_{\mathrm{prim}}(X)\big)$. The degree 4 form 
\[
s_2=\sum_{i=n+1}^{2n+3}(t^i)^4
\]
has only finitely many automorphisms; this is a classical theorem of Jordan \cite{Jor1880}. So we can expect that the information of  correlators of length 4 implies the semisimplicity. To capture the full information of  correlators of length 4  we need only the cutoff of the matrix $\widetilde{E}$  at an order $\geq 2$.

In view of Kuznetsov's results \cite[Corollary 5.7]{Kuz08}, the semisimplicity statement above confirms a part of Dubrovin's conjecture \cite{Dub98} which predicts the equivalence of the existence of full exceptional collections in $D^b(\mathrm{Coh}(X))$ and the generic semisimplicity of the Frobenius manifold of the quantum cohomology of $X$.\\

By an analogy to a vanishing conjecture \cite[Conjecture 10.26]{Hu15} of the coefficient of $s^{n+1}$ in $F$ for cubic hypersurfaces,  we make the following conjecture on the special correlator. For details see Section \ref{sec:conjecturesOnSpecialCorrelator}. 
\begin{conjecture}\label{conj-unknownCorrelator-Even(2,2)-intro}
For an even $n$ dimensional complete intersection of two quadrics in $\mathbb{P}^{n+2}$, let $\varepsilon_1,\dots,\varepsilon_{n+3}$ be the basis of $H^n_{\mathrm{prim}}(X)$ defined in Section \ref{sec:explictD-Lattice}. Then 
\begin{equation}\label{eq-unknownCorrelator-Even(2,2)-intro}
	\langle \varepsilon_{1},\dots,\varepsilon_{n+3}\rangle_{0,n+3,\frac{n}{2}}=\frac{(-1)^{\frac{n}{2}}}{2}.
\end{equation}
\end{conjecture}

I must confess that this analogy is not so strong to make us confident about the validity. We can only show Conjecture \ref{conj-unknownCorrelator-Even(2,2)-intro} in dimension 4.
\begin{theorem}\label{thm-unknownCorrelator-Even(2,2)-4dim-intro}(= Corollary \ref{cor-unknownCorrelator-Even(2,2)-4dim})
Let $\varepsilon_1,\dots,\varepsilon_{n+3}$ be the basis of $H^{n}_{\mathrm{prim}}(X)$ defined as in Section \ref{sec:explictD-Lattice}, and  $\epsilon_1,\dots,\epsilon_{n+3}$ be defined as in (\ref{eq-normalizedOrthonormalBasis}). 
Then for any  $4$-dimensional complete intersections of two quadrics in $\mathbb{P}^{6}$, 
\begin{equation}\label{eq-unknownCorrelator-Even(2,2)-4dim-intro}
	\langle \epsilon_{1},\dots,\epsilon_{7}\rangle_{0,7,2}=\frac{1}{2}.
\end{equation}
\end{theorem}
The proof is essentially a computation \emph{from the first principle}. Namely, we find explicit cycles $S_{0,1},S_{1,2},\dots,S_{5,6},S_{6,0}$ in $X$ such that  (\ref{eq-unknownCorrelator-Even(2,2)-4dim-intro}) is equivalent to
\begin{equation}\label{eq-enumerativeCorrelator-Dim4-intro}
	 \langle [S_{0,1}],[S_{1,2}],[S_{2,3}],[S_{3,4}],[S_{4,5}],[S_{5,6}],[S_{6,0}]\rangle_{0,7,2}=1.
\end{equation}
Then we show (\ref{eq-enumerativeCorrelator-Dim4-intro}) by the definition of Gromov-Witten invariants: we count the curves that intersect with the given cycles, and use deformation theory to  compute the intersection multiplicity of the image of the virtual fundamental cycle and cycles in $X^7$ (see Definition \ref{def-enumerativeCorrelator}). As a byproduct we obtain an enumerative result with classical flavor.
\begin{theorem}\label{thm-countingConics-4dim-general-intro}(= Theorem \ref{thm-countingConics-4dim-general})
For general 4-dimensional  smooth complete intersections $X$ of two quadrics in $\mathbb{P}^6$, there exists exactly one smooth conic that meets each of the 2-planes $S_{i,i+1}$ in $X$ for $0\leq i\leq 6$.
\end{theorem}

Let us emphasize the significance of the computation of the special correlator. In view of our previous results and conjectures (see \cite[Conjecture 1.15 and Table 1]{Hu15}), in the algorithmic sense the special correlator  is the only genus 0  Gromov-Witten invariant of complete intersections, at least the even dimensional ones,  that we  have no definite way to compute. But see Remark \ref{rem:compute-specialCorrelator-fromHigherGenusInv} for a possible approach via higher genus Gromov-Witten invariants with ambient insertions, as a byproduct of Theorem \ref{thm-intro-convergence-semisimplicity}. Shortly after the first version of this paper, Argüz, Bousseau, Pandharipande, and Zvonkine  (\cite{ABPZ21}) showed that one can compute all genera Gromov-Witten invariants of the even dimensional Gromov-Witten invariants by the degeneration formula. However, the computation might be rather complicated to carry out.\\

We are now in a position to make an additional remark on the semisimplicity.
There seems to be a folklore conjecture that if a smooth projective variety has generically semisimple quantum cohomology, then there exist finitely many correlators that can determine all correlators by WDVV and the Euler field. One might be tempted to make a stronger and, quantitative, conjecture: if the values of the correlators of length $\leq k$ suffice to imply the generic semisimplicity, then they determine all  correlators. Our Theorem \ref{thm-reconstruction-even(2,2)-intro} confirms the first conjecture in the case of the even $(2,2)$-complete intersections. On the other hand, by the 
claim (\ref{eq-statement-cutoffEulerField-distinctEigenValues}), Lemma \ref{lem-specialCorrelator-freedomOfSign}, and Theorem \ref{thm-unknownCorrelator-Even(2,2)-4dim-intro}, the 4-dimensional $(2,2)$-complete intersections are  counterexamples of the latter one. From Example \ref{example-f(6)}-\ref{example-f(10)} and integrality, one sees that the special correlator of the  6, 8, and 10 dimensional $(2,2)$-complete intersections also do not vanish, so they are also counterexamples. We summarize our observation as the following assertion.

\begin{corollary}
There exist  a  formal Frobenius manifold $\mathcal{M}$ with a potential  $\Phi$ and a natural number $k$, such that the expansion $\check{\Phi}$ of $\Phi$ at $\mathbf{t}=0$ (where $\mathbf{t}$ is a system of flat coordinates) at order $k$ yields that $\mathcal{M}$ is generically semisimple, but $\mathcal{M}$ is not uniquely determined by $\check{\Phi}$ from the WDVV equation and the Euler field.
\end{corollary}


The paper is organized as follows. In Section \ref{sec:preliminaries} we recall the properties of Gromov-Witten invariants we need, and the primitive cohomology of an even $(2,2)$-complete intersection $X$ as a representation of type-$D$ Weyl group. In Section \ref{sec:correlators-lengt-atMost4} we compute the genus zero GW invariants of $X$ of length at most 4. Then in Section \ref{sec:reconstructionTheorem} we show the reconstruction theorem from the data in  Section \ref{sec:correlators-lengt-atMost4}, and propose the conjecture on the special correlator. In Section \ref{sec:convergence} we show the convergence of the generating function in an open neighborhood of the origin, which makes the study of the semisimplicity easier, than only a formal Frobenius manifold. In Section \ref{sec:semisimplicity} we show that the Euler field has distinct eigenvalues at a general point of the Frobenius manifold. In Section \ref{sec:EnumerativeGeometry-Even(2,2)} we relate Conjecture \ref{conj-unknownCorrelator-Even(2,2)-intro} to a curve counting problem on $X$ and solve it, and thus show Conjecture \ref{conj-unknownCorrelator-Even(2,2)-intro}, in dimension 4.

In the Appendix, we present an algorithm based on the proof of Theorem \ref{thm-reconstruction-even(2,2)-intro} to compute the primary genus 0 Gromov-Witten invariants of an even $(2,2)$-complete intersection of dimension $\geq 4$, with the special correlator (\ref{eq-intro-length(n+3)Invariant-even(2,2)}) as an unknown. The algorithm is implemented in our Macaulay2 package 
\texttt{QuantumCohomologyFanoCompleteIntersection}.
We write also a Macaulay2 package
 \texttt{ConicsOn4DimIntersectionOfTwoQuadrics} 
 for the computations in Section \ref{sec:EnumerativeGeometry-Even(2,2)}. The reader can find the packages at

\url{https://github.com/huxw06/Quantum-cohomology-of-Fano-complete-intersections}

\vspace{0.2cm}

\emph{Acknowledgement}\quad 
 I am grateful to  Hua-Zhong Ke for enlightening discussions. I also thank  Huai-Liang Chang, Weiqiang He, Christopher Lyons, Giosuè Muratore, Maxim Smirnov, and  Jinxing Xu  for  discussions on various related topics. 
  This work is supported by  NSFC 11701579.

\section{Preliminaries}\label{sec:preliminaries}
\subsection{Properties of Gromov-Witten invariants}
We recall the definition of the Gromov-Witten invariants and their properties that we need to use in this paper. Our main reference is
\cite[Chapter VI]{Man99}. All schemes in this paper are over $\mathbb{C}$.

Let $X$ be a smooth projective scheme over $\mathbb{C}$ of dimension $n$. Let $k\in \mathbb{Z}_{\geq 0}$, and $\beta\in H_2(X;\mathbb{Z})/\mathrm{tor}$. The stack $\Mbar_{g,k}(X,\beta)$ of stable maps of degree $\beta$ from $k$-point genus $g$ marked semistable curves to $X$  is a proper Deligne-Mumford stack and carries a virtual fundamental class (\cite{BF97}, \cite{LT98}) $[\Mbar_{g,k}(X,\beta)]^{\mathrm{vir}}$ of dimension $(1-g)(n-3)+k+c_1(T_X)\cdot \beta$. For each $1\leq i\leq k$, the section $\sigma_i$ pulls back the relative cotangent line bundle of the universal curve to form a line bundle on $\Mbar_{g,k}(X,\beta)$, whose first Chern class is denoted by $\psi_i$; moreover there is an associated \emph{evaluation map} $\mathrm{ev}_i=f\circ \sigma_i$, where $f$ is the universal stable map. For $\gamma_1,\dots,\gamma_k\in H^*(X;\mathbb{Q})$ and $a_1,\dots,a_k\in \mathbb{Z}_{\geq 0}$, there is an associated \emph{Gromov-Witten invariant}
\begin{equation*}
	\langle \psi_1^{a_1}\gamma_1,\dots,\psi_k^{a_k}\gamma_k\rangle_{g,k,\beta}^X:=\int_{[\Mbar_{g,k}(X,\beta)]^{\mathrm{vir}}}\prod_{i=1}^{k}\psi_i^{a_i}\mathrm{ev}_i^*\gamma_i\in \mathbb{Q}.
\end{equation*}
A term like $\psi_i^{a_i}\gamma_i$ in $\langle \psi_1^{a_1}\gamma_1,\dots,\psi_k^{a_k}\gamma_k\rangle_{g,k,\beta}$ is called an \emph{insertion} of this invariant. We say $(g,k,\beta)$ is in the \emph{stable range} if either $2g-2+k>0$ or $\beta$ is a nonzero effective curve class.
It is convenient to use simplified notations on the following occasions:
\begin{enumerate}
	\item[(i)] The superscript $X$ will be omitted when it is obvious;
	\item[(iii)] the subscript $k$  might be dropped when $k$ is obvious from the expression;
	\item[(iv)] the subscript  $\beta$ might be dropped when it can be uniquely determined by the insertions and the following condition (\ref{eq-Dim}), which is always the case for Fano complete intersections in projective spaces.
\end{enumerate}

The GW invariants $\langle \psi_1^{a_1}\gamma_1,\dots,\psi_k^{a_k}\gamma_k\rangle_{g,k,\beta}$ with $a_1=\dots=a_k=0$ are called \emph{primary}. In this paper we only deal with primary genus 0 invariants, although some results depend on the computation of certain non-primary invariants and genus 1 invariants in \cite{Hu15}. For brevity we call the  genus 0 primary Gromov-Witten invariants \emph{correlators}.

For simplicity we assume $H^{\mathrm{odd}}(X;\mathbb{Q})=0$, which is the case for even dimensional $(2,2)$-complete intersections. We denote $H^*(X)=H^*(X;\mathbb{C})$. 
For  two cohomology classes $\gamma_1$ and $\gamma_2$, we denote the Poincaré pairing by
\begin{equation*}
	(\gamma_1,\gamma_2):=\int_X \gamma_1\cup \gamma_2.
\end{equation*}
If a basis $\gamma_0,\dots,\gamma_N$ of $H^*(X)$ is given, we denote $g_{a,b}=(\gamma_a,\gamma_b)$ for $0\leq a,b\leq N$, and set $(g^{a,b})_{0\leq a,b\leq N}$ to be the inverse matrix of $(g_{a,b})_{0\leq a,b\leq N}$.
Let $T^0,\dots,T^N$ be dual basis with respect to $\gamma_0,\dots,\gamma^N$, then the genus $g$ generating function is defined as
 \begin{equation}\label{eq-def-generatingFunction}
 	\mathcal{F}_g(T^0,\dots,T^{N})=\sum_{k\geq 0} \sum_{\beta} \frac{1}{k!}\big\langle \sum_{i=0}^N \gamma_i T^i,\dots,\sum_{i=0}^N \gamma_i T^i\big\rangle_{g,k,\beta},
 \end{equation}
where the invariants outside of the stable range are defined to be zero, by convention. Here we have implicitly used that $X$ is Fano, so that with fixed insertions there are only finitely many $\beta$ such that the invariant does not vanish. We denote $F= \mathcal{F}_0$.  

We record some standard properties of GW invariants as follows.
\begin{enumerate}
	\item[(i)] Degree 0 correlators:
		\begin{gather}\label{eq-Deg0}\tag{Deg0}
		\langle \gamma_1,\dots,\gamma_k\rangle_{g,k,0}=
		\begin{cases}
		\int_{X}\gamma_1\cup \gamma_2\cup \gamma_3,& \mbox{if}\ g=0, k=3;\nn\\
		-\frac{1}{24}\int_X \gamma_1\cup c_{n-1}(T_X), & \mbox{if}\ g=1, k=1,\\
		0, & \mbox{if}\ 2g-2+k\geq 2.
		\end{cases}
		\end{gather}
	\item[(ii)]	Suppose each $\gamma_i$ has pure real degree $|\gamma_i|$, then there is the dimension constraint:
		\begin{equation}\label{eq-Dim}\tag{Dim}
			\langle \gamma_{1},\dots, \gamma_{k}\rangle_{g,k,\beta}=0\ \mbox{unless}\
			\sum_{i=1}^k \frac{1}{2}|\gamma_{b_k}|=(1-g)(n-3)+k+c_1(T_X)\cap \beta.
		\end{equation}
	\item[(iii)] The $S_n$-equivariance:
		\begin{equation}\label{eq-Sym}\tag{Sym}
			\langle \gamma_{1},\dots, \gamma_{i-1}, \gamma_{i},\dots, \gamma_{b_k}\rangle_{g,k,\beta}
			=\langle \gamma_{1},\dots, \gamma_{i}, \gamma_{i-1},\dots,\gamma_{b_k}\rangle_{g,k,\beta}.
		\end{equation}
	\item[(iv)]
			The divisor equation: for $\gamma\in H^2(X)$,
			\begin{gather}\label{eq-Div}\tag{Div}
				\langle \gamma_1,\dots,\gamma_k,\gamma\rangle_{g,k+1,\beta}=(\gamma\cap \beta)\langle \gamma_1,\dots,\gamma_k\rangle_{g,k,\beta}.
			\end{gather}
	\item[(v)] Fundamental class axiom:
			\begin{gather}\label{eq-FCA}\tag{FCA}
				\langle 1, \gamma_{1},\dots,\gamma_{k-1}\rangle_{g,k,\beta}=
				\begin{cases}
				(\gamma_1,\gamma_2),& \mbox{if}\ g=0, k=3, \beta=0;\\
				0, & \mbox{if}\ 3g-3+k\geq 1\ \mbox{or}\ \beta\neq 0.
				\end{cases}
			\end{gather}
	\item[(vi)] If a basis $\gamma_0,\dots,\gamma_N$ of $H^*(X)$ and $T^0,\dots,T^N$ is the dual basis, then we have the WDVV equations, for $0\leq a,b\leq N$:
		\begin{gather}\label{eq-WDVV}\tag{WDVV}
		\sum_{e=0}^N \sum_{f=0}^N \frac{\partial^3 F}{\partial T^a \partial T^b\partial T^e}g^{ef}\frac{\partial^3 F}{\partial T^f \partial T^c\partial T^d}
		=\sum_{e=0}^N \sum_{f=0}^N \frac{\partial^3 F}{\partial T^a \partial T^c\partial T^e}g^{ef}\frac{\partial^3 F}{\partial T^f \partial T^b\partial T^d}.
		\end{gather}

\end{enumerate}

Now suppose  $\gamma_0,\dots,\gamma_N$ is a basis of $H^*(X)$ such that each $\gamma_i$ has a pure degree,  and let $T^0,\dots,T^N$ be the dual basis have pure degrees. Let
\[
c_1(T_X)=\sum_{i=0}^N a_i \gamma_i.
\]
Of course $a_i=0$ unless $|\gamma_i|=2$. The Euler vector field is defined as
\begin{equation}\label{eq-EV-0}
	E=\sum_{i=0}^N(1-\frac{|\gamma_i|}{2})\frac{\partial }{\partial T^i}+\sum_{i=0}^N a_i \frac{\partial }{\partial T^i}.
\end{equation}
Then (\ref{eq-Dim}) and the divisor equation (\ref{eq-Div}) imply 
\begin{gather}\label{eq-EulerVectorField}\tag{EV}
	EF=(3-n)F+  \sum_{i=0}^N a_i \frac{\partial }{\partial T^i} c,
\end{gather}
where  $c$ is the classical cubic intersection form 
\begin{eqnarray}\label{eq-cubicIntersectionForm}
c(t_0, \cdots, t^{n+m})=\sum_a\sum_b\sum_c\frac{t^{a}t^{b}t^{c}}{6}\int_{X}\gamma_a \gamma_b \gamma_c.
\end{eqnarray}

The \emph{big quantum product} is defined as
\begin{equation*}
 	\gamma_a\bqp \gamma_b=\sum_{e}\sum_f\frac{\partial^3 \mathsf{F}}{\partial T^a \partial T^b\partial T^e}g^{ef}\gamma_f,
 \end{equation*} 
 and the \emph{small quantum product} is defined as
\begin{equation*}
 	\gamma_a\sqp \gamma_b=\gamma_a\bqp \gamma_b|_{T^0=\dots=T^N=0}.
 \end{equation*}

\subsection{Monodromy group and the \texorpdfstring{$D_{n+3}$}{D{n+3}} lattice}\label{sec:monodromy-lattice}
Let $n\geq 4$ be even. 
Let $X$ be a complete intersection of two quadric hypersurfaces in $\mathbb{P}^{n+2}$. Then $H^*(X)=H^{\mathrm{even}}(X)$, and the  Fano index of $X$ is $n-1$. Denote the hyperplane class of $X$ by $\sfh$, and
\begin{equation*}
	\sfh_{i}=\underbrace{\sfh\cup\dots\cup\sfh}_{i\ \sfh's}.
\end{equation*}

Let $V=\mathbb{R}^{n+3}$ be the Euclidean space with the standard inner product. Let $\varepsilon_1,\dots,\varepsilon_{n+3}$ be an orthonormal basis, and let
\begin{equation}\label{eq-roots-D}
	\begin{cases}
	\alpha_i=\varepsilon_{i}-\varepsilon_{i+1}\ \mbox{for}\ 1\leq i\leq n+2,\\
	\alpha_{n+3}=\varepsilon_{n+2}+\varepsilon_{n+3}.
	\end{cases}
\end{equation} 
Given a vector $\alpha\in V$, the \emph{reflection with respect to $\alpha$} is the linear automorphism of $V$ defined by
\[
\gamma\mapsto \gamma-\frac{2(\gamma,\alpha)}{(\alpha,\alpha)}\alpha.
\]
The Weyl group $D_{n+3}\subset \mathrm{GL}(n+3,\mathbb{R})$ is generated the reflections with respect to the $\alpha_i$'s. If one writes vectors in $\mathbb{R}^{n+3}$ in terms of the coordinates according to the basis $\varepsilon_1,\dots,\varepsilon_{n+3}$, i.e. 
\[
\mathbf{v}=(v_1,\dots,v_{n+3})=\sum_{i=1}^{n+3}v_i \varepsilon_i,
\]
then  the group $D_{n+3}$ coincides with the group generated by the permutations of the coordinates, and the change of signs
\[
(v_1,\dots,v_{n+1},v_{n+2},v_{n+3})\mapsto (v_1,\dots,v_{n+1},-v_{n+2},-v_{n+3}).
\]
By \cite[\S 5]{Del73}, the monodromy group of the whole family of smooth complete intersections of two quadrics in $\mathbb{P}^{n+2}$ is isomorphic to $D_{n+3}$, and in this way the primitive cohomology $H^n_{\mathrm{prim}}(X)$ is a standard representation of $D_{n+3}$, and the integral lattice $H^n_{\mathrm{prim}}(X)\cap H^n(X;\mathbb{Z})$ is generated by the roots $\alpha_i$'s of $D_{n+3}$. In other words, there is an isomorphism
\begin{equation}\label{eq-isomorphism-primCoh-even(2,2)}
	V\otimes_{\mathbb{R}}\mathbb{C}\xrightarrow{\sim} H^n_{\mathrm{prim}}(X),
\end{equation}
via which lattice in $V$ generated by $\alpha_1,\dots,\alpha_{n+3}$ is mapped onto $H^n_{\mathrm{prim}}(X)\cap H^n(X;\mathbb{Z})$.
\begin{proposition}\label{prop-lattice-sign}
 We equip  $H^n_{\mathrm{prim}}(X)$ with the inner product $(-1)^{\frac{n}{2}}(.,.)$, where  $(.,.)$ is the Poincaré pairing. Then (\ref{eq-isomorphism-primCoh-even(2,2)}) becomes an isometry. 
\end{proposition}
\begin{proof}
The Poincaré pairing is invariant under monodromies. As we will recall in the following, the degree 2 polynomial invariant of the standard representation of $D_{n+3}$ is generated by the Poincaré pairing. By the Hodge index theorem, $(-1)^{\frac{n}{2}}(.,.)$ is positive definite. So the inner product on $H^n_{\mathrm{prim}}(X)$ induced from $V$ by (\ref{eq-isomorphism-primCoh-even(2,2)}) coincides with $(-1)^{\frac{n}{2}}(.,.)$ up to a positive constant multiple. Since the discriminant of the lattice $\alpha_1,\dots,\alpha_{n+3}$ equals 4, we are left to show that the discriminant of the lattice $H^n_{\mathrm{prim}}(X)\cap H^n(X;\mathbb{Z})$ is $\pm 4$. First we note that there exists $\frac{n}{2}$-planes in $X$ (\cite[Corollary 3.3]{Rei72}), whose intersection number with $\sfh_{n/2}$ is 1, from which it follows that  the sub-lattice $L=H^*_{\mathrm{amb}}(X)\cap H^*(X;\mathbb{Z})$ is generated by $\sfh_{n/2}$. 
Since $H^n(X;\mathbb{Z})$ is unimodular and   $(\sfh_{n/2},\sfh_{n/2})=4$, by \cite[Prop. 5.3.3]{Kit93} the discriminant of $L^{\perp}=H^n_{\mathrm{prim}}(X)\cap H^n(X;\mathbb{Z})$ is $\pm 4$. So we are done.
\end{proof}

\begin{corollary}
$H^{n}(X;\mathbb{Z})$ is generated by the  classes of $\frac{n}{2}$-planes in $X$.
\end{corollary}
\begin{proof}
By \cite[Theorem 3.14]{Rei72}, the discriminant of the lattice of the classes of $\frac{n}{2}$-planes in $X$ is $\pm 4$, which coincides with that of the lattice $H^{n}(X;\mathbb{Z})$ (see the proof of Proposition \ref{prop-lattice-sign}).
\end{proof}

 Via (\ref{eq-isomorphism-primCoh-even(2,2)}), we identify the vectors $\varepsilon_i$ and $\alpha_i$ with their images in $H^n_{\mathrm{prim}}(X)$. Moreover we define, for $1\leq i\leq n+3$,
\begin{equation}\label{eq-normalizedOrthonormalBasis}
	\epsilon_i=\begin{cases}
	\varepsilon_i, & \mbox{if}\ n\equiv 0 \mod 4;\\
	\sqrt{-1}\varepsilon_i, & \mbox{if}\ n\equiv 2 \mod 4.
	\end{cases}
\end{equation}
Then $\epsilon_1,\dots,\epsilon_{n+3}$ is an \emph{orthonormal basis} of $H^n_{\mathrm{prim}}(X)$.

 Let $t^{n+1},\dots,t^{2n+3}$ be the basis of $H^*_{\mathrm{prim}}(X)^{\vee}$ dual to $\epsilon_1,\dots,\epsilon_{n+3}$.  By the invariant theory of Weyl groups \cite[\S 3.12]{Hum90}, 
the polynomial invariants of $D_{n+3}$ are generated by $s_{1},\dots,s_{n+3}$, where
\begin{equation}\label{eq-invariantsOf-typeD-1}
	s_{i}=\frac{1}{(2i)!}\sum_{j=n+1}^{2n+3}(t^j)^{2i},\ \mbox{for}\ 1\leq i\leq n+2,
\end{equation}
and
\begin{equation}\label{eq-invariantsOf-typeD-2}
	s_{n+3}=\prod_{j=n+1}^{2n+3}t^j.
\end{equation}
Moreover, $s_1,\dots,s_{n+3}$ are algebraically independent. As a consequence of the deformation invariance of Gromov-Witten invariants, we have:
\begin{theorem}\label{thm-monodromy-evenDim(2,2)}
The genus $g$ generating function $\mathcal{F}_g$  of $X$  can be written in a unique way as a series of  $s_1,\dots,s_{n+3}$.
\end{theorem}
One can see this by directly quoting the definition of Gromov-Witten invariants via symplectic geometry. For an algebraic proof using \cite[Theorem 4.2]{LT98}, see \cite[Corollary 3.2]{Hu15}.

We introduce some notations that will be used throughout this paper. We denote the genus 0 generating function by $F$. The basis dual to $1,\sfh,\dots,\sfh_{n},\epsilon_1,\dots,\epsilon_{n+3}$ is denoted by $t^0,\dots,t^{2n+3}$. Denote the small quantum multiplication by $\diamond$. Let 
\begin{equation*}
	\tsfh_{i}=	\underbrace{\sfh\diamond\dots\diamond\sfh}_{i\ \sfh's}.
\end{equation*}
The basis dual to $1,\tsfh,\dots,\tsfh_{n},\epsilon_1,\dots,\epsilon_{n+3}$ is denoted by $\tau^0,\dots,\tau^{n+3}$.

\begin{remark}\label{rem:choiceOfBasis}
The choice of the isomorphism from a $D_{n+3}$-lattice to $\big(H^n_{\mathrm{prim}}(X)\cap H^n(X;\mathbb{Z}),(-1)^{\frac{n}{2}}(.,.)\big)$, or equivalently, the choice of $\varepsilon_1,\dots,\varepsilon_{n+3}$ in 
$H^n_{\mathrm{prim}}(X)$ as above, is not unique. Different choices differ by automorphisms of the $D_{n+3}$-lattice.

By \cite[Theorem 1]{KM83}, the automorphism group, which we denote by $G_{n+3}$,  of a $D_{n+3}$-lattice is generated by $D_{n+3}$ and the automorphism group of the $D_{n+3}$ Dynkin diagram. Since $n>0$, the latter group is $\mathbb{Z}/2 \mathbb{Z}$. Then $G_{n+3}$ is the semidirect product $D_{n+3}\rtimes \mathbb{Z}/2 \mathbb{Z}$, and $G_{n+3}$ is generated by $D_{n+3}$ and the map which fixes $\alpha_i$ for $1\leq i\leq n+1$ and interchanges $\alpha_{n+2}$ and $\alpha_{n+3}$. Combined with the description of the $D_{n+3}$ action on $V$ recalled at the beginning of this section, we can also describe $G_{n+3}$ as generated by $D_{n+3}$ and the map $-1$ which sends $\varepsilon_i$ to $-\varepsilon_i$ for $1\leq i\leq n+3$.

As a consequence of Theorem \ref{thm-monodromy-evenDim(2,2)}, two choices of the basis $\varepsilon_1,\dots,\varepsilon_{n+3}$ which differ by an action by  $g\in D_{n+3}$ do not affect the values of the correlators and also the generating functions $\mathcal{F}_g$. But two choices  which differ by an action by the map $-1$ do affect the value of the correlator with an odd number of primitive insertions,  for example
\[
\langle \varepsilon_1,\dots,\varepsilon_{n+3}\rangle_{0,n+3,\frac{n}{2}}.
\]
In Section \ref{sec:explictD-Lattice} we will give an explicit choice of the $D_{n+3}$-roots $\alpha_i$ and the basis $\varepsilon_i$.
\end{remark}

\section{Correlators of length at most 4}\label{sec:correlators-lengt-atMost4}
In this section, we fix a smooth complete intersection $X$ of two quadrics in $\mathbb{P}^{n+2}$, where $n$ is even and $\geq 4$. 

\subsection{Recap of known results}\label{sec:knownResults-correlators}
By Theorem \ref{thm-monodromy-evenDim(2,2)}, we expand $F$ as a series of $s_1,\dots,s_{n+3}$; in particular we denote the constant term of this expansion by $F^{(0)}$, and the coefficient of $s_1$ by $F^{(1)}$. Then $F^{(0)}$ is the generating function of genus 0 primary GW invariants with only \emph{ambient} insertions, and $F^{(1)}$ is the  generating function of genus 0 primary GW invariants with exactly two primitive insertions. By \cite[Example 6.11]{Hu15},
\begin{equation}\label{eq-tauTot}
\begin{cases}
\tau^0=t^0-4t^{n-1},\\
\tau^1=t^1-12t^{n},\\
\tau^i=t^i\ \mbox{for}\ i\geq 2.
\end{cases}
\end{equation}
Then the information of the correlators of length 3 and length 4 with at most two primitive insertions are encoded in (see \cite[Lemma 6.5]{Hu15})
\begin{eqnarray}\label{eq-qp1.5}
\frac{\partial F^{(0)}}{\partial \tau^a \partial\tau^b \partial\tau^c}(0)=\left\{
\begin{array}{cc}
16^{\frac{a+b+c-n}{n-1}} \times 4, & \mathrm{if}\ \frac{a+b+c-n}{n-1}\in \mathbb{Z}_{\geq 0}; \\
0, & \mathrm{otherwise}.
\end{array}
\right.
\end{eqnarray}
and (see  \cite[Example 6.11]{Hu15})
\begin{equation}\label{eq-F122-tau}
F^{(1)}(\tau)
=\tau^0-2\sum_{\begin{subarray}{c}1\leq i,j\leq n\\
i+j=n\end{subarray}
}
\tau^i \tau^{j}
-16\tau^{n-1}\tau^{n}
+O\big((\tau)^3\big),
\end{equation}
or equivalently
\begin{eqnarray}\label{eq-F^122}
F^{(1)}
=t^0-4t^{n-1}-2\sum_{i=1}^{n-1}t^{i}t^{n-i}-16t^{n-1}t^{n}+O((t)^3).
\end{eqnarray}
The leading terms of (\ref{eq-F^122}) is contained in \cite{Bea95}, and (\ref{eq-F^122}) is obtained by WDVV and deformation invariance.
As in \cite[Section 6.1]{Hu15}, we define two transition matrices $M$ and $W$ by
\begin{eqnarray}\label{eq-transform1}
\sfh_i=\sum_{j=0}^n M_{i}^{j}\tsfh_{j},\
\tsfh_i=\sum_{j=0}^n W_{i}^{j}\sfh_{j},\ \mbox{for } 0\leq i\leq n,
\label{eq-transform1-2}
\end{eqnarray}
or equivalently
\begin{eqnarray}\label{eq-transform2}
\tau^i=\sum_{j=0}^n M_j^i t^j,\
t^i=\sum_{j=0}^n W_j^i \tau^j.
\end{eqnarray}
Then by (\ref{eq-tauTot}), 
\begin{equation}\label{eq-tTotau}
\begin{cases}
t^0=\tau^0+4 \tau^{n-1},\\
t^1=\tau^1+12 \tau^{n},\\
t^i=\tau^i\ \mbox{for}\ i\geq 2.
\end{cases}
\end{equation}
Let $\gamma_0,\dots,\gamma_{2n+3}$ be the basis $1,\sfh,\dots,\sfh_n,\epsilon_1,\dots,\epsilon_{n+3}$. For $0\leq e,f\leq 2n+3$, we define
\[
g_{ef}=(\gamma_e,\gamma_f),
\]
and define $(g^{ef})$ to be the inverse matrix of $(g_{ef})_{0\leq e,f\leq 2n+3}$. Let $\widetilde{\gamma}_0,\dots,\widetilde{\gamma}_{2n+3}$ be the basis
\[
1,\tsfh,\dots,\tsfh_n,\epsilon_1,\dots,\epsilon_{n+3}.
 \] For $0\leq e,f\leq 2n+3,$ we define
\[
\eta_{ef}=(\widetilde{\gamma}_e,\widetilde{\gamma}_f),
\]
and define $(\eta^{ef})$ to be the inverse matrix of $(\eta_{ef})_{0\leq e,f\leq 2n+3}$. Then by \cite[Lemma 6.2]{Hu15},
\begin{equation}\label{eq-etaInversePairing-even(2,2)}
	\eta^{ef}=\begin{cases}
	-4,& \mbox{if}\ e+f=1;\\
	\frac{1}{4},& \mbox{if}\ e+f=n,\\
	\delta_{e,f},& \mbox{if}\ n+1\leq e,f\leq 2n+3,\\
	0,& \mbox{otherwise}.
	\end{cases}
\end{equation}

\subsection{Some preparatory computations}\label{sec:preparatoryComputation-even(2,2)}
In this section we compute several genus 0 GW invariants of $X$ that are needed in the proof of Theorem \ref{thm-4points-fanoIndex-even(2,2)}. We follow the notations in Section \ref{sec:monodromy-lattice}. For convenience in summations we denote the basis $\gamma_0,\dots,\gamma_{2n+3}$ also by $1,\sfh,\dots,\sfh_{n},\epsilon_1,\dots,\epsilon_{n+3}$, and  use Einstein's summation convention, where the summation is from $0$ to $2n+3$.

In principle one can do the symmetric reduction of the WDVV equation for $X$. But the computation is quite complicated and we have not completed it. The following computations use the deformation invariance and the monodromy group in the same spirit as the symmetric reduction.

 \begin{lemma}\label{lem-3point-inv-even(2,2)}
\begin{equation}\label{eq-3point-inv-even(2,2)}
	\langle \sfh_{n-1},\sfh_{n-1},\sfh_n\rangle_{0,2}=192.
\end{equation}
\end{lemma}
\begin{proof}
By (\ref{eq-tauTot}),
\begin{eqnarray*}
&& \langle \sfh_{n-1},\sfh_{n-1},\sfh_n\rangle_{0,2}=\frac{\partial^3 F}{\partial t^{n-1}\partial t^{n-1}\partial t^n}(0)\\
&=& \sum_{0\leq i,j,k\leq n}\frac{\partial \tau^i}{\partial t^{n-1}}\frac{\partial \tau^j}{\partial t^{n-1}}\frac{\partial \tau^k}{\partial t^{n}}\frac{\partial^3 F}{\partial \tau^i \partial \tau^j\partial \tau^k}(0)\\
&=&\Big(\big(-4\frac{\partial}{\partial \tau^0}+\frac{\partial}{\partial \tau^{n-1}}\big)
\big(-4\frac{\partial}{\partial \tau^0}+\frac{\partial}{\partial \tau^{n-1}}\big)
\big(-12\frac{\partial}{\partial \tau^1}+\frac{\partial}{\partial \tau^{n}}\big)F\Big)|_{\tau=0}.
\end{eqnarray*}
Then from (\ref{eq-qp1.5}) we get (\ref{eq-3point-inv-even(2,2)}).
\end{proof}

\begin{lemma}\label{lem-5point-withTwoPrim-inv-even(2,2)}
\begin{equation}\label{eq-5point-withTwoPrim-inv-even(2,2)}
	\langle \epsilon_a,\epsilon_a,\sfh_n,\sfh_{n-1},\sfh_n\rangle_{0,3}=-192.
\end{equation}
\end{lemma}
\begin{proof}
By (\ref{eq-WDVV}) for $1\leq a\neq b\leq n+3$,
\begin{eqnarray}\label{lem-5point-withTwoPrim-inv-even(2,2)-1}
&& \langle \epsilon_a,\epsilon_a,\sfh_n,\sfh_{n-1},\gamma_e\rangle_0 g^{ef}\langle \gamma_f,\epsilon_b,\epsilon_b\rangle_0
+\langle \epsilon_a,\epsilon_a,\sfh_n,\gamma_e\rangle_0 g^{ef}\langle \gamma_f,\sfh_{n-1},\epsilon_b,\epsilon_b\rangle_0\nn\\
&&+\langle \epsilon_a,\epsilon_a,\sfh_{n-1},\gamma_e\rangle_0 g^{ef}\langle \gamma_f,\sfh_n,\epsilon_b,\epsilon_b\rangle_0
+\langle \epsilon_a,\epsilon_a,\gamma_e\rangle_0 g^{ef}\langle \gamma_f,\sfh_n,\sfh_{n-1},\epsilon_b,\epsilon_b\rangle_0\nn\\
&=& \langle \epsilon_a,\epsilon_b,\sfh_n,\sfh_{n-1},\gamma_e\rangle_0 g^{ef}\langle \gamma_f,\epsilon_a,\epsilon_b\rangle_0
+\langle \epsilon_a,\epsilon_b,\sfh_n,\gamma_e\rangle_0 g^{ef}\langle \gamma_f,\sfh_{n-1},\epsilon_a,\epsilon_b\rangle_0\nn\\
&&+\langle \epsilon_a,\epsilon_b,\sfh_{n-1},\gamma_e\rangle_0 g^{ef}\langle \gamma_f,\sfh_n,\epsilon_a,\epsilon_b\rangle_0
+\langle \epsilon_a,\epsilon_b,\gamma_e\rangle_0 g^{ef}\langle \gamma_f,\sfh_n,\sfh_{n-1},\epsilon_a,\epsilon_b\rangle_0.
\end{eqnarray}
By  (\ref{eq-Dim}), (\ref{eq-FCA}) and Theorem \ref{thm-monodromy-evenDim(2,2)}, the LHS of (\ref{lem-5point-withTwoPrim-inv-even(2,2)-1}) equals
\begin{eqnarray}\label{lem-5point-withTwoPrim-inv-even(2,2)-2}
&& \frac{1}{4}\langle \epsilon_a,\epsilon_a,\sfh_n,\sfh_{n-1},\sfh_n\rangle_{0,3} 
+\frac{1}{4}\langle \epsilon_a,\epsilon_a,\sfh_n,\sfh_{n-1},\sfh\rangle_{0,2} \langle \sfh_{n-1},\epsilon_b,\epsilon_b\rangle_{0,1}\nn\\
&&+\frac{1}{4}\langle \epsilon_a,\epsilon_a,\sfh_n,\sfh_{n-1}\rangle_{0,2} \langle \sfh,\sfh_{n-1},\epsilon_b,\epsilon_b\rangle_{0,1}
+\frac{1}{4}\langle \epsilon_a,\epsilon_a,\sfh_{n-1},\sfh\rangle_{0,1} \langle \sfh_{n-1},\sfh_n,\epsilon_b,\epsilon_b\rangle_{0,2}\nn\\
&&+\frac{1}{4}\langle \sfh_n,\sfh_n,\sfh_{n-1},\epsilon_b,\epsilon_b\rangle_{0,3}
+\frac{1}{4}\langle \epsilon_a,\epsilon_a,\sfh_{n-1}\rangle_{0,1}\langle \sfh,\sfh_n,\sfh_{n-1},\epsilon_b,\epsilon_b\rangle_{0,2},
\end{eqnarray}
and the RHS of (\ref{lem-5point-withTwoPrim-inv-even(2,2)-1}) is 0. 
By (\ref{eq-F^122}) and (\ref{eq-Div}),
\begin{equation}\label{lem-5point-withTwoPrim-inv-even(2,2)-3}
	\langle \epsilon_a,\epsilon_a,\sfh_{n-1},\sfh\rangle_{0,1}=\langle \epsilon_a,\epsilon_a,\sfh_{n-1}\rangle_{0,1}=-4,
\end{equation}
and
\begin{equation}\label{lem-5point-withTwoPrim-inv-even(2,2)-4}
	\langle \epsilon_a,\epsilon_a,\sfh_n,\sfh_{n-1}\rangle_{0,2}=-16.
\end{equation}
Then by (\ref{eq-3point-inv-even(2,2)}), (\ref{lem-5point-withTwoPrim-inv-even(2,2)-3}) and (\ref{lem-5point-withTwoPrim-inv-even(2,2)-4}), (\ref{lem-5point-withTwoPrim-inv-even(2,2)-1}) reads
\begin{eqnarray*}
&& \frac{1}{4}\langle \epsilon_a,\epsilon_a,\sfh_n,\sfh_{n-1},\sfh_n\rangle_{0,3} 
+\frac{1}{4}\times 2\times (-16)\times (-4)
+\frac{1}{4}\times (-16)\times (-4)\\
&&+\frac{1}{4}\times (-4)\times(-16)
+\frac{1}{4}\langle \sfh_n,\sfh_n,\sfh_{n-1},\epsilon_b,\epsilon_b\rangle_{0,3}
+\frac{1}{4}\times(-4)\times 2\times(-16)=0.
\end{eqnarray*}
By Theorem \ref{thm-monodromy-evenDim(2,2)}, 
\begin{equation*}
	\langle \epsilon_a,\epsilon_a,\sfh_n,\sfh_{n-1},\sfh_n\rangle_{0,3}=\langle \sfh_n,\sfh_n,\sfh_{n-1},\epsilon_b,\epsilon_b\rangle_{0,3}.
\end{equation*}
So we obtain (\ref{eq-5point-withTwoPrim-inv-even(2,2)}).
\end{proof}

\begin{lemma}\label{lem-5point-inv-even(2,2)-1}
\begin{equation}\label{eq-5point-inv-even(2,2)-1}
	\langle \epsilon_a,\epsilon_a,\epsilon_b,\epsilon_b,\sfh_n\rangle_{0,2}=4-4\langle \epsilon_a,\epsilon_b,\epsilon_a,\epsilon_b\rangle_{0,1}.
\end{equation}
\end{lemma}
\begin{proof}
By (\ref{eq-WDVV}), for $1\leq a\neq b\leq n+3$,
\begin{eqnarray*}
&& \langle \epsilon_a,\sfh_n,\epsilon_a,\epsilon_b,\gamma_e\rangle_0 g^{ef}\langle \gamma_f,\epsilon_b,\sfh_{n-1}\rangle_0
+\langle \epsilon_a,\sfh_n,\epsilon_a,\gamma_e\rangle_0 g^{ef}\langle \gamma_f,\epsilon_b,\epsilon_b,\sfh_{n-1}\rangle_0\\
&&+\langle \epsilon_a,\sfh_n,\epsilon_b,\gamma_e\rangle_0 g^{ef}\langle \gamma_f,\epsilon_a,\epsilon_b,\sfh_{n-1}\rangle_0
+\langle \epsilon_a,\sfh_n,\gamma_e\rangle_0 g^{ef}\langle \gamma_f,\epsilon_a,\epsilon_b,\epsilon_b,\sfh_{n-1}\rangle_0\\
&=& \langle \epsilon_a,\epsilon_b,\epsilon_a,\epsilon_b,\gamma_e\rangle_0 g^{ef}\langle \gamma_f,\sfh_n,\sfh_{n-1}\rangle_0
+\langle \epsilon_a,\epsilon_b,\epsilon_a,\gamma_e\rangle_0 g^{ef}\langle \gamma_f,\epsilon_b,\sfh_n,\sfh_{n-1}\rangle_0\\
&&+\langle \epsilon_a,\epsilon_b,\epsilon_b,\gamma_e\rangle_0 g^{ef}\langle \gamma_f,\epsilon_a,\sfh_n,\sfh_{n-1}\rangle_0
+\langle \epsilon_a,\epsilon_b,\gamma_e\rangle_0 g^{ef}\langle \gamma_f,\epsilon_a,\epsilon_b,\sfh_n,\sfh_{n-1}\rangle_0.
\end{eqnarray*}
Then (\ref{eq-Dim}), (\ref{eq-FCA}) and Theorem \ref{thm-monodromy-evenDim(2,2)} yield
\begin{eqnarray}\label{eq-lem-5point-inv-even(2,2)-1-1}
&& \langle \epsilon_a,\sfh_n,\epsilon_a,\epsilon_b,\epsilon_b\rangle_{0,2} \langle \epsilon_b,\epsilon_b,\sfh_{n-1}\rangle_{0,1}
+\frac{1}{4}\langle \epsilon_a,\sfh_n,\epsilon_a,\sfh_{n-1}\rangle_{0,2} \langle \sfh,\epsilon_b,\epsilon_b,\sfh_{n-1}\rangle_{0,1}\nn\\
&=&\frac{1}{4} \langle \epsilon_a,\epsilon_b,\epsilon_a,\epsilon_b,\sfh\rangle_{0,1} \langle \sfh_{n-1},\sfh_n,\sfh_{n-1}\rangle_{0,2}
+\langle \epsilon_a,\epsilon_b,\epsilon_a,\epsilon_b\rangle_{0,1} \langle \epsilon_b,\epsilon_b,\sfh_n,\sfh_{n-1}\rangle_{0,2}\nn\\
&&+\langle \epsilon_a,\epsilon_b,\epsilon_b,\epsilon_a\rangle_{0,1} \langle \epsilon_a,\epsilon_a,\sfh_n,\sfh_{n-1}\rangle_{0,2}.
\end{eqnarray}
Then by (\ref{eq-3point-inv-even(2,2)}), (\ref{lem-5point-withTwoPrim-inv-even(2,2)-3}) and (\ref{lem-5point-withTwoPrim-inv-even(2,2)-4}), (\ref{eq-lem-5point-inv-even(2,2)-1-1}) reads
\begin{eqnarray*}
&& -4\langle \epsilon_a,\sfh_n,\epsilon_a,\epsilon_b,\epsilon_b\rangle_{0,2}
+\frac{1}{4} (-16)(-4)\\
&=&\frac{1}{4} \langle \epsilon_a,\epsilon_b,\epsilon_a,\epsilon_b\rangle_{0,1} \times 192
-16\langle \epsilon_a,\epsilon_b,\epsilon_a,\epsilon_b\rangle_{0,1} -16\langle \epsilon_a,\epsilon_b,\epsilon_b,\epsilon_a\rangle_{0,1}.
\end{eqnarray*}
So we get (\ref{eq-5point-inv-even(2,2)-1}).
\end{proof}

\begin{lemma}\label{lem-5point-inv-even(2,2)-2}
\begin{equation}\label{eq-5point-inv-even(2,2)-2}
	 \langle \epsilon_a,\sfh_n,\epsilon_a,\epsilon_a,\epsilon_a\rangle_{0,2} 
= 12-4\langle \epsilon_a,\epsilon_a,\epsilon_a,\epsilon_a\rangle_{0,1}.
\end{equation}
\end{lemma}
\begin{proof}
By (\ref{eq-WDVV}), for $1\leq a\leq n+3$,
\begin{eqnarray}\label{eq-lem-5point-inv-even(2,2)-1}
&& \langle \epsilon_a,\sfh_n,\epsilon_a,\epsilon_a,\gamma_e\rangle_0 g^{ef}\langle \gamma_f,\epsilon_a,\sfh_{n-1}\rangle_0
+2\langle \epsilon_a,\sfh_n,\epsilon_a,\gamma_e\rangle_0 g^{ef}\langle \gamma_f,\epsilon_a,\epsilon_a,\sfh_{n-1}\rangle_0\nn\\
&&+\langle \epsilon_a,\sfh_n,\gamma_e\rangle_0 g^{ef}\langle \gamma_f,\epsilon_a,\epsilon_a,\epsilon_a,\sfh_{n-1}\rangle_0\nn\\
&=& \langle \epsilon_a,\epsilon_a,\epsilon_a,\epsilon_a,\gamma_e\rangle_0 g^{ef}\langle \gamma_f,\sfh_n,\sfh_{n-1}\rangle_0
+2\langle \epsilon_a,\epsilon_a,\epsilon_a,\gamma_e\rangle_0 g^{ef}\langle \gamma_f,\epsilon_a,\sfh_n,\sfh_{n-1}\rangle_0\nn\\
&&+\langle \epsilon_a,\epsilon_a,\gamma_e\rangle_0 g^{ef}\langle \gamma_f,\epsilon_a,\epsilon_a,\sfh_n,\sfh_{n-1}\rangle_0,
\end{eqnarray}
Then (\ref{eq-Dim}), (\ref{eq-FCA}) and Theorem \ref{thm-monodromy-evenDim(2,2)} yield
\begin{eqnarray*}
&& \langle \epsilon_a,\sfh_n,\epsilon_a,\epsilon_a,\epsilon_a\rangle_{0,2} \langle \epsilon_a,\epsilon_a,\sfh_{n-1}\rangle_{0,1}
+2\cdot\frac{1}{4}\langle \epsilon_a,\sfh_n,\epsilon_a,\sfh_{n-1}\rangle_{0,2} \langle \sfh,\epsilon_a,\epsilon_a,\sfh_{n-1}\rangle_{0,1}\\
&=&\frac{1}{4} \langle \epsilon_a,\epsilon_a,\epsilon_a,\epsilon_a,\sfh\rangle_{0,1}
\langle \sfh_{n-1},\sfh_n,\sfh_{n-1}\rangle_{0,2}
+2\langle \epsilon_a,\epsilon_a,\epsilon_a,\epsilon_a\rangle_{0,1} \langle \epsilon_a,\epsilon_a,\sfh_n,\sfh_{n-1}\rangle_{0,2}\\
&&+\frac{1}{4}\langle \epsilon_a,\epsilon_a,1\rangle_{0,0} \langle \sfh_n,\epsilon_a,\epsilon_a,\sfh_n,\sfh_{n-1}\rangle_{0,3}
+\frac{1}{4}\langle \epsilon_a,\epsilon_a,\sfh_{n-1}\rangle_{0,1} \langle \sfh,\epsilon_a,\epsilon_a,\sfh_n,\sfh_{n-1}\rangle_{0,2}.
\end{eqnarray*}
By (\ref{eq-3point-inv-even(2,2)}), (\ref{lem-5point-withTwoPrim-inv-even(2,2)-3}), (\ref{lem-5point-withTwoPrim-inv-even(2,2)-4}), (\ref{eq-5point-withTwoPrim-inv-even(2,2)})  and (\ref{eq-Div}),  (\ref{eq-lem-5point-inv-even(2,2)-1}) reads
\begin{eqnarray*}
&& -4\langle \epsilon_a,\sfh_n,\epsilon_a,\epsilon_a,\epsilon_a\rangle_{0,2} 
+2\cdot\frac{1}{4} (-16)(-4)\\
&=& 16 \langle \epsilon_a,\epsilon_a,\epsilon_a,\epsilon_a\rangle_{0,1}
+\frac{1}{4}\cdot (-192)
+\frac{1}{4}\times (-4)\times 2\times (-16),
\end{eqnarray*}
so we obtain (\ref{eq-5point-inv-even(2,2)-2}).
\end{proof}

\subsection{Correlators of length  4}\label{sec:correlators-length4}
By Theorem \ref{thm-monodromy-evenDim(2,2)} and the results recalled in Section \ref{sec:knownResults-correlators}, all the correlators of  length 4 other than the ones with only primitive insertions are computed. To compute the latter ones, we begin to recall a result on the length 4 correlators with only primitive insertions. Let $X$ be a smooth complete intersection of dimension $n$ and  Fano index $n-1$. Then $X$ is a cubic hypersurface i.e. the multidegree $\mathbf{d}=3$, or a complete intersection of two quadrics i.e. the multidegree $\mathbf{d}=(2,2)$. Let $m=\dim H^*_{\mathrm{prim}}(X)$, and  $\gamma_{n+1},\dots,\gamma_{n+m}$ be a basis of $H^*_{\mathrm{prim}}(X)$. Then by \cite[Proposition 9.12]{Hu15}, for $n+1\leq b,c\leq n+m$,
\begin{equation}\label{eq-4points-sum}
\sum_{e=n+1}^{n+m}\sum_{f=n+1}^{n+m}\langle \gamma_b, \gamma_e, g^{ef}\gamma_f,\gamma_c\rangle_{0,1}
=g_{bc}\cdot\begin{cases}
\frac{(-2)^{n+2}+8}{3},&
 \mbox{if}\ \mathbf{d}=3;\\
 (-1)^n(n+1)+2,
 &
 \mbox{if}\ \mathbf{d}=(2,2).
\end{cases}
\end{equation}

\begin{remark}\label{rem:eq-4points-sum}
For the reader's convenience, we briefly recall the proof of (\ref{eq-4points-sum}). Denote the multidegree by $\mathbf{d}=(d_1,\dots,d_r)$. Firstly, by the topological recursion relation in genus 1 \cite{Get98} and certain vanishing results from deformation invariance, we get
\begin{eqnarray}\label{eq-apply-TRR1-FanoIndex=n-1-2}
&&\langle \psi\gamma_b,\gamma_c\rangle_{1,1}\nn\\
&=&\frac{1}{\prod_{i=1}^r d_i}\langle \gamma_b,\gamma_c,\sfh_{n-1}\rangle_{0,1}\langle \sfh \rangle_{1,0}
+\frac{1}{\prod_{i=1}^r d_i}\langle \gamma_b,\gamma_c,1\rangle_{0,3,0}\langle \sfh_n \rangle_{1,1,1}\nn\\
&&+\frac{1}{24}\sum_{e=0}^{n+m}\sum_{f=0}^{n+m}\langle \gamma_b, \gamma_e, g^{ef}\gamma_f,\gamma_c\rangle_{0,1}.
\end{eqnarray}
Secondly, by Zinger's comparison formula \cite{Zin08} between genus 1 GW invariants, and reduced genus 1 GW invariants, and a vanishing result of reduced genus 1 GW invariants of degree 1, one can compute the genus 1 invariants in (\ref{eq-apply-TRR1-FanoIndex=n-1-2}). Then a lengthy computation in \cite[\S 9.3]{Hu15} yields (\ref{eq-4points-sum}).  This approach is essentially symplectic, for it uses Zinger's series of deep work on reduced genus 1 GW invariants.

An algebraic approach to (\ref{eq-4points-sum}) is to study the Fano variety of lines in $X$. We carry out this approach for cubic hypersurfaces in \cite[Chap. 8]{Hu15}. In principle one can do this for $\mathbf{d}=(2,2)$  in a similar way.
\end{remark}
When $X$ is a cubic hypersurface, or $X$ is an odd dimensional complete intersection of two quadrics, the Zariski closure of the  monodromy group is the orthogonal group or the symplectic group, according to the parity of the dimension. Then as we have seen in  \cite[Theorem 9.13]{Hu15}, (\ref{eq-4points-sum}) suffices to give all the length 4 correlators of $X$. In the exceptional case, i.e. $X$ is an even dimensional complete intersection of two quadrics, we need some ad hoc computations. 
\begin{proposition}\label{prop-4points-fanoIndex-even(2,2)}
Let   $X=X_n(2,2)$  of even dimension $n\geq 4$. Then
\begin{equation}\label{eq-4points-fanoIndex-even(2,2)-0}
	(n+5)\frac{\partial^2 F}{(\partial s_1)^2}(0)+\frac{\partial F}{\partial s_2}(0)
	=n+3.
\end{equation}
\end{proposition}
\begin{proof}
Take $\gamma_{n+1},\dots,\gamma_{2n+3}$ in (\ref{eq-4points-sum}) to be the orthonormal basis $\epsilon_1,\dots,\epsilon_{n+3}$ of $H^*_{\mathrm{prim}}(X)$. Then
\[
\sum_{e=n+1}^{n+m}\sum_{f=n+1}^{n+m}\langle \gamma_b, \gamma_e, g^{ef}\gamma_f,\gamma_c\rangle_{0,1}
=\sum_{e=1}^{n+3}\langle \epsilon_b, \epsilon_c,\epsilon_e,\epsilon_e\rangle_{0,1}.
\]
By (\ref{eq-invariantsOf-typeD-1}) and the definition (\ref{eq-def-generatingFunction}) of the generating function $F$, one finds
\begin{equation}\label{eq-prop-4points-fanoIndex-even(2,2)-1}
\langle \epsilon_b,\epsilon_b,\epsilon_e,\epsilon_e\rangle=\begin{cases}
\frac{\partial^2 F}{(\partial s_1)^2}(0),& \mbox{if}\ 1\leq e\leq n+3\ \mbox{and}\ e\neq b,\\
3\frac{\partial^2 F}{(\partial s_1)^2}(0)+\frac{\partial F}{\partial s_2}(0), &  \mbox{if}\ e=b.
\end{cases}
\end{equation}
So
\begin{equation}\label{eq-4points-sum-even(2,2)-1}
\sum_{e=n+1}^{n+m}\sum_{f=n+1}^{n+m}\langle \gamma_b, \gamma_e, g^{ef}\gamma_f,\gamma_c\rangle_{0,1}
=(n+5)\frac{\partial^2 F}{(\partial s_1)^2}(0) \delta_{b,c}+\frac{\partial F}{\partial s_2}(0)\delta_{b,c}.
\end{equation}
On the other hand by  (\ref{eq-4points-sum})  one finds
\begin{equation}\label{eq-4points-sum-even(2,2)-2}
\sum_{e=n+1}^{n+m}\sum_{f=n+1}^{n+m}\langle \gamma_b, \gamma_e, g^{ef}\gamma_f,\gamma_c\rangle_{0,1}
=(n+3)\delta_{b,c}.
\end{equation}
Comparing (\ref{eq-4points-sum-even(2,2)-1}) and (\ref{eq-4points-sum-even(2,2)-2}) we get (\ref{eq-4points-fanoIndex-even(2,2)-0}).
\end{proof}

Using (\ref{eq-WDVV}) and Theorem \ref{thm-monodromy-evenDim(2,2)}, we can get a quadratic equation for $\frac{\partial^2 F}{(\partial s_1)^2}(0)$ and $\frac{\partial F}{\partial s_2}(0)$. Then we determine the correct solution by the integrality of certain invariants.

\begin{lemma}\label{lem-4points-even(2,2)-integrality}
Let $X=X_n(2,2)$. Let $\alpha_i\in H^n_{\mathrm{prim}}(X)$ be the primitive classes corresponding to simple roots of $D_{n+3}$ as (\ref{eq-roots-D}), for $1\leq i\leq n+3$. Then
\begin{equation}\label{eq-4points-even(2,2)-integrality}
\langle \alpha_{i},\alpha_{j},\alpha_{k},\alpha_{l}\rangle_{0,1}\in \mathbb{Z},\ \mbox{for}\ 1\leq i,j,k,l\leq n+3.
\end{equation}
\end{lemma}
\begin{proof}
By \cite[Proposition 13.9]{Lew99}, for a general member $X$ in the family of $n$-dimensional smooth complete intersections of two quadrics, the Fano variety of lines $\Mbar_{0,0}(X,1)$
is a smooth scheme of dimension $2n-4$. It follows that $\Mbar_{0,0}(X,1)$ has the expected dimension. Moreover one directly sees that $\Mbar_{0,0}(X,1)$ has no stacky point. Since $\alpha_i\in H^n(X;\mathbb{Z})$, we have (\ref{eq-4points-even(2,2)-integrality}) for such $X$. Then (\ref{eq-4points-even(2,2)-integrality}) holds for all $n$-dimensional smooth complete intersections of two quadrics, by the deformation invariance. 
\end{proof}

\begin{remark}\label{rem:integreality-symplectic}
One can also directly apply the integrality of genus 0 Gromov-Witten invariants with integral classes as insertions, of semipositive symplectic manifolds (e.g. \cite[Theorem A]{Ruan96}).
\end{remark}

\begin{theorem}\label{thm-4points-fanoIndex-even(2,2)}
Let   $X$ be an even dimensional complete intersection of two quadrics in $\mathbb{P}^{n+2}$,  with $n\geq 4$. Then
\begin{equation}\label{eq-4points-fanoIndex-even(2,2)}
	\frac{\partial^2 F}{(\partial s_1)^2}(0)=1,\ \frac{\partial F}{\partial s_2}(0)=-2.
\end{equation}
Equivalently, for $1\leq a,b\leq n+3$,
\begin{equation}\label{eq-4points-fanoIndex-even(2,2)-ab}
	\langle \epsilon_a,\epsilon_a,\epsilon_b,\epsilon_b\rangle_{0,1}=1.
\end{equation}
\end{theorem}
\begin{proof}
By (\ref{eq-WDVV}), for $n+1\leq a\neq b\leq 2n+3$,
\begin{eqnarray}\label{eq-thm-4points-fanoIndex-even(2,2)-1}
&&\sum_{e=0}^{2n+3}\sum_{f=0}^{2n+3}\big(\langle \epsilon_a,\epsilon_a,\epsilon_a,\epsilon_a,\gamma_e\rangle_0 g^{ef}\langle \gamma_f,\epsilon_b,\epsilon_b\rangle_0+2\langle \epsilon_a,\epsilon_a,\epsilon_a,\gamma_e\rangle_0 g^{ef}\langle \gamma_f,\epsilon_a,\epsilon_b,\epsilon_b\rangle_0
\nn\\
&&+\langle \epsilon_a,\epsilon_a,\gamma_e\rangle_0 g^{ef}\langle \gamma_f,\epsilon_a,\epsilon_a,\epsilon_b,\epsilon_b\rangle_0\big)\nn\\
&=&\sum_{e=0}^{2n+3}\sum_{f=0}^{2n+3}\big(2 \langle \epsilon_a,\epsilon_b,\epsilon_a,\gamma_e\rangle_0 g^{ef}\langle \gamma_f,\epsilon_a,\epsilon_a,\epsilon_b\rangle_0
+2 \langle \epsilon_a,\epsilon_b,\epsilon_a,\epsilon_a,\gamma_e\rangle_0 g^{ef}\langle \gamma_f,\epsilon_a,\epsilon_b\rangle_0\big).\nn\\
\end{eqnarray}
By Theorem \ref{thm-monodromy-evenDim(2,2)} and (\ref{eq-prop-4points-fanoIndex-even(2,2)-1}), the RHS of (\ref{eq-thm-4points-fanoIndex-even(2,2)-1})  equals
\begin{eqnarray}\label{eq-thm-4points-fanoIndex-even(2,2)-2}
2 \langle \epsilon_a,\epsilon_b,\epsilon_a,\epsilon_b\rangle_{0,1} \langle \epsilon_b,\epsilon_a,\epsilon_a,\epsilon_b\rangle_{0,1}=2\big(\frac{\partial^2 F}{\partial s_1^2}(0)
\big)^2.
\end{eqnarray}
By (\ref{eq-Dim}), (\ref{eq-FCA}) and Theorem \ref{thm-monodromy-evenDim(2,2)}, the LHS (\ref{eq-thm-4points-fanoIndex-even(2,2)-1})  equals
\begin{eqnarray}\label{eq-thm-4points-fanoIndex-even(2,2)-3}
&&\frac{1}{4}\langle \epsilon_a,\epsilon_a,\epsilon_a,\epsilon_a,\sfh\rangle_{0,1} \langle \sfh_{n-1},\epsilon_b,\epsilon_b\rangle_{0,1}
+\frac{1}{4}\langle \epsilon_a,\epsilon_a,\epsilon_a,\epsilon_a,\sfh_n\rangle_{0,2} \langle 1,\epsilon_b,\epsilon_b\rangle_{0,0}\nn\\
&&+2\langle \epsilon_a,\epsilon_a,\epsilon_a,\epsilon_a\rangle_{0,1} \langle \epsilon_a,\epsilon_a,\epsilon_b,\epsilon_b\rangle_{0,1}
+\frac{1}{4}\langle \epsilon_a,\epsilon_a,1\rangle_{0,0} \langle \sfh_{n},\epsilon_a,\epsilon_a,\epsilon_b,\epsilon_b\rangle_{0,2}\nn\\
&&+\frac{1}{4}\langle \epsilon_a,\epsilon_a,\sfh_{n-1}\rangle_{0,1} \langle \sfh,\epsilon_a,\epsilon_a,\epsilon_b,\epsilon_b\rangle_{0,1}.
\end{eqnarray}
 Since $\epsilon_a$ and $\epsilon_b$ are chosen to be orthonormal, we have
\begin{equation}\label{eq-thm-4points-fanoIndex-even(2,2)-4}
	1=\langle 1,\epsilon_a,\epsilon_a\rangle_{0,0}=\langle 1,\epsilon_b,\epsilon_b\rangle_{0,0}.
\end{equation}
By (\ref{eq-F^122}),
\begin{equation}\label{eq-thm-4points-fanoIndex-even(2,2)-5}
	\langle \epsilon_a,\epsilon_a,\sfh_{n-1}\rangle_{0,1}=\langle \epsilon_b,\epsilon_b,\sfh_{n-1}\rangle_{0,1}=-4.
\end{equation}
By (\ref{eq-prop-4points-fanoIndex-even(2,2)-1}),
\begin{equation}\label{eq-thm-4points-fanoIndex-even(2,2)-6}
	\langle \epsilon_a,\epsilon_a,\epsilon_a,\epsilon_a\rangle_{0,1}=3\frac{\partial^2 F}{\partial s_1^2}(0)+\frac{\partial F}{\partial s_2}(0).
\end{equation}
Using (\ref{eq-Div}),  by (\ref{eq-thm-4points-fanoIndex-even(2,2)-4}), (\ref{eq-thm-4points-fanoIndex-even(2,2)-5}), (\ref{eq-thm-4points-fanoIndex-even(2,2)-6}),
 and  (\ref{eq-5point-inv-even(2,2)-1}), (\ref{eq-5point-inv-even(2,2)-2}) in Section \ref{sec:preparatoryComputation-even(2,2)}, (\ref{eq-thm-4points-fanoIndex-even(2,2)-3}) equals
\begin{eqnarray*}
 6\big(\frac{\partial^2 F}{\partial s_1^2}(0)\big)^2+2\frac{\partial F}{\partial s_2}(0)\frac{\partial^2 F}{\partial s_1^2}(0)-8\frac{\partial^2 F}{\partial s_1^2}(0)
-2 \frac{\partial F}{\partial s_2}(0)+4.
\end{eqnarray*}
So (\ref{eq-thm-4points-fanoIndex-even(2,2)-1})  yields
\begin{equation}\label{eq-thm-4points-fanoIndex-even(2,2)-7}
	4\big(\frac{\partial^2 F}{\partial s_1^2}(0)\big)^2+2\frac{\partial F}{\partial s_2}(0)\frac{\partial^2 F}{\partial s_1^2}(0)-8\frac{\partial^2 F}{\partial s_1^2}(0)
-2 \frac{\partial F}{\partial s_2}(0)+4=0.
\end{equation}
Substituting (\ref{eq-4points-fanoIndex-even(2,2)-0}) into (\ref{eq-thm-4points-fanoIndex-even(2,2)-7}),
we obtain
\begin{equation*}
	\big((n+3)\frac{\partial^2 F}{\partial s_1^2}(0)-(n+1)\big)(\frac{\partial^2 F}{\partial s_1^2}(0)-1)=0.
\end{equation*}
Thus
\begin{equation*}
	\frac{\partial^2 F}{\partial s_1^2}(0)=\frac{n+1}{n+3}\ \mbox{or}\ 1.
\end{equation*}
Recall the simple roots (\ref{eq-roots-D}). We have
\begin{equation*}
\langle \alpha_{1},\alpha_{1},\alpha_{n+3},\alpha_{n+3}\rangle_{0,1}=\frac{\partial^2 F}{\partial s_1^2}(0)\cdot (\alpha_{1},\alpha_{1})(\alpha_{n+3},\alpha_{n+3})=4 \frac{\partial^2 F}{\partial s_1^2}(0).
\end{equation*}
But by Lemma \ref{lem-4points-even(2,2)-integrality}, $\langle \alpha_{1},\alpha_{1},\alpha_{n+3},\alpha_{n+3}\rangle_{0,1}\in \mathbb{Z}$. Since $n$ is even, $\frac{4(n+1)}{n+3}$ is never an integer. Hence
\begin{equation*}
	\frac{\partial^2 F}{\partial s_1^2}(0)=1,
\end{equation*}
and by (\ref{eq-4points-fanoIndex-even(2,2)-0}) we get $\frac{\partial F}{\partial s_2}(0)=-2$.
The formula (\ref{eq-4points-fanoIndex-even(2,2)-ab}) follows then from (\ref{eq-prop-4points-fanoIndex-even(2,2)-1}).
\end{proof}

\section{A reconstruction theorem}\label{sec:reconstructionTheorem}
In this section, we fix a smooth complete intersection $X$ of two quadrics in $\mathbb{P}^{n+2}$, where $n$ is even and $\geq 4$. 
This section aims to show that, with the results in Section \ref{sec:correlators-lengt-atMost4}, besides  a special GW invariant,  we can compute all genus zero GW invariants of $X$ by (\ref{eq-Dim}), (\ref{eq-EulerVectorField}), (\ref{eq-WDVV}), and the deformation invariance. We begin with an easy observation.
By (\ref{eq-Dim}), a length $k$ genus 0 GW invariant of $X$ with only primitive insertions is zero unless
\begin{equation*}
	k\cdot \frac{n}{2}=n-3+k+\beta\cdot(n-1),
\end{equation*}
i.e.
\[
\beta=\beta(k):=\frac{\frac{k(n-2)}{2}-n+3}{n-1}.
\]
This is an integer if and only if $n-1$ divides $k-4$. In particular
\begin{equation*}
	 \beta(4)=1,\ \beta(n+3)=\frac{n}{2},\ \beta(2n+2)=n-1.
\end{equation*}

For the brevity of expressions, we introduce some notations. For $0\leq j\leq 2n+3$, we set 
\[
\partial_{t^j}=\frac{\partial}{\partial t^j}.
\]
For $I=\{i_0,i_1,\dots,i_{2n+3}\}\in \mathbb{Z}_{\geq 0}^{2n+4}$, we define
\[
\partial_{t^I}=(\partial_{t^0})^{i_0}\circ\dots\circ (\partial_{t^{2n+3}})^{i_{2n+3}}.
\]
For $0\leq j\leq 2n+3$, let $e_j$ be the $(j+1)$-th unit vector in $\mathbb{Z}_{\geq 0}^{2n+4}$. So
\[
\partial_{t^{I+e_j}}=\partial_{t^I}\circ \partial_{t^j}.
\]
We apply similar notations to the coordinates $\tau^0,\dots,\tau^{2n+3}$. 

For $I=(i_0,\dots,i_{2n+3})$ in $\mathbb{Z}_{\geq 0}^{2n+4}$, we set $|I|=\sum_{k=0}^{2n+3}i_k$, and $I!=\prod_{k=0}^{2n+3}i_k!$. For $I=(i_0,\dots,i_{2n+3})$ and $J=(j_0,\dots,j_{2n+3})$ in $\mathbb{Z}_{\geq 0}^{2n+4}$, we say $I\leq J$ if and only if  $i_k\leq j_k$ for $0\leq k\leq 2n+3$. We denote an zero vector $(0,\dots,0)$ also by $0$ when no confusion arises. Moreover we define
\begin{equation}\label{eq-binomOfLists}
\binom{I}{J}=\prod_{k=0}^{2n+3} \binom{i_k}{j_k}.
 \end{equation}

In the following subsections we adopt Einstein's summation convention, where the range of the indices runs over $0,\dots,2n+3$.

\subsection{Elimination of ambient classes}
  By \cite[Theorem 3.1]{KM94}, a correlator with only ambient insertions can be computed from the length 3 correlators with only ambient correlators and (\ref{eq-WDVV}). In  \cite[Appendix D]{Hu15} we present an explicit algorithm in the $\tau$-coordinates for any Fano complete intersections in projective spaces.
  The system of $\tau$-coordinates has the advantage that the linear recursion of the highest order terms in the WDVV equations is quite simple. The cost is that the expression of the Euler field becomes complicated. For the even dimensional intersections of two quadrics in consideration, by (\ref{eq-tauTot}), the $\tau$-coordinates are very close to the $t$-coordinates. Nevertheless we still work in  $\tau$-coordinates in this section. 
  We will show how to eliminate an ambient class in any correlator by (\ref{eq-EulerVectorField}) and (\ref{eq-WDVV}). 
\begin{lemma}\label{lem-recursion-EulerVecField-even(2,2)}
Let $I\in \mathbb{Z}_{\geq 0}^{2n+4}$, and suppose $|I|\geq 4$. Then
\begin{eqnarray}\label{eq-recursion-EulerVecField-even(2,2)}
\partial_{\tau^1}\partial_{\tau^I}F(0)
&=& \frac{\sum_{j=0}^n(j-1)i_j+(\frac{n}{2}-1)\sum_{j=n+1}^{2n+3}i_j+3-n}{n-1}\partial_{\tau^I}F(0)\nn\\
&&-12i_n\partial_{\tau^1}\partial_{\tau^{I-e_n}}F(0).
\end{eqnarray}
\end{lemma}
\begin{proof}
By (\ref{eq-tauTot}) and (\ref{eq-tTotau}), we write the Euler vector field (\ref{eq-EV-0}) in the $\tau$-coordinates:
\begin{eqnarray}\label{eq-EulerField-tau-Coordinates}
E&=& \sum_{i=0}^{n}(1-i)t^{i}\frac{\partial}{\partial t^i}+\sum_{i=n+1}^{2n+3}(1-\frac{n}{2})t^{i}\frac{\partial}{\partial t^i}+(n-1)\frac{\partial}{\partial t^1}\nn\\
&=&\sum_{i=0}^{n}(1-i)\tau^i\frac{\partial}{\partial \tau^i}+(4n-4) \tau^{n-1}\frac{\partial}{\partial \tau^{0}}
	+(12n-12) \tau^n\frac{\partial}{\partial \tau^1}\nn\\
&&+\sum_{i=n+1}^{2n+3}(1-\frac{n}{2})\tau^{i}\frac{\partial}{\partial \tau^i}
+(n-1)\frac{\partial}{\partial \tau^1}.
\end{eqnarray}
Recall (\ref{eq-EulerVectorField})
\[
E F=(3-n)F+(n-1)\partial_{t^1}c.
\]
Since $|I|\geq 4$, we have
\begin{eqnarray*}
&&(n-1)\partial_{\tau^1}\partial_{\tau^I}F(0)\\
&=&-\sum_{j=0}^n(1-j)i_j\partial_{\tau^I}F(0)-(12n-12)i_n\partial_{\tau^1}\partial_{\tau^{I-e_n}}F(0)\\
&&-\sum_{j=n+1}^{2n+3}(1-\frac{n}{2})i_j \partial_{\tau^I}F(0)+(3-n)\partial_{\tau^I}F(0)\\
&=& \big(-\sum_{j=0}^n(1-j)i_j-\sum_{j=n+1}^{2n+3}(1-\frac{n}{2})i_j+3-n\big)\partial_{\tau^I}F(0)\\
&&-(12n-12)i_n\partial_{\tau^1}\partial_{\tau^{I-e_n}}F(0).
\end{eqnarray*}
So we obtain (\ref{eq-recursion-EulerVecField-even(2,2)}).
\end{proof}

\begin{lemma}\label{lem-recursion-ambient-simplified-even(2,2)}
Let $I\in \mathbb{Z}^{2n+4}_{\geq 0}$.
Let $n+1\leq a,b\leq 2n+3$, and $2\leq i\leq n$. Then
\begin{eqnarray}\label{eq-recursion-ambient-simplified-even(2,2)}
&&\partial_{\tau^i}\partial_{\tau^I}\partial_{\tau^a}\partial_{\tau^b}F(0)\nn\\
&=& 
4 \delta_{a,b}\partial_{\tau^1}^2\partial_{\tau^{i-1}}\partial_{\tau^I}F(0)	\nn\\
&&-\frac{1}{4} \sum_{\begin{subarray}{c}0\leq J\leq I\\ 1\leq |J|\leq |I|\end{subarray}}\sum_{e=0}^n
\binom{I}{J}\partial_{\tau^1}\partial_{\tau^{i-1}}\partial_{\tau^J}\partial_{\tau^e}F(0)
	\partial_{\tau^{n-e}}\partial_{\tau^{I-J}}\partial_{\tau^a}\partial_{\tau^b}F(0)\nn\\
&&- \sum_{\begin{subarray}{c}0\leq J\leq I\\ 1\leq |J|\leq |I|\end{subarray}}\sum_{e=n+1}^{2n+3}
\binom{I}{J}\partial_{\tau^1}\partial_{\tau^{i-1}}\partial_{\tau^J}\partial_{\tau^e}F(0)
	\partial_{\tau^{e}}\partial_{\tau^{I-J}}\partial_{\tau^a}\partial_{\tau^b}F(0)\nn\\	
&&+	\frac{1}{4}\sum_{\begin{subarray}{c}0\leq J\leq I\\ 1\leq |J|\leq |I|-1\end{subarray}}
	\sum_{e=0}^n
\binom{I}{J}\partial_{\tau^1}\partial_{\tau^{a}}\partial_{\tau^J}\partial_{\tau^e}F(0)
	\partial_{\tau^{n-e}}\partial_{\tau^{I-J}}\partial_{\tau^{i-1}}\partial_{\tau^b}F(0)\nn\\
&&+\sum_{\begin{subarray}{c}0\leq J\leq I\\ 1\leq |J|\leq |I|-1\end{subarray}}
	\sum_{e=n+1}^{2n+3}
\binom{I}{J}\partial_{\tau^1}\partial_{\tau^{a}}\partial_{\tau^J}\partial_{\tau^e}F(0)
	\partial_{\tau^{e}}\partial_{\tau^{I-J}}\partial_{\tau^{i-1}}\partial_{\tau^b}F(0).	
\end{eqnarray}
\end{lemma}
\begin{proof}
Consider the WDVV equation 
\begin{equation}\label{eq-lem-recursion-ambient-simplified-even(2,2)-WDVV-1}
	(\partial_{\tau^1}\partial_{\tau^{i-1}}\partial_{\tau^{e}}F)\eta^{ef}(\partial_{\tau^f}\partial_{\tau^{a}}\partial_{\tau^{b}}F)
	=(\partial_{\tau^1}\partial_{\tau^{a}}\partial_{\tau^{e}}F)\eta^{ef}(\partial_{\tau^f}\partial_{\tau^{i-1}}\partial_{\tau^{b}}F).
\end{equation}
We apply the differential operator $\partial_{\tau^I}$ to (\ref{eq-lem-recursion-ambient-simplified-even(2,2)-WDVV-1}), and then take the constant terms of both sides. We obtain
\begin{eqnarray}\label{eq-lem-recursion-ambient-simplified-even(2,2)-WDVV-2}
 &&\sum_{0\leq J\leq I}\binom{I}{J}\partial_{\tau^1}\partial_{\tau^{i-1}}\partial_{\tau^J}\partial_{\tau^e}F(0)\eta^{ef}
	\partial_{\tau^{f}}\partial_{\tau^{I-J}}\partial_{\tau^a}\partial_{\tau^b}F(0)\nn\\
&=&\sum_{0\leq J\leq I}\binom{I}{J}\partial_{\tau^1}\partial_{\tau^{a}}\partial_{\tau^J}\partial_{\tau^e}F(0)\eta^{ef}
	\partial_{\tau^{f}}\partial_{\tau^{I-J}}\partial_{\tau^{i-1}}\partial_{\tau^b}F(0).
\end{eqnarray}
 Note that 
\[
\partial_{\tau^j}\partial_{\tau^{k}}\partial_{\tau^e}F(0)\eta^{ef}
\]
are the structure constants of the small quantum multiplications. By definition of $\tsfh_{i}$ we have
\[
\tsfh_{1}\diamond \tsfh_{i-1}=\tsfh_{i},
\]
and thus
\[
\partial_{\tau^1}\partial_{\tau^{i-1}}\partial_{\tau^e}F(0)\eta^{ef}=\delta_{i,f}.
\]
So the LHS of (\ref{eq-lem-recursion-ambient-simplified-even(2,2)-WDVV-2}) equals
\begin{eqnarray}\label{eq-lem-recursion-ambient-simplified-even(2,2)-WDVV-3}
&&  \partial_{\tau^1}\partial_{\tau^{i-1}}\partial_{\tau^e}F(0)\eta^{ef}\partial_{\tau^f}\partial_{\tau^I}\partial_{\tau^a}\partial_{\tau^b}F(0)\nn\\
&&+ \sum_{\begin{subarray}{c}0\leq J\leq I\\ 1\leq |J|\leq |I|\end{subarray}}\binom{I}{J}\partial_{\tau^1}\partial_{\tau^{i-1}}\partial_{\tau^J}\partial_{\tau^e}F(0)\eta^{ef}
	\partial_{\tau^{f}}\partial_{\tau^{I-J}}\partial_{\tau^a}\partial_{\tau^b}F(0)\nn\\
&=& \partial_{\tau^i}\partial_{\tau^I}\partial_{\tau^a}\partial_{\tau^b}F(0)
	+ \sum_{\begin{subarray}{c}0\leq J\leq I\\ 1\leq |J|\leq |I|\end{subarray}}\binom{I}{J}\partial_{\tau^1}\partial_{\tau^{i-1}}\partial_{\tau^J}\partial_{\tau^e}F(0)\eta^{ef}
	\partial_{\tau^{f}}\partial_{\tau^{I-J}}\partial_{\tau^a}\partial_{\tau^b}F(0).
\end{eqnarray}
The RHS of (\ref{eq-lem-recursion-ambient-simplified-even(2,2)-WDVV-2}) equals
\begin{eqnarray*}
&& \partial_{\tau^1}\partial_{\tau^{a}}\partial_{\tau^e}F(0)\eta^{ef}\partial_{\tau^f}\partial_{\tau^I}\partial_{\tau^b}\partial_{\tau^{i-1}}F(0)
	+\partial_{\tau^1}\partial_{\tau^{a}}\partial_{\tau^I}\partial_{\tau^e}F(0)\eta^{ef}\partial_{\tau^f}\partial_{\tau^b}\partial_{\tau^{i-1}}F(0)\\
&& +	\sum_{\begin{subarray}{c}0\leq J\leq I\\ 1\leq |J|\leq |I|-1\end{subarray}}\binom{I}{J}\partial_{\tau^1}\partial_{\tau^{a}}\partial_{\tau^J}\partial_{\tau^e}F(0)\eta^{ef}
	\partial_{\tau^{f}}\partial_{\tau^{I-J}}\partial_{\tau^{i-1}}\partial_{\tau^b}F(0)\\
&\stackrel{\mbox{\footnotesize{by} }(\ref{eq-etaInversePairing-even(2,2)})}{=}&
\partial_{\tau^1}\partial_{\tau^{a}}^2F(0)\partial_{\tau^I}\partial_{\tau^a}\partial_{\tau^b}\partial_{\tau^{i-1}}F(0)
	+\partial_{\tau^1}\partial_{\tau^{a}}\partial_{\tau^b}\partial_{\tau^I}F(0)\partial_{\tau^b}^2\partial_{\tau^{i-1}}F(0)\\
&& +	\sum_{\begin{subarray}{c}0\leq J\leq I\\ 1\leq |J|\leq |I|-1\end{subarray}}\binom{I}{J}\partial_{\tau^1}\partial_{\tau^{a}}\partial_{\tau^J}\partial_{\tau^e}F(0)\eta^{ef}
	\partial_{\tau^{f}}\partial_{\tau^{I-J}}\partial_{\tau^{i-1}}\partial_{\tau^b}F(0).
\end{eqnarray*}
By (\ref{eq-F122-tau}),
\[
\partial_{\tau^1}\partial_{\tau^{a}}^2F(0)=0=\partial_{\tau^b}^2\partial_{\tau^{i-1}}F(0).
\]	
So the RHS of (\ref{eq-lem-recursion-ambient-simplified-even(2,2)-WDVV-2}) equals
\begin{equation}\label{eq-lem-recursion-ambient-simplified-even(2,2)-WDVV-4}
\sum_{\begin{subarray}{c}0\leq J\leq I\\ 1\leq |J|\leq |I|-1\end{subarray}}\binom{I}{J}\partial_{\tau^1}\partial_{\tau^{a}}\partial_{\tau^J}\partial_{\tau^e}F(0)\eta^{ef}
	\partial_{\tau^{f}}\partial_{\tau^{I-J}}\partial_{\tau^{i-1}}\partial_{\tau^b}F(0).
\end{equation}
Putting together (\ref{eq-lem-recursion-ambient-simplified-even(2,2)-WDVV-2}), (\ref{eq-lem-recursion-ambient-simplified-even(2,2)-WDVV-3}) and (\ref{eq-lem-recursion-ambient-simplified-even(2,2)-WDVV-4}) yields
\begin{eqnarray}\label{eq-recursion-ambient-even(2,2)}
\partial_{\tau^i}\partial_{\tau^I}\partial_{\tau^a}\partial_{\tau^b}F(0)
&=&- \sum_{\begin{subarray}{c}0\leq J\leq I\\ 1\leq |J|\leq |I|\end{subarray}}\binom{I}{J}\partial_{\tau^1}\partial_{\tau^{i-1}}\partial_{\tau^J}\partial_{\tau^e}F(0)\eta^{ef}
	\partial_{\tau^{f}}\partial_{\tau^{I-J}}\partial_{\tau^a}\partial_{\tau^b}F(0)\nn\\
&&+	\sum_{\begin{subarray}{c}0\leq J\leq I\\ 1\leq |J|\leq |I|-1\end{subarray}}\binom{I}{J}\partial_{\tau^1}\partial_{\tau^{a}}\partial_{\tau^J}\partial_{\tau^e}F(0)\eta^{ef}
	\partial_{\tau^{f}}\partial_{\tau^{I-J}}\partial_{\tau^{i-1}}\partial_{\tau^b}F(0).
\end{eqnarray}
Finally applying (\ref{eq-etaInversePairing-even(2,2)}) to the RHS of (\ref{eq-recursion-ambient-even(2,2)}) we obtain (\ref{eq-recursion-ambient-simplified-even(2,2)}).
\end{proof}

\subsection{Elimination of primitive classes}

\begin{lemma}\label{lem-recursion-primitive-abcc-even(2,2)}
Let $I=(i_{n+1},\dots,i_{2n+3})$ be a list  indicating the number of insertions of $\epsilon_{n+1},\dots,\epsilon_{2n+3}$. Suppose $n+1\leq a,b\leq 2n+3$.
Then
\begin{eqnarray}\label{eq-recursion-primitive-abcc-even(2,2)-3}
&&\partial_{\tau^n}\partial_{\tau^I}\partial_{\tau^a}\partial_{\tau^b}F(0)
=4|I|\partial_{\tau^{a}}\partial_{\tau^b}\partial_{\tau^I}F(0)\nn\\
&&-\sum_{\begin{subarray}{c}0\leq J\leq I\\ 2\leq |J|\leq |I|\end{subarray}}\binom{I}{J}\partial_{\tau^1}\partial_{\tau^{n-1}}\partial_{\tau^J}\partial_{\tau^e}F(0)\eta^{ef}
	\partial_{\tau^{f}}\partial_{\tau^{I-J}}\partial_{\tau^a}\partial_{\tau^b}F(0) \nn\\
&&+	\sum_{\begin{subarray}{c}0\leq J\leq I\\ 1\leq |J|\leq |I|-1\end{subarray}}\binom{I}{J}\partial_{\tau^1}\partial_{\tau^{a}}\partial_{\tau^J}\partial_{\tau^e}F(0)\eta^{ef}
	\partial_{\tau^{f}}\partial_{\tau^{I-J}}\partial_{\tau^{n-1}}\partial_{\tau^b}F(0).
\end{eqnarray}
\end{lemma}
\begin{proof}
We make use of the WDVV equation
\begin{equation}\label{eq-WDVV-recursion-primitive-abcc-even(2,2)-3}
	(\partial_{\tau^1}\partial_{\tau^{n-1}}\partial_{\tau^e}F)\eta^{ef}(\partial_{\tau^f}\partial_{\tau^a}\partial_{\tau^b}F)=(\partial_{\tau^1}\partial_{\tau^a}\partial_{\tau^e}F)\eta^{ef}(\partial_{\tau^f}\partial_{\tau^{n-1}}\partial_{\tau^b}F).
\end{equation}
The coefficient of $\tau^I$ of the LHS of (\ref{eq-WDVV-recursion-primitive-abcc-even(2,2)-3}), after multiplying $I!$, is
\begin{eqnarray*}
&&  \sum_{0\leq J\leq I}\binom{I}{J}\partial_{\tau^1}\partial_{\tau^{n-1}}\partial_{\tau^J}\partial_{\tau^e}F(0)\eta^{ef}
	\partial_{\tau^{f}}\partial_{\tau^{I-J}}\partial_{\tau^a}\partial_{\tau^b}F(0)\\
&\stackrel{\mbox{by (\ref{eq-F122-tau}) and (\ref{eq-etaInversePairing-even(2,2)})} }{=}& \partial_{\tau^1}\partial_{\tau^{n-1}}\partial_{\tau^e}F(0)\eta^{ef}\partial_{\tau^f}\partial_{\tau^I}\partial_{\tau^a}\partial_{\tau^b}F(0)\\
&&	+\sum_{j=n+1}^{2n+3}i_j \partial_{\tau^1}\partial_{\tau^{n-1}}\partial_{\tau^j}^2F(0)\partial_{\tau^I}\partial_{\tau^a}\partial_{\tau^b}F(0)\\
&&+ \sum_{\begin{subarray}{c}0\leq J\leq I\\ 2\leq |J|\leq |I|\end{subarray}}\binom{I}{J}\partial_{\tau^1}\partial_{\tau^{n-1}}\partial_{\tau^J}\partial_{\tau^e}F(0)\eta^{ef}
	\partial_{\tau^{f}}\partial_{\tau^{I-J}}\partial_{\tau^a}\partial_{\tau^b}F(0)\\
&\stackrel{\mbox{by (\ref{eq-F122-tau}) and $\tsfh_1\sqp \tsfh_{n-1}=\tsfh_{n}$ }}{=}& \partial_{\tau^n}\partial_{\tau^I}\partial_{\tau^a}\partial_{\tau^b}F(0)-4|I|\partial_{\tau^{a}}\partial_{\tau^b}\partial_{\tau^I}F(0)\\
&&+ \sum_{\begin{subarray}{c}0\leq J\leq I\\ 2\leq |J|\leq |I|\end{subarray}}\binom{I}{J}\partial_{\tau^1}\partial_{\tau^{n-1}}\partial_{\tau^J}\partial_{\tau^e}F(0)\eta^{ef}
	\partial_{\tau^{f}}\partial_{\tau^{I-J}}\partial_{\tau^a}\partial_{\tau^b}F(0).
\end{eqnarray*}
The coefficient of $\tau^I$ of the RHS of (\ref{eq-WDVV-recursion-primitive-abcc-even(2,2)-3}), after multiplying $I!$, is
\begin{eqnarray*}
&&  \sum_{0\leq J\leq I}\binom{I}{J}\partial_{\tau^1}\partial_{\tau^{a}}\partial_{\tau^J}\partial_{\tau^e}F(0)\eta^{ef}
	\partial_{\tau^{f}}\partial_{\tau^{I-J}}\partial_{\tau^{n-1}}\partial_{\tau^b}F(0)\\
&=& \partial_{\tau^1}\partial_{\tau^{a}}\partial_{\tau^e}F(0)\eta^{ef}\partial_{\tau^f}\partial_{\tau^I}\partial_{\tau^b}\partial_{\tau^{n-1}}F(0)
	+\partial_{\tau^1}\partial_{\tau^{a}}\partial_{\tau^I}\partial_{\tau^e}F(0)\eta^{ef}\partial_{\tau^f}\partial_{\tau^b}\partial_{\tau^{n-1}}F(0)\\
&& +	\sum_{\begin{subarray}{c}0\leq J\leq I\\ 1\leq |J|\leq |I|-1\end{subarray}}\binom{I}{J}\partial_{\tau^1}\partial_{\tau^{a}}\partial_{\tau^J}\partial_{\tau^e}F(0)\eta^{ef}
	\partial_{\tau^{f}}\partial_{\tau^{I-J}}\partial_{\tau^{n-1}}\partial_{\tau^b}F(0)\\
&\stackrel{\mbox{by (\ref{eq-F122-tau})}}{=}& \partial_{\tau^1}\partial_{\tau^{a}}^2F(0)\partial_{\tau^I}\partial_{\tau^a}\partial_{\tau^b}\partial_{\tau^{n-1}}F(0)
	+\partial_{\tau^1}\partial_{\tau^{a}}\partial_{\tau^b}\partial_{\tau^I}F(0)\partial_{\tau^b}^2\partial_{\tau^{n-1}}F(0)\\
&& +	\sum_{\begin{subarray}{c}0\leq J\leq I\\ 1\leq |J|\leq |I|-1\end{subarray}}\binom{I}{J}\partial_{\tau^1}\partial_{\tau^{a}}\partial_{\tau^J}\partial_{\tau^e}F(0)\eta^{ef}
	\partial_{\tau^{f}}\partial_{\tau^{I-J}}\partial_{\tau^{n-1}}\partial_{\tau^b}F(0)\\
&\stackrel{\mbox{by (\ref{eq-F122-tau})}}{=}&	\sum_{\begin{subarray}{c}0\leq J\leq I\\ 1\leq |J|\leq |I|-1\end{subarray}}\binom{I}{J}\partial_{\tau^1}\partial_{\tau^{a}}\partial_{\tau^J}\partial_{\tau^e}F(0)\eta^{ef}
	\partial_{\tau^{f}}\partial_{\tau^{I-J}}\partial_{\tau^{n-1}}\partial_{\tau^b}F(0).
\end{eqnarray*}
So (\ref{eq-recursion-primitive-abcc-even(2,2)-3}) follows.
\end{proof}

\begin{proposition}\label{proposition-recursion-primitive-aabb-even(2,2)}
Let $I=(i_{n+1},\dots,i_{2n+3})$ be a list  indicating the number of insertions of $\epsilon_{n+1},\dots,\epsilon_{2n+3}$. Suppose $n+1\leq a, b\leq 2n+3$ and $a\neq b$.
Then
\begin{eqnarray}\label{eq-recursion-primitive-aabb-even(2,2)}
&&(\frac{2|I|-4}{n-1}-2i_b)\partial_{\tau^{a}}^2\partial_{\tau^I}F(0)
	+(\frac{2|I|-4}{n-1}-2i_a)\partial_{\tau^{b}}^2\partial_{\tau^I}F(0)\nn\\
&=&\frac{1}{4}\sum_{\begin{subarray}{c}0\leq J\leq I\\ 2\leq |J|\leq |I|\end{subarray}}\binom{I}{J}\partial_{\tau^1}\partial_{\tau^{n-1}}\partial_{\tau^J}\partial_{\tau^e}F(0)\eta^{ef}
	\partial_{\tau^{f}}\partial_{\tau^{I-J}}\partial_{\tau^a}^2F(0) \nn\\
&&- \frac{1}{4}\sum_{\begin{subarray}{c}0\leq J\leq I\\ 1\leq |J|\leq |I|-1\end{subarray}}\binom{I}{J}\partial_{\tau^1}\partial_{\tau^{a}}\partial_{\tau^J}\partial_{\tau^e}F(0)\eta^{ef}
	\partial_{\tau^{f}}\partial_{\tau^{I-J}}\partial_{\tau^{n-1}}\partial_{\tau^a}F(0)\nn\\	
&&+\frac{1}{4}\sum_{\begin{subarray}{c}0\leq J\leq I\\ 2\leq |J|\leq |I|\end{subarray}}\binom{I}{J}\partial_{\tau^1}\partial_{\tau^{n-1}}\partial_{\tau^J}\partial_{\tau^e}F(0)\eta^{ef}
	\partial_{\tau^{f}}\partial_{\tau^{I-J}}\partial_{\tau^b}^2F(0) \nn\\
&&- \frac{1}{4}\sum_{\begin{subarray}{c}0\leq J\leq I\\ 1\leq |J|\leq |I|-1\end{subarray}}\binom{I}{J}\partial_{\tau^1}\partial_{\tau^{b}}\partial_{\tau^J}\partial_{\tau^e}F(0)\eta^{ef}
	\partial_{\tau^{f}}\partial_{\tau^{I-J}}\partial_{\tau^{n-1}}\partial_{\tau^b}F(0)\nn\\
&&- \sum_{\begin{subarray}{c}0\leq J\leq I\\ 2\leq |J|\leq |I|-2\end{subarray}}\binom{I}{J}\partial_{\tau^a}^2\partial_{\tau^J}\partial_{\tau^e}F(0)\eta^{ef}
	\partial_{\tau^{f}}\partial_{\tau^{I-J}}\partial_{\tau^b}^2F(0)\nn\\
&&+ \sum_{\begin{subarray}{c}0\leq J\leq I\\ 2\leq |J|\leq |I|-2\end{subarray}}\binom{I}{J}\partial_{\tau^a}\partial_{\tau^b}\partial_{\tau^J}\partial_{\tau^e}F(0)\eta^{ef}
	\partial_{\tau^{f}}\partial_{\tau^{I-J}}\partial_{\tau^a}\partial_{\tau^b}F(0).	
\end{eqnarray}	
\end{proposition}
\begin{proof}
We make use of the WDVV equation
\begin{equation}\label{eq-WDVV-aabb}
	(\partial_{\tau^a}\partial_{\tau^a}\partial_{\tau^e}F)\eta^{ef}(\partial_{\tau^f}\partial_{\tau^b}\partial_{\tau^b}F)=(\partial_{\tau^a}\partial_{\tau^b}\partial_{\tau^e}F)\eta^{ef}(\partial_{\tau^f}\partial_{\tau^a}\partial_{\tau^b}F).
\end{equation}
The coefficient of $\tau^I$ of the LHS of (\ref{eq-WDVV-aabb}), after multiplying $I!$, is
\begin{eqnarray*}
&& \sum_{0\leq J\leq I}\binom{I}{J}\partial_{\tau^a}^2\partial_{\tau^J}\partial_{\tau^e}F(0)\eta^{ef}
	\partial_{\tau^{f}}\partial_{\tau^{I-J}}\partial_{\tau^b}^2F(0)\\
&\stackrel{\mbox{by (\ref{eq-F122-tau})}}{=}& \partial_{\tau^a}^2\partial_{\tau^I}\partial_{\tau^e}F(0)\eta^{ef}
	\partial_{\tau^{f}}\partial_{\tau^b}^2F(0)
	+\partial_{\tau^a}^2\partial_{\tau^e}F(0)\eta^{ef}
	\partial_{\tau^{f}}\partial_{\tau^{I}}\partial_{\tau^b}^2F(0)\\
&&+ \sum_{k=n+1}^{2n+3}	i_k\partial_{\tau^a}^2\partial_{\tau^I}F(0)\partial_{\tau^{k}}^2\partial_{\tau^b}^2F(0)
	+ \sum_{k=n+1}^{2n+3}	i_k\partial_{\tau^a}^2\partial_{\tau^{k}}^2F(0)\partial_{\tau^I}\partial_{\tau^b}^2F(0)\\
&&+ \sum_{\begin{subarray}{c}0\leq J\leq I\\ 2\leq |J|\leq |I|-2\end{subarray}}\binom{I}{J}\partial_{\tau^a}^2\partial_{\tau^J}\partial_{\tau^e}F(0)\eta^{ef}
	\partial_{\tau^{f}}\partial_{\tau^{I-J}}\partial_{\tau^b}^2F(0)	\\
&\stackrel{\mbox{by (\ref{eq-F122-tau}), 
				(\ref{eq-etaInversePairing-even(2,2)}), 
				and (\ref{eq-4points-fanoIndex-even(2,2)-ab})}
				}{=}& -4\partial_{\tau^a}^2\partial_{\tau^I}\partial_{\tau^1}F(0)+
	\frac{1}{4}\partial_{\tau^a}^2\partial_{\tau^I}\partial_{\tau^n}F(0)
	+\frac{1}{4}\partial_{\tau^{n}}\partial_{\tau^I}\partial_{\tau^b}^2 F(0)
	-4\partial_{\tau^{1}}\partial_{\tau^I}\partial_{\tau^b}^2F(0)\\
&&+ |I|\partial_{\tau^a}^2\partial_{\tau^I}F(0)
	+ |I|\partial_{\tau^I}\partial_{\tau^b}^2F(0)\\
&&+ \sum_{\begin{subarray}{c}0\leq J\leq I\\ 2\leq |J|\leq |I|-2\end{subarray}}\binom{I}{J}\partial_{\tau^a}^2\partial_{\tau^J}\partial_{\tau^e}F(0)\eta^{ef}
	\partial_{\tau^{f}}\partial_{\tau^{I-J}}\partial_{\tau^b}^2F(0).	
\end{eqnarray*}
The coefficient of $\tau^I$ of the RHS of (\ref{eq-WDVV-aabb}), after multiplying $I!$, is
\begin{eqnarray*}
&& \sum_{0\leq J\leq I}\binom{I}{J}\partial_{\tau^a}\partial_{\tau^b}\partial_{\tau^J}\partial_{\tau^e}F(0)\eta^{ef}
	\partial_{\tau^{f}}\partial_{\tau^{I-J}}\partial_{\tau^a}\partial_{\tau^b}F(0)\\
&\begin{subarray}{c}\mbox{by (\ref{eq-F122-tau}) and} \\ \mbox{Theorem \ref{thm-monodromy-evenDim(2,2)}}\\ =\end{subarray}& i_a\partial_{\tau^a}\partial_{\tau^b}\partial_{\tau^{I-e_a}}\partial_{\tau^b}F(0)\partial_{\tau^{b}}\partial_{\tau^a}\partial_{\tau^a}\partial_{\tau^b}F(0)
	+i_b\partial_{\tau^a}\partial_{\tau^b}\partial_{\tau^{I-e_b}}\partial_{\tau^a}F(0)\partial_{\tau^{a}}\partial_{\tau^b}\partial_{\tau^a}\partial_{\tau^b}F(0)\\
&&+ i_a\partial_{\tau^a}\partial_{\tau^b}\partial_{\tau^a}\partial_{\tau^b}F(0)\partial_{\tau^{b}}\partial_{\tau^{I-e_a}}\partial_{\tau^a}\partial_{\tau^b}F(0)	
	+ i_b\partial_{\tau^a}\partial_{\tau^b}\partial_{\tau^b}\partial_{\tau^a}F(0)\partial_{\tau^{a}}\partial_{\tau^{I-e_b}}\partial_{\tau^a}\partial_{\tau^b}F(0)\\
&&+ \sum_{\begin{subarray}{c}0\leq J\leq I\\ 2\leq |J|\leq |I|-2\end{subarray}}\binom{I}{J}\partial_{\tau^a}\partial_{\tau^b}\partial_{\tau^J}\partial_{\tau^e}F(0)\eta^{ef}
	\partial_{\tau^{f}}\partial_{\tau^{I-J}}\partial_{\tau^a}\partial_{\tau^b}F(0)\\	
&\stackrel{\mbox{by  (\ref{eq-4points-fanoIndex-even(2,2)-ab})}}{=}& i_a\partial_{\tau^{I}}\partial_{\tau^b}^2F(0)
	+i_b\partial_{\tau^{I}}\partial_{\tau^a}^2F(0)
	+ i_a\partial_{\tau^{I}}\partial_{\tau^b}^2F(0)	
	+ i_b\partial_{\tau^{I}}\partial_{\tau^a}^2F(0)\\
&&+ \sum_{\begin{subarray}{c}0\leq J\leq I\\ 2\leq |J|\leq |I|-2\end{subarray}}\binom{I}{J}\partial_{\tau^a}\partial_{\tau^b}\partial_{\tau^J}\partial_{\tau^e}F(0)\eta^{ef}
	\partial_{\tau^{f}}\partial_{\tau^{I-J}}\partial_{\tau^a}\partial_{\tau^b}F(0)\\	
&=& 2 i_a\partial_{\tau^{I}}\partial_{\tau^b}^2F(0)
	+2i_b\partial_{\tau^{I}}\partial_{\tau^a}^2F(0)\\
&&	+ \sum_{\begin{subarray}{c}0\leq J\leq I\\ 2\leq |J|\leq |I|-2\end{subarray}}\binom{I}{J}\partial_{\tau^a}\partial_{\tau^b}\partial_{\tau^J}\partial_{\tau^e}F(0)\eta^{ef}
	\partial_{\tau^{f}}\partial_{\tau^{I-J}}\partial_{\tau^a}\partial_{\tau^b}F(0).
\end{eqnarray*}
So (\ref{eq-WDVV-aabb}) yields
\begin{eqnarray}\label{proposition-recursion-primitive-aabb-even(2,2)-1}
&&  -4\partial_{\tau^a}^2\partial_{\tau^I}\partial_{\tau^1}F(0)+
	\frac{1}{4}\partial_{\tau^a}^2\partial_{\tau^I}\partial_{\tau^n}F(0)
	+\frac{1}{4}\partial_{\tau^{n}}\partial_{\tau^I}\partial_{\tau^b}^2 F(0)
	-4\partial_{\tau^{1}}\partial_{\tau^I}\partial_{\tau^b}^2F(0)\nn\\
&&+ |I|\partial_{\tau^a}^2\partial_{\tau^I}F(0)
	+ |I|\partial_{\tau^I}\partial_{\tau^b}^2F(0)
	-2 i_a\partial_{\tau^{I}}\partial_{\tau^b}^2F(0)
	-2i_b\partial_{\tau^{I}}\partial_{\tau^a}^2F(0)\nn\\
&=&- \sum_{\begin{subarray}{c}0\leq J\leq I\\ 2\leq |J|\leq |I|-2\end{subarray}}\binom{I}{J}\partial_{\tau^a}^2\partial_{\tau^J}\partial_{\tau^e}F(0)\eta^{ef}
	\partial_{\tau^{f}}\partial_{\tau^{I-J}}\partial_{\tau^b}^2F(0)\nn\\
&&+ \sum_{\begin{subarray}{c}0\leq J\leq I\\ 2\leq |J|\leq |I|-2\end{subarray}}\binom{I}{J}\partial_{\tau^a}\partial_{\tau^b}\partial_{\tau^J}\partial_{\tau^e}F(0)\eta^{ef}
	\partial_{\tau^{f}}\partial_{\tau^{I-J}}\partial_{\tau^a}\partial_{\tau^b}F(0).
\end{eqnarray}
We make  manipulations on the LHS of (\ref{proposition-recursion-primitive-aabb-even(2,2)-1}). Since
\[
\frac{n}{2}\times (2+|I|)-(n-3+2+|I|)=(\frac{n}{2}-1)|I|+1,
\]
by (\ref{eq-tTotau}) and (\ref{eq-Div}) we have
\begin{equation}\label{proposition-recursion-primitive-aabb-even(2,2)-2}
\partial_{\tau^a}^2\partial_{\tau^I}\partial_{\tau^1}F(0)
=\partial_{t^a}^2\partial_{t^I}\partial_{t^1}F(0)
=\frac{(\frac{n}{2}-1)|I|+1}{n-1}\partial_{\tau^a}^2\partial_{\tau^I}F(0).
\end{equation}
So from (\ref{eq-recursion-primitive-abcc-even(2,2)-3}) we obtain
\begin{eqnarray}\label{proposition-recursion-primitive-aabb-even(2,2)-4}
&&  -4\partial_{\tau^a}^2\partial_{\tau^I}\partial_{\tau^1}F(0)+
	\frac{1}{4}\partial_{\tau^a}^2\partial_{\tau^I}\partial_{\tau^n}F(0)
	+\frac{1}{4}\partial_{\tau^{n}}\partial_{\tau^I}\partial_{\tau^b}^2 F(0)
	-4\partial_{\tau^{1}}\partial_{\tau^I}\partial_{\tau^b}^2F(0)\nn\\
&&+ |I|\partial_{\tau^a}^2\partial_{\tau^I}F(0)
	+ |I|\partial_{\tau^I}\partial_{\tau^b}^2F(0)
	-2 i_a\partial_{\tau^{I}}\partial_{\tau^b}^2F(0)
	-2i_b\partial_{\tau^{I}}\partial_{\tau^a}^2F(0)\nn\\
&=&	-4\partial_{\tau^a}^2\partial_{\tau^I}\partial_{\tau^1}F(0)+
	|I|\partial_{\tau^{a}}^2\partial_{\tau^I}F(0)\nn\\
&&-\frac{1}{4}\sum_{\begin{subarray}{c}0\leq J\leq I\\ 2\leq |J|\leq |I|\end{subarray}}\binom{I}{J}\partial_{\tau^1}\partial_{\tau^{n-1}}\partial_{\tau^J}\partial_{\tau^e}F(0)\eta^{ef}
	\partial_{\tau^{f}}\partial_{\tau^{I-J}}\partial_{\tau^a}^2F(0) \nn\\
&&+ \frac{1}{4}\sum_{\begin{subarray}{c}0\leq J\leq I\\ 1\leq |J|\leq |I|-1\end{subarray}}\binom{I}{J}\partial_{t^1}\partial_{\tau^{a}}\partial_{\tau^J}\partial_{\tau^e}F(0)\eta^{ef}
	\partial_{\tau^{f}}\partial_{\tau^{I-J}}\partial_{\tau^{n-1}}\partial_{\tau^a}F(0)\nn\\	
&&	+|I|\partial_{\tau^{b}}^2\partial_{\tau^I}F(0)
-\frac{1}{4}\sum_{\begin{subarray}{c}0\leq J\leq I\\ 2\leq |J|\leq |I|\end{subarray}}\binom{I}{J}\partial_{\tau^1}\partial_{\tau^{n-1}}\partial_{\tau^J}\partial_{\tau^e}F(0)\eta^{ef}
	\partial_{\tau^{f}}\partial_{\tau^{I-J}}\partial_{\tau^b}^2F(0) \nn\\
&&+ \frac{1}{4}\sum_{\begin{subarray}{c}0\leq J\leq I\\ 1\leq |J|\leq |I|-1\end{subarray}}\binom{I}{J}\partial_{\tau^1}\partial_{\tau^{b}}\partial_{\tau^J}\partial_{\tau^e}F(0)\eta^{ef}
	\partial_{\tau^{f}}\partial_{\tau^{I-J}}\partial_{\tau^{n-1}}\partial_{\tau^b}F(0)\nn\\
&&	-4\partial_{\tau^{1}}\partial_{\tau^I}\partial_{\tau^b}^2F(0)
	+ |I|\partial_{\tau^a}^2\partial_{\tau^I}F(0)
	+ |I|\partial_{\tau^I}\partial_{\tau^b}^2F(0)
	-2 i_a\partial_{\tau^{I}}\partial_{\tau^b}^2F(0)
	-2i_b\partial_{\tau^{I}}\partial_{\tau^a}^2F(0)\nn\\
&=&	2(|I|-i_b)\partial_{\tau^{a}}^2\partial_{\tau^I}F(0)-4\partial_{\tau^1}\partial_{\tau^{a}}^2\partial_{\tau^I}F(0)
	+2(|I|-i_a)\partial_{\tau^{b}}^2\partial_{\tau^I}F(0)-4\partial_{\tau^1}\partial_{\tau^{b}}^2\partial_{\tau^I}F(0)\nn\\
&&-\frac{1}{4}\sum_{\begin{subarray}{c}0\leq J\leq I\\ 2\leq |J|\leq |I|\end{subarray}}\binom{I}{J}\partial_{\tau^1}\partial_{\tau^{n-1}}\partial_{\tau^J}\partial_{\tau^e}F(0)\eta^{ef}
	\partial_{\tau^{f}}\partial_{\tau^{I-J}}\partial_{\tau^a}^2F(0) \nn\\
&&+ \frac{1}{4}\sum_{\begin{subarray}{c}0\leq J\leq I\\ 1\leq |J|\leq |I|-1\end{subarray}}\binom{I}{J}\partial_{\tau^1}\partial_{\tau^{a}}\partial_{\tau^J}\partial_{\tau^e}F(0)\eta^{ef}
	\partial_{\tau^{f}}\partial_{\tau^{I-J}}\partial_{\tau^{n-1}}\partial_{\tau^a}F(0)\nn\\	
&&-\frac{1}{4}\sum_{\begin{subarray}{c}0\leq J\leq I\\ 2\leq |J|\leq |I|\end{subarray}}\binom{I}{J}\partial_{\tau^1}\partial_{\tau^{n-1}}\partial_{\tau^J}\partial_{\tau^e}F(0)\eta^{ef}
	\partial_{\tau^{f}}\partial_{\tau^{I-J}}\partial_{\tau^b}^2F(0) \nn\\
&&+ \frac{1}{4}\sum_{\begin{subarray}{c}0\leq J\leq I\\ 1\leq |J|\leq |I|-1\end{subarray}}\binom{I}{J}\partial_{\tau^1}\partial_{\tau^{b}}\partial_{\tau^J}\partial_{\tau^e}F(0)\eta^{ef}
	\partial_{\tau^{f}}\partial_{\tau^{I-J}}\partial_{\tau^{n-1}}\partial_{\tau^b}F(0)\nn\\	
&\stackrel{\mbox{by (\ref{proposition-recursion-primitive-aabb-even(2,2)-2})}}{=}&(2|I|-2i_b-4\times\frac{(\frac{n}{2}-1)|I|+1}{n-1})\partial_{\tau^{a}}^2\partial_{\tau^I}F(0)\nn\\
&&+(2|I|-2i_a-4\times\frac{(\frac{n}{2}-1)|I|+1}{n-1})\partial_{\tau^{b}}^2\partial_{\tau^I}F(0)\nn\\
&&-\frac{1}{4}\sum_{\begin{subarray}{c}0\leq J\leq I\\ 2\leq |J|\leq |I|\end{subarray}}\binom{I}{J}\partial_{\tau^1}\partial_{\tau^{n-1}}\partial_{\tau^J}\partial_{\tau^e}F(0)\eta^{ef}
	\partial_{\tau^{f}}\partial_{\tau^{I-J}}\partial_{\tau^a}^2F(0) \nn\\
&&+ \frac{1}{4}\sum_{\begin{subarray}{c}0\leq J\leq I\\ 1\leq |J|\leq |I|-1\end{subarray}}\binom{I}{J}\partial_{\tau^1}\partial_{\tau^{a}}\partial_{\tau^J}\partial_{\tau^e}F(0)\eta^{ef}
	\partial_{\tau^{f}}\partial_{\tau^{I-J}}\partial_{\tau^{n-1}}\partial_{\tau^a}F(0)\nn\\	
&&-\frac{1}{4}\sum_{\begin{subarray}{c}0\leq J\leq I\\ 2\leq |J|\leq |I|\end{subarray}}\binom{I}{J}\partial_{\tau^1}\partial_{\tau^{n-1}}\partial_{\tau^J}\partial_{\tau^e}F(0)\eta^{ef}
	\partial_{\tau^{f}}\partial_{\tau^{I-J}}\partial_{\tau^b}^2F(0) \nn\\
&&+ \frac{1}{4}\sum_{\begin{subarray}{c}0\leq J\leq I\\ 1\leq |J|\leq |I|-1\end{subarray}}\binom{I}{J}\partial_{\tau^1}\partial_{\tau^{b}}\partial_{\tau^J}\partial_{\tau^e}F(0)\eta^{ef}
	\partial_{\tau^{f}}\partial_{\tau^{I-J}}\partial_{\tau^{n-1}}\partial_{\tau^b}F(0)\nn\\	
&=& (\frac{2|I|-4}{n-1}-2i_b)\partial_{\tau^{a}}^2\partial_{\tau^I}F(0)
	+(\frac{2|I|-4}{n-1}-2i_a)\partial_{\tau^{b}}^2\partial_{\tau^I}F(0)\nn\\
&&-\frac{1}{4}\sum_{\begin{subarray}{c}0\leq J\leq I\\ 2\leq |J|\leq |I|\end{subarray}}\binom{I}{J}\partial_{\tau^1}\partial_{\tau^{n-1}}\partial_{\tau^J}\partial_{\tau^e}F(0)\eta^{ef}
	\partial_{\tau^{f}}\partial_{\tau^{I-J}}\partial_{\tau^a}^2F(0) \nn\\
&&+ \frac{1}{4}\sum_{\begin{subarray}{c}0\leq J\leq I\\ 1\leq |J|\leq |I|-1\end{subarray}}\binom{I}{J}\partial_{\tau^1}\partial_{\tau^{a}}\partial_{\tau^J}\partial_{\tau^e}F(0)\eta^{ef}
	\partial_{\tau^{f}}\partial_{\tau^{I-J}}\partial_{\tau^{n-1}}\partial_{\tau^a}F(0)\nn\\	
&&-\frac{1}{4}\sum_{\begin{subarray}{c}0\leq J\leq I\\ 2\leq |J|\leq |I|\end{subarray}}\binom{I}{J}\partial_{\tau^1}\partial_{\tau^{n-1}}\partial_{\tau^J}\partial_{\tau^e}F(0)\eta^{ef}
	\partial_{\tau^{f}}\partial_{\tau^{I-J}}\partial_{\tau^b}^2F(0) \nn\\
&&+ \frac{1}{4}\sum_{\begin{subarray}{c}0\leq J\leq I\\ 1\leq |J|\leq |I|-1\end{subarray}}\binom{I}{J}\partial_{\tau^1}\partial_{\tau^{b}}\partial_{\tau^J}\partial_{\tau^e}F(0)\eta^{ef}
	\partial_{\tau^{f}}\partial_{\tau^{I-J}}\partial_{\tau^{n-1}}\partial_{\tau^b}F(0).	
\end{eqnarray}
Then (\ref{eq-recursion-primitive-aabb-even(2,2)}) follows from (\ref{proposition-recursion-primitive-aabb-even(2,2)-1}) and (\ref{proposition-recursion-primitive-aabb-even(2,2)-4}).
\end{proof}

\begin{proposition}\label{proposition-recursion-primitive-abcc-even(2,2)}
Let $I=(i_{n+1},\dots,i_{2n+3})$ be a list  indicating the number of insertions of $\epsilon_{n+1},\dots,\epsilon_{2n+3}$. Suppose $n+1\leq a, b,c\leq 2n+3$ and $a,b,c$ are pairwise distinct. Then
\begin{eqnarray}\label{eq-recursion-primitive-abcc-even(2,2)}
 &&(\frac{2|I|-4}{n-1}-2i_c)\partial_{\tau^{a}}\partial_{\tau^b}\partial_{\tau^I}F(0)\nn\\	
&=&\frac{1}{4}\sum_{\begin{subarray}{c}0\leq J\leq I\\ 2\leq |J|\leq |I|\end{subarray}}\binom{I}{J}\partial_{\tau^1}\partial_{\tau^{n-1}}\partial_{\tau^J}\partial_{\tau^e}F(0)\eta^{ef}
	\partial_{\tau^{f}}\partial_{\tau^{I-J}}\partial_{\tau^a}\partial_{\tau^b}F(0) \nn\\
&&-\frac{1}{4}	\sum_{\begin{subarray}{c}0\leq J\leq I\\ 1\leq |J|\leq |I|-1\end{subarray}}\binom{I}{J}\partial_{\tau^1}\partial_{\tau^{a}}\partial_{\tau^J}\partial_{\tau^e}F(0)\eta^{ef}
	\partial_{\tau^{f}}\partial_{\tau^{I-J}}\partial_{\tau^{n-1}}\partial_{\tau^b}F(0)\nn\\
&&-\sum_{\begin{subarray}{c}0\leq J\leq I\\ 2\leq |J|\leq |I|-2\end{subarray}}\binom{I}{J}\partial_{\tau^a}\partial_{\tau^b}\partial_{\tau^J}\partial_{\tau^e}F(0)\eta^{ef}
	\partial_{\tau^{f}}\partial_{\tau^{I-J}}\partial_{\tau^c}^2F(0)\nn\\
&& +\sum_{\begin{subarray}{c}0\leq J\leq I\\ 2\leq |J|\leq |I|-2\end{subarray}}\binom{I}{J}\partial_{\tau^a}\partial_{\tau^c}\partial_{\tau^J}\partial_{\tau^e}F(0)\eta^{ef}
	\partial_{\tau^{f}}\partial_{\tau^{I-J}}\partial_{\tau^b}\partial_{\tau^c}F(0).
\end{eqnarray}
\end{proposition}
\begin{proof}
We make use of the WDVV equation
\begin{equation}\label{eq-WDVV-abcc}
	(\partial_{\tau^a}\partial_{\tau^b}\partial_{\tau^e}F)\eta^{ef}(\partial_{\tau^f}\partial_{\tau^c}\partial_{\tau^c}F)=(\partial_{\tau^a}\partial_{\tau^c}\partial_{\tau^e}F)\eta^{ef}(\partial_{\tau^f}\partial_{\tau^b}\partial_{\tau^c}F).
\end{equation}
The coefficient of $\tau^I$ of the LHS of (\ref{eq-WDVV-abcc}), after multiplying $I!$, is
\begin{eqnarray*}
&& \sum_{0\leq J\leq I}\binom{I}{J}\partial_{\tau^a}\partial_{\tau^b}\partial_{\tau^J}\partial_{\tau^e}F(0)\eta^{ef}
	\partial_{\tau^{f}}\partial_{\tau^{I-J}}\partial_{\tau^c}^2F(0)\\
&\begin{subarray}{c}\mbox{by (\ref{eq-F122-tau}) and} \\ \mbox{Theorem \ref{thm-monodromy-evenDim(2,2)}}\\ =\end{subarray}& \partial_{\tau^a}\partial_{\tau^b}\partial_{\tau^I}\partial_{\tau^e}F(0)\eta^{ef}
	\partial_{\tau^{f}}\partial_{\tau^c}^2F(0)
	+\partial_{\tau^a}\partial_{\tau^b}\partial_{\tau^e}F(0)\eta^{ef}
	\partial_{\tau^{f}}\partial_{\tau^{I}}\partial_{\tau^c}^2F(0)\\
&&+ \sum_{k=n+1}^{2n+3}	i_k\partial_{\tau^a}\partial_{\tau^b}\partial_{\tau^I}F(0)\partial_{\tau^{k}}^2\partial_{\tau^c}^2F(0)
	+ i_a\partial_{\tau^a}^2\partial_{\tau^{b}}^2F(0)\partial_{\tau^{I-e_a}}\partial_{\tau^b}\partial_{\tau^c}^2F(0)\\
&&+  i_b\partial_{\tau^a}^2\partial_{\tau^{b}}^2F(0)\partial_{\tau^{I-e_b}}\partial_{\tau^a}\partial_{\tau^c}^2F(0)	\\
&&+ \sum_{\begin{subarray}{c}0\leq J\leq I\\ 2\leq |J|\leq |I|-2\end{subarray}}\binom{I}{J}\partial_{\tau^a}\partial_{\tau^b}\partial_{\tau^J}\partial_{\tau^e}F(0)\eta^{ef}
	\partial_{\tau^{f}}\partial_{\tau^{I-J}}\partial_{\tau^c}^2F(0)	\\	
&=& \partial_{t^a}\partial_{t^b}\partial_{t^I}\partial_{t^1}F(0)\frac{1}{4}\partial_{t^{n-1}}\partial_{t^c}^2F(0)+
	\partial_{t^a}\partial_{t^b}\partial_{t^I}\partial_{t^n}F(0)\frac{1}{4}\partial_{t^{0}}\partial_{t^c}^2F(0)\\
&&+ \sum_{j=n+1}^{2n+3}	i_j\partial_{t^a}\partial_{t^b}\partial_{t^I}F(0)\partial_{t^{j}}^2\partial_{t^c}^2F(0)
	+ i_a\partial_{t^a}^2\partial_{t^{b}}^2F(0)\partial_{t^{I-e_a}}\partial_{t^b}\partial_{t^c}^2F(0)\\
&&+  i_b\partial_{t^a}^2\partial_{t^{b}}^2F(0)\partial_{t^{I-e_b}}\partial_{t^a}\partial_{t^c}^2F(0)	\\
&& +\sum_{\begin{subarray}{c}0\leq J\leq I\\ 2\leq |J|\leq |I|-2\end{subarray}}\binom{I}{J}\partial_{t^a}\partial_{t^b}\partial_{t^J}\partial_{t^e}F(0)g^{ef}
	\partial_{t^{f}}\partial_{t^{I-J}}\partial_{t^c}^2F(0)\\
&\stackrel{\mbox{by (\ref{eq-F122-tau}), 
				(\ref{eq-etaInversePairing-even(2,2)}), 
				and (\ref{eq-4points-fanoIndex-even(2,2)-ab})}
				}{=}& - 4\partial_{\tau^a}\partial_{\tau^b}\partial_{\tau^I}\partial_{\tau^1}F(0)
	+\frac{1}{4}\partial_{\tau^a}\partial_{\tau^b}\partial_{\tau^I}\partial_{\tau^n}F(0)\\
&&+ |I|\partial_{\tau^a}\partial_{\tau^b}\partial_{\tau^I}F(0)
	+ i_a\partial_{\tau^{I-e_a}}\partial_{\tau^b}\partial_{\tau^c}^2F(0)
	+  i_b\partial_{\tau^{I-e_b}}\partial_{\tau^a}\partial_{\tau^c}^2F(0)\\
&& +\sum_{\begin{subarray}{c}0\leq J\leq I\\ 2\leq |J|\leq |I|-2\end{subarray}}\binom{I}{J}\partial_{\tau^a}\partial_{\tau^b}\partial_{\tau^J}\partial_{\tau^e}F(0)\eta^{ef}
	\partial_{\tau^{f}}\partial_{\tau^{I-J}}\partial_{\tau^c}^2F(0).	
\end{eqnarray*}
The coefficient of $\tau^I$ of the RHS of (\ref{eq-WDVV-abcc}), after multiplying $I!$, is
\begin{eqnarray*}
&& \sum_{0\leq J\leq I}\binom{I}{J}\partial_{\tau^a}\partial_{\tau^c}\partial_{\tau^J}\partial_{\tau^e}F(0)\eta^{ef}
	\partial_{\tau^{f}}\partial_{\tau^{I-J}}\partial_{\tau^b}\partial_{\tau^c}F(0)\\
&\begin{subarray}{c}\mbox{by (\ref{eq-F122-tau}) and} \\ \mbox{Theorem \ref{thm-monodromy-evenDim(2,2)}}\\ =\end{subarray}& i_b\partial_{\tau^a}\partial_{\tau^c}\partial_{\tau^{I-e_b}}\partial_{\tau^c}F(0)\partial_{\tau^{c}}\partial_{\tau^b}\partial_{\tau^b}\partial_{\tau^c}F(0)
	+i_c\partial_{\tau^a}\partial_{\tau^c}\partial_{\tau^{I-e_c}}\partial_{\tau^b}F(0)\partial_{\tau^{b}}\partial_{\tau^c}\partial_{\tau^b}\partial_{\tau^c}F(0)\\
&&+ i_a\partial_{\tau^a}\partial_{\tau^c}\partial_{\tau^a}\partial_{\tau^c}F(0)\partial_{\tau^{c}}\partial_{\tau^{I-e_a}}\partial_{\tau^b}\partial_{\tau^c}F(0)	
	+ i_c\partial_{\tau^a}\partial_{\tau^c}\partial_{\tau^c}\partial_{\tau^a}F(0)\partial_{\tau^{a}}\partial_{\tau^{I-e_c}}\partial_{\tau^b}\partial_{\tau^c}F(0)\\
&& +\sum_{\begin{subarray}{c}0\leq J\leq I\\ 2\leq |J|\leq |I|-2\end{subarray}}\binom{I}{J}\partial_{\tau^a}\partial_{\tau^c}\partial_{\tau^J}\partial_{\tau^e}F(0)\eta^{ef}
	\partial_{\tau^{f}}\partial_{\tau^{I-J}}\partial_{\tau^b}\partial_{\tau^c}F(0)	\\
&\stackrel{\mbox{by  (\ref{eq-4points-fanoIndex-even(2,2)-ab})}}{=}& i_b\partial_{\tau^a}\partial_{\tau^c}^2\partial_{\tau^{I-e_b}}F(0)
	+2i_c\partial_{\tau^a}\partial_{\tau^b}\partial_{\tau^{I}}F(0)
	+ i_a\partial_{\tau^b}\partial_{\tau^{c}}^2\partial_{\tau^{I-e_a}}F(0)\\
&& +\sum_{\begin{subarray}{c}0\leq J\leq I\\ 2\leq |J|\leq |I|-2\end{subarray}}\binom{I}{J}\partial_{\tau^a}\partial_{\tau^c}\partial_{\tau^J}\partial_{\tau^e}F(0)\eta^{ef}
	\partial_{\tau^{f}}\partial_{\tau^{I-J}}\partial_{t^b}\partial_{\tau^c}F(0).	
\end{eqnarray*}
So (\ref{eq-WDVV-abcc}) yields
\begin{eqnarray}\label{proposition-recursion-primitive-abcc-even(2,2)-1}
&&-4\partial_{\tau^a}\partial_{\tau^b}\partial_{\tau^I}\partial_{\tau^1}F(0)
	+\frac{1}{4}\partial_{\tau^a}\partial_{\tau^b}\partial_{\tau^I}\partial_{\tau^n}F(0)
	+ (|I|-2i_c)\partial_{\tau^a}\partial_{\tau^b}\partial_{\tau^I}F(0)\nn\\
&=&-\sum_{\begin{subarray}{c}0\leq J\leq I\\ 2\leq |J|\leq |I|-2\end{subarray}}\binom{I}{J}\partial_{\tau^a}\partial_{\tau^b}\partial_{\tau^J}\partial_{\tau^e}F(0)\eta^{ef}
	\partial_{\tau^{f}}\partial_{\tau^{I-J}}\partial_{\tau^c}^2F(0)\nn\\
&& +\sum_{\begin{subarray}{c}0\leq J\leq I\\ 2\leq |J|\leq |I|-2\end{subarray}}\binom{I}{J}\partial_{\tau^a}\partial_{\tau^c}\partial_{\tau^J}\partial_{\tau^e}F(0)\eta^{ef}
	\partial_{\tau^{f}}\partial_{\tau^{I-J}}\partial_{\tau^b}\partial_{\tau^c}F(0).
\end{eqnarray}
By  (\ref{eq-tTotau})  and (\ref{eq-Div}),
\begin{equation}\label{proposition-recursion-primitive-abcc-even(2,2)-2}
\partial_{\tau^a}\partial_{\tau^b}\partial_{\tau^I}\partial_{\tau^1}F(0)
=\partial_{t^a}\partial_{t^b}\partial_{t^I}\partial_{t^1}F(0)
=\frac{(\frac{n}{2}-1)|I|+1}{n-1}\partial_{\tau^a}\partial_{\tau^b}\partial_{\tau^I}F(0).
\end{equation}
So  (\ref{eq-recursion-primitive-abcc-even(2,2)-3}) we obtain
\begin{eqnarray}\label{proposition-recursion-primitive-abcc-even(2,2)-4}
&&-4\partial_{\tau^a}\partial_{\tau^b}\partial_{\tau^I}\partial_{\tau^1}F(0)
	+\frac{1}{4}\partial_{\tau^a}\partial_{\tau^b}\partial_{\tau^I}\partial_{\tau^n}F(0)
	+ (|I|-2i_c)\partial_{\tau^a}\partial_{\tau^b}\partial_{\tau^I}F(0)\nn\\
&=&-4\partial_{\tau^a}\partial_{\tau^b}\partial_{\tau^I}\partial_{\tau^1}F(0)
	+|I|\partial_{\tau^{a}}\partial_{\tau^b}\partial_{\tau^I}F(0)\nn\\
&&-\frac{1}{4}\sum_{\begin{subarray}{c}0\leq J\leq I\\ 2\leq |J|\leq |I|\end{subarray}}\binom{I}{J}\partial_{\tau^1}\partial_{\tau^{n-1}}\partial_{\tau^J}\partial_{\tau^e}F(0)\eta^{ef}
	\partial_{\tau^{f}}\partial_{\tau^{I-J}}\partial_{\tau^a}\partial_{\tau^b}F(0) \nn\\
&&+\frac{1}{4}	\sum_{\begin{subarray}{c}0\leq J\leq I\\ 1\leq |J|\leq |I|-1\end{subarray}}\binom{I}{J}\partial_{\tau^1}\partial_{\tau^{a}}\partial_{\tau^J}\partial_{\tau^e}F(0)\eta^{ef}
	\partial_{\tau^{f}}\partial_{\tau^{I-J}}\partial_{\tau^{n-1}}\partial_{\tau^b}F(0)\nn\\
&& + (|I|-2i_c)\partial_{\tau^a}\partial_{\tau^b}\partial_{\tau^I}F(0)\nn\\	
&\stackrel{\mbox{by (\ref{proposition-recursion-primitive-abcc-even(2,2)-2})}}{=}& (2|I|-2i_c-4\times\frac{(\frac{n}{2}-1)|I|+1}{n-1})\partial_{\tau^a}\partial_{\tau^b}\partial_{\tau^I}F(0)\nn\\	
&&-\frac{1}{4}\sum_{\begin{subarray}{c}0\leq J\leq I\\ 2\leq |J|\leq |I|\end{subarray}}\binom{I}{J}\partial_{\tau^1}\partial_{\tau^{n-1}}\partial_{\tau^J}\partial_{\tau^e}F(0)\eta^{ef}
	\partial_{\tau^{f}}\partial_{\tau^{I-J}}\partial_{\tau^a}\partial_{\tau^b}F(0) \nn\\
&&+\frac{1}{4}	\sum_{\begin{subarray}{c}0\leq J\leq I\\ 1\leq |J|\leq |I|-1\end{subarray}}\binom{I}{J}\partial_{\tau^1}\partial_{\tau^{a}}\partial_{\tau^J}\partial_{\tau^e}F(0)\eta^{ef}
	\partial_{\tau^{f}}\partial_{\tau^{I-J}}\partial_{\tau^{n-1}}\partial_{\tau^b}F(0)\nn\\
&=& (\frac{2|I|-4}{n-1}-2i_c)\partial_{\tau^{a}}\partial_{\tau^b}\partial_{\tau^I}F(0)\nn\\	
&&-\frac{1}{4}\sum_{\begin{subarray}{c}0\leq J\leq I\\ 2\leq |J|\leq |I|\end{subarray}}\binom{I}{J}\partial_{\tau^1}\partial_{\tau^{n-1}}\partial_{\tau^J}\partial_{\tau^e}F(0)\eta^{ef}
	\partial_{\tau^{f}}\partial_{\tau^{I-J}}\partial_{\tau^a}\partial_{\tau^b}F(0) \nn\\
&&+\frac{1}{4}	\sum_{\begin{subarray}{c}0\leq J\leq I\\ 1\leq |J|\leq |I|-1\end{subarray}}\binom{I}{J}\partial_{\tau^1}\partial_{\tau^{a}}\partial_{\tau^J}\partial_{\tau^e}F(0)\eta^{ef}
	\partial_{\tau^{f}}\partial_{\tau^{I-J}}\partial_{\tau^{n-1}}\partial_{\tau^b}F(0).
\end{eqnarray}
Then (\ref{eq-recursion-primitive-abcc-even(2,2)}) follows from (\ref{proposition-recursion-primitive-abcc-even(2,2)-1}) and (\ref{proposition-recursion-primitive-abcc-even(2,2)-4}).
\end{proof}

\subsection{A recursion with an unknown correlator}

\begin{theorem}\label{thm-reconstruction-even(2,2)}
With the knowledge of the 4-point invariants, all the invariants can be reconstructed from the WDVV, the deformation invariance, and the correlator
\begin{equation}\label{eq-specialLength(n+3)Invariant-even(2,2)}
	\langle \epsilon_{1},\dots,\epsilon_{n+3}\rangle_{0,n+3,\frac{n}{2}}.
\end{equation}
\end{theorem}
\begin{proof}
By Theorem \ref{thm-monodromy-evenDim(2,2)}, a  correlator with exactly one primitive insertion vanishes. By (\ref{eq-FCA}) a correlator of length $\geq 4$ having  an insertion $1$ vanishes. So using Lemma \ref{lem-recursion-EulerVecField-even(2,2)} and  \ref{lem-recursion-ambient-simplified-even(2,2)}, one can explicitly write a correlator of length $\geq 5$  with at least one ambient insertion and one primitive insertion as a combination of correlators of smaller lengths.

Suppose $I=(i_0,\dots,i_{2n+3})\in \mathbb{Z}_{\geq 0}^{2n+3}$ and  $|I|>4$. What we said above means that  if there exists $i_j>0$ for some $j\leq n$, then $\partial_{t^I}F(0)$ can be expressed as a combination of correlators with no larger length and  fewer ambient insertions. Repeating this process our task is  reduced to computing invariants with only primitive insertions. So we can assume that $I=(0,\dots,0,i_{n+1},\dots,i_{2n+3})$, and $|I|>4$. If there exists pairwise distinct $a,b,c$ in $\{n+1,\dots,2n+3\}$ such that 
\begin{equation}\label{eq-condition-recursion-even(2,2)}
	i_a,i_b>0, \mbox{and}\ |I|-4-(n-1)i_c\neq 0,
\end{equation}
then using (\ref{eq-recursion-primitive-abcc-even(2,2)}), $\partial_{t^I}F(0)$ can be  expressed as a combination of invariants with no larger length and  fewer primitive insertions. If there does not exist such $a,b,c$, then one of the following two cases happens:
\begin{enumerate}
 	\item[(i)] all $i_j$ but one vanish;
 	\item[(ii)] all $i_j>0$ and  are equal, and $|I|-4-(n-1)i_j=0$ for all $n+1\leq j\leq 2n+3$.
 \end{enumerate}  
 In the case  (i) one uses (\ref{eq-recursion-primitive-aabb-even(2,2)}) to reduce to the case that there exists at least two nonvanishing components in $I$. In the  case (ii), 
\[
i_{j}(n+3)-4=i_{j}(n-1)
\]
for all $n+1\leq j\leq 2n+3$, which implies $i_{n+1}=i_{n+2}=\dots=i_{2n+3}=1$, and thus corresponds to (\ref{eq-specialLength(n+3)Invariant-even(2,2)}).
\end{proof}

The following observation will be used in the  proof of Theorem \ref{thm-convergence}.
\begin{remark}\label{rem:recursion-boundOfIndex}
In the proof of Theorem \ref{thm-reconstruction-even(2,2)}, if the cases (i) and (ii) do not happen, we can take $a,b,c$  such that  $i_a$ and $i_b$ are the biggest two components in $I$ and  $i_c$ is the smallest one in $I$, and thus
\begin{equation}\label{eq-requirementOnic}
	i_c\leq \frac{|I|}{n+3}. 
\end{equation}
In fact, when $|I|>n+3$, 
\[
|I|-4>(n-1)\cdot \frac{|I|}{n+3},
\]
so the condition (\ref{eq-condition-recursion-even(2,2)}) is satisfied. 
When $|I|\leq n+3$, from Theorem \ref{thm-monodromy-evenDim(2,2)} it follows that the  correlator $\partial_{\tau^I}F(0)$ has $i_c=0$ unless it is zero or $I=(1,\dots,1)$.  Then the condition (\ref{eq-condition-recursion-even(2,2)}) is again satisfied. The requirement (\ref{eq-requirementOnic}) will be used in the proof of Theorem \ref{thm-convergence}.
\end{remark}

\begin{remark}{}
We call (\ref{eq-specialLength(n+3)Invariant-even(2,2)}) the \emph{special correlator} of $X$.
The proof of Theorem \ref{thm-reconstruction-even(2,2)} provides an effective algorithm to compute any correlators of $X$, where we regard the special correlator  as an indeterminate number. 
We implement the algorithm as a Macaulay2 package. For further information see Appendix \ref{sec:algorithm}.
\end{remark}

\subsection{Conjectures on the special correlator}\label{sec:conjecturesOnSpecialCorrelator}
We made some attempts to extract equations on the special correlator of $X$ from the WDVV equation and Theorem \ref{thm-monodromy-evenDim(2,2)}. All the relations that we found involving the special correlator turn out to be trivial ones. I guess that this is always true (unfortunately!). More precisely:
\begin{conjecture}\label{conj-specialCorrelator-free}
Set the special correlator to be an indeterminate $z$. Let $F(t_0,\dots,t_{2n+3};z)$ be the generating function of primary genus 0 Gromov-Witten invariants of $X$ computed by Algorithm \ref{algorithm-correlator-even(2,2)}  induced by the proof of Theorem \ref{thm-reconstruction-even(2,2)}. Then $F(t_0,\dots,t_{2n+3};z)$ satisfies (\ref{eq-WDVV}) and the conclusion of Theorem \ref{thm-monodromy-evenDim(2,2)}.
\end{conjecture}
This means that to compute the special correlator of $X$ one needs to introduce new tools. In the following of this section I try to guess the value of the special correlators, by making a comparison to the non-exceptional complete intersections of Fano index $n-1$, i.e. the cubic hypersurfaces and the odd dimensional complete intersections of two quadrics as we recalled at the beginning of Section \ref{sec:correlators-length4}. Recall that for the non-exceptional complete intersections, a correlator of odd length with only primitive insertions vanishes. So to make comparisons we need first to find an appropriate correlator of $X$, which has an even length, and express it in terms of the special correlator.
\begin{conjecture}\label{conj-unknownCorrelator-Even(2,2)-quadraticEquation}
For even $n$ dimensional complete intersections of two quadrics in $\mathbb{P}^{n+2}$, 
\begin{equation}\label{eq-unknownCorrelator-Even(2,2)-quadraticEquation-normalized}
\langle \epsilon_{1},\epsilon_1,\dots,\epsilon_{n+1},\epsilon_{n+1}\rangle_{0,2n+2,n-1}=2^{n-3}\big((\langle \epsilon_{1},\dots,\epsilon_{n+3}\rangle_{0,n+3,\frac{n}{2}})^2-\frac{1}{4}\big).
\end{equation}
Equivalently,
\begin{equation}\label{eq-unknownCorrelator-Even(2,2)-quadraticEquation}
\langle \varepsilon_{1},\varepsilon_1,\dots,\varepsilon_{n+1},\varepsilon_{n+1}\rangle_{0,2n+2,n-1}=2^{n-3}\big((\langle \varepsilon_{1},\dots,\varepsilon_{n+3}\rangle_{0,n+3,\frac{n}{2}})^2-\frac{(-1)^{\frac{n}{2}}}{4}\big).
\end{equation}
\end{conjecture}
We verify this conjecture of $4\leq n\leq 10$; see  Appendix \ref{sec:algorithm}. In \cite[Conjecture 10.26]{Hu15} we conjectured that for cubic hypersurfaces the correlators of the same form as the LHS of (\ref{eq-unknownCorrelator-Even(2,2)-quadraticEquation-normalized}) (equivalently, the LHS of (\ref{eq-unknownCorrelator-Even(2,2)})) vanish.  I guess that the same vanishing holds for $X$ when $n\equiv 0\mod 4$.
 Note that $\langle \varepsilon_{1},\dots,\varepsilon_{n+3}\rangle_{0,n+3,\frac{n}{2}}\in \mathbb{Q}$ and thus the  RHS of (\ref{eq-unknownCorrelator-Even(2,2)-quadraticEquation}) cannot vanish when $n\equiv 2\mod 4$.
 This is the reason that we limit the guess to the dimension $n\equiv 0\mod 4$, and leads to the following conjecture. 

Before stating the conjecture, we need to recall  Remark \ref{rem:choiceOfBasis} that there are choices of the basis $\varepsilon_1,\dots,\varepsilon_{n+3}$. All the previous statements are independent of such choices, while the following conjecture on the value of the special correlator does depend on.  

\begin{conjecture}\label{conj-unknownCorrelator-Even(2,2)}
For an even $n$ dimensional complete intersection of two quadrics in $\mathbb{P}^{n+2}$, let $\varepsilon_1,\dots,\varepsilon_{n+3}$ be the basis of $H^n_{\mathrm{prim}}(X)$ defined in Section \ref{sec:explictD-Lattice}. Then 
\begin{equation}\label{eq-unknownCorrelator-Even(2,2)}
	\langle \varepsilon_{1},\dots,\varepsilon_{n+3}\rangle_{0,n+3,\frac{n}{2}}=\frac{(-1)^{\frac{n}{2}}}{2},
\end{equation}
and 
\begin{equation*}
	\langle \varepsilon_{1},\varepsilon_1,\dots,\varepsilon_{n+1},\varepsilon_{n+1}\rangle_{0,2n+2,n-1}=
	\begin{cases}
	0,& \mbox{if}\ n\equiv 0 \mod 4,\\
	2^{n-4}, & \mbox{if}\ n\equiv 2 \mod 4.
	\end{cases}
\end{equation*}
\end{conjecture}
We will prove the $n=4$ case in Section \ref{sec:EnumerativeGeometry-Even(2,2)}.
At this stage we have no further evidence for Conjecture \ref{conj-unknownCorrelator-Even(2,2)}.
The values in (\ref{eq-unknownCorrelator-Even(2,2)}) are quite speculative, at least when $n\equiv 2 \mod 4$. For the reason for our choice of the sign in (\ref{eq-unknownCorrelator-Even(2,2)}), we refer the reader to Example \ref{example-f(6)}.

\section{Convergence of the generating function}\label{sec:convergence}
In this section, we still fix a smooth complete intersection $X$ of two quadrics in $\mathbb{P}^{n+2}$, where $n$ is even and $\geq 4$. 
Using Algorithm \ref{algorithm-correlator-even(2,2)} induced by the proof of Theorem \ref{thm-reconstruction-even(2,2)} we will show that the generating function $F$ has a positive convergence radius. So $F$ is an analytic function rather than only a formal series. The following theorem is a verification of \cite[Conjecture 1]{Zin14} for $X$.
\begin{theorem}\label{thm-convergence}
Let $\gamma_0,\dots,\gamma_{n},\gamma_{n+1},\dots,\gamma_{2n+3}$ be a basis of $H^*(X)$. Then there exists a constant $C>0$ such that 
\begin{equation}\label{eq-convergence-0}
	\langle \gamma_{i_1},\dots,\gamma_{i_k}\rangle_{0,k,\beta}\leq k! C^{k}
\end{equation}
for all $k\geq 0$, and $0\leq i_1,\dots,i_k\leq 2n+3$.
\end{theorem}

\begin{lemma}\label{lem-inequalityOfBinomialOfLists}
For any $M\leq |I|$, 
\begin{equation}\label{eq-inequalityOfBinomialOfLists}
	\sum_{\begin{subarray}{c}0\leq J\leq I\\ |J|\leq M\end{subarray}}\binom{I}{J}\leq \binom{|I|}{M}.
\end{equation}
\end{lemma}
\begin{proof}
This follows  from the enumerative meaning of both sides.
\end{proof}
\begin{lemma}\label{lem-inequality-1}
For $n\geq 4$,
\begin{equation}\label{eq-inequality-1}
	\sum_{k=2}^{n-2}\frac{n(n-1)}{k(k-1)(n-k)(n-k-1)}<4.
\end{equation}
\end{lemma}
\begin{proof}
We compute
\begin{eqnarray*}
&&\sum_{k=2}^{n-2}\frac{n(n-1)}{k(k-1)(n-k)(n-k-1)}\\
&=& \sum_{k=2}^{n-2}\frac{n-1}{(k-1)(n-k)(n-k-1)}+\sum_{k=2}^{n-2}\frac{n-1}{k(k-1)(n-k-1)}\\
&=& \sum_{k=2}^{n-2}\frac{1}{(n-k)(n-k-1)}+\sum_{k=2}^{n-2}\frac{1}{(k-1)(n-k-1)}\\
&&+\sum_{k=2}^{n-2}\frac{1}{(k-1)(n-k-1)}+\sum_{k=2}^{n-2}\frac{1}{k(k-1)}\\
&=& 2\sum_{k=2}^{n-2}\frac{1}{k(k-1)}+2\sum_{k=2}^{n-2}\frac{1}{(k-1)(n-k-1)}.
\end{eqnarray*}
Then we estimate the two sums separately. For the first sum,
\[
\sum_{k=2}^{n-2}\frac{1}{k(k-1)}=\sum_{k=2}^{n-2}(\frac{1}{k-1}-\frac{1}{k})=1-\frac{1}{n-2}.
\]
For the second sum we have
\begin{eqnarray*}
&&\sum_{k=2}^{n-2}\frac{1}{(k-1)(n-k-1)}\\
&=& 2\sum_{k=2}^{\frac{n}{2}-1}\frac{1}{(k-1)(n-k-1)}+\frac{1}{(\frac{n}{2}-1)^2}\\
&\leq & 2\sum_{k=2}^{\frac{n}{2}-1}\frac{1}{(k-1)\frac{n}{2}}+\frac{1}{(\frac{n}{2}-1)^2}\\
&\leq& \frac{4}{n}\big(1+\log(\frac{n}{2}-2)\big)+\frac{1}{(\frac{n}{2}-1)^2}.
\end{eqnarray*}
It follows that
\begin{eqnarray*}
\sum_{k=2}^{n-2}\frac{n(n-1)}{k(k-1)(n-k)(n-k-1)}
\leq  2-\frac{2}{n-2}+\frac{8}{n}\big(1+\log(\frac{n}{2}-2)\big)+\frac{2}{(\frac{n}{2}-1)^2}<4.
\end{eqnarray*}
\end{proof}

\begin{proof}[Proof of Theorem \ref{thm-convergence}]
The statement is independent of the choice of the basis $\gamma_0,\dots,\gamma_{2n+3}$. 
We take the basis $1,\tsfh_1,\dots,\tsfh_n,\epsilon_1,\dots,\epsilon_{n+3}$. Then (\ref{eq-convergence-0}) is equivalent to the existence of $C>0$ such that
\begin{equation}\label{eq-convergence-01}
	|\partial_{\tau^I}F(0)|\leq |I|! C^{|I|}
\end{equation}
for all $I\in \mathbb{Z}_{\geq 0}^{2n+4}$. 
Without loss of generality, we can assume that (\ref{eq-convergence-0}) holds for $k\leq K$, where $K$ is an arbitrary chosen natural number, and prove (\ref{eq-convergence-0}) inductively for all $k$. We note that the wanted statement  is equivalent to the existence of $C>0$ such that 
\begin{equation}\label{eq-convergence-02}
	|\partial_{\tau^I}F(0)|\leq (|I|-5)! C^{|I|-5}
\end{equation}
 for all $I\in \mathbb{Z}_{\geq 0}^{2n+4}$. By \cite[Theorem 1]{Zin14}, (\ref{eq-convergence-0}), equivalently (\ref{eq-convergence-02}), holds for correlators with only ambient classes; one can also find a simple proof of this fact in \cite[Remark D.13]{Hu15}.  In the following we show  (\ref{eq-convergence-02}) by induction on $|I|$ and the number of primitive insertions. Suppose (\ref{eq-convergence-02}) holds for $5\leq |I|\leq k$. 
By (\ref{eq-recursion-EulerVecField-even(2,2)}),
\begin{eqnarray*}
|\partial_{\tau^1}\partial_{\tau^I}F(0)|
&=&\big| \frac{\sum_{j=0}^n(j-1)i_j+(\frac{n}{2}-1)\sum_{j=n+1}^{2n+3}i_j+3-n}{n-1}\partial_{\tau^I}F(0)\nn\\
&&-12i_n\partial_{\tau^1}\partial_{\tau^{I-e_n}}F(0)\big|\\
&\leq & |I|\cdot  (|I|-5)! C^{|I|-5} +12|I|\cdot (|I|-5)! C^{|I|-5}.
\end{eqnarray*}
Thus replacing $C$ by a constant $C>65$ if necessary, we have
\begin{eqnarray*}
|\partial_{\tau^1}\partial_{\tau^I}F(0)|\leq (|I|-4)! C^{|I|-4}.
\end{eqnarray*}
In the following of the proof we use a temporary convention
\begin{equation*}
	k!=1\ \mbox{for}\ k<0.
\end{equation*}
Then by (\ref{eq-recursion-ambient-even(2,2)}) and (\ref{eq-etaInversePairing-even(2,2)}), for $2\leq i\leq n$ and $n+1\leq a,b\leq 2n+3$,
\begin{eqnarray*}
&&|\partial_{\tau^i}\partial_{\tau^I}\partial_{\tau^a}\partial_{\tau^b}F(0)|\\
&\leq &\Big|\sum_{\begin{subarray}{c}0\leq J\leq I\\ 1\leq |J|\leq |I|\end{subarray}}\binom{I}{J}\partial_{\tau^1}\partial_{\tau^{i-1}}\partial_{\tau^J}\partial_{\tau^e}F(0)\eta^{ef}
	\partial_{\tau^{f}}\partial_{\tau^{I-J}}\partial_{\tau^a}\partial_{\tau^b}F(0)\Big|\nn\\
&&+	\Big|\sum_{\begin{subarray}{c}0\leq J\leq I\\ 1\leq |J|\leq |I|-1\end{subarray}}\binom{I}{J}\partial_{\tau^1}\partial_{\tau^{a}}\partial_{\tau^J}\partial_{\tau^e}F(0)\eta^{ef}
	\partial_{\tau^{f}}\partial_{\tau^{I-J}}\partial_{\tau^{i-1}}\partial_{\tau^b}F(0)\Big|\\
&\leq & (2n+6)\times 2\times	4\sum_{\begin{subarray}{c}0\leq J\leq I\\ 1\leq |J|\leq |I|\end{subarray}}\binom{I}{J} (|J|-2)! C^{|J|-2} (|I|-|J|-2)! C^{|I|-|J|-2}\\
&=& 8(2n+6) (|I|-2)! C^{|I|-4} \sum_{\begin{subarray}{c}0\leq J\leq I\\ 1\leq |J|\leq |I|\end{subarray}}\binom{I}{J}\frac{|J|-2)!(|I|-|J|-2)!}{(|I|-2)!}.
\end{eqnarray*}
By Lemma \ref{lem-inequalityOfBinomialOfLists} and Lemma \ref{lem-inequality-1},
\begin{eqnarray*}
&&\sum_{\begin{subarray}{c}0\leq J\leq I\\ 1\leq |J|\leq |I|\end{subarray}}\binom{I}{J}\frac{|J|-2)!(|I|-|J|-2)!}{(|I|-2)!}\\
&\leq & \sum_{M=1}^{|I|}\binom{|I|}{M}\frac{(M-2)!(|I|-M-2)!}{(|I|-2)!}\\
&=& 2\times \frac{|I|}{|I|-2}+1+\sum_{M=2}^{|I|-2}\frac{|I|(|I|-1)}{M(M-1)(|I|-M)(|I|-M-1)}<9.
\end{eqnarray*}
So if $C^2>72(2n+6)$, we have
\begin{eqnarray*}
&&|\partial_{\tau^i}\partial_{\tau^I}\partial_{\tau^a}\partial_{\tau^b}F(0)|\\
&<& 72(2n+6) (|I|-2)! C^{|I|-4} <(|I|-2)! C^{|I|-2}.
\end{eqnarray*}

By (\ref{eq-recursion-primitive-abcc-even(2,2)}), 
\begin{eqnarray*}
 &&\big|(\frac{2|I|-4}{n-1}-2i_c)\partial_{\tau^{a}}\partial_{\tau^b}\partial_{\tau^I}F(0)\big|\nn\\	
&=&\Big|\frac{1}{4}\sum_{\begin{subarray}{c}0\leq J\leq I\\ 2\leq |J|\leq |I|\end{subarray}}\binom{I}{J}\partial_{\tau^1}\partial_{\tau^{n-1}}\partial_{\tau^J}\partial_{\tau^e}F(0)\eta^{ef}
	\partial_{\tau^{f}}\partial_{\tau^{I-J}}\partial_{\tau^a}\partial_{\tau^b}F(0) \nn\\
&&-\frac{1}{4}	\sum_{\begin{subarray}{c}0\leq J\leq I\\ 1\leq |J|\leq |I|-1\end{subarray}}\binom{I}{J}\partial_{\tau^1}\partial_{\tau^{a}}\partial_{\tau^J}\partial_{\tau^e}F(0)\eta^{ef}
	\partial_{\tau^{f}}\partial_{\tau^{I-J}}\partial_{\tau^{n-1}}\partial_{\tau^b}F(0)\nn\\
&&-\sum_{\begin{subarray}{c}0\leq J\leq I\\ 2\leq |J|\leq |I|-2\end{subarray}}\binom{I}{J}\partial_{\tau^a}\partial_{\tau^b}\partial_{\tau^J}\partial_{\tau^e}F(0)\eta^{ef}
	\partial_{\tau^{f}}\partial_{\tau^{I-J}}\partial_{\tau^c}^2F(0)\nn\\
&& +\sum_{\begin{subarray}{c}0\leq J\leq I\\ 2\leq |J|\leq |I|-2\end{subarray}}\binom{I}{J}\partial_{\tau^a}\partial_{\tau^c}\partial_{\tau^J}\partial_{\tau^e}F(0)\eta^{ef}
	\partial_{\tau^{f}}\partial_{\tau^{I-J}}\partial_{\tau^b}\partial_{\tau^c}F(0)\Big|\\
&\leq & 4\times (2n+6)\times 4 	\sum_{\begin{subarray}{c}0\leq J\leq I\\ 1\leq |J|\leq |I|\end{subarray}}\binom{I}{J} (|J|-2)!C^{|J|-2} (|I|-|J|-2)!C^{|I|-|J|-2}\\
&= & 16(2n+6) (|I|-2)! C^{|I|-4}	\sum_{\begin{subarray}{c}0\leq J\leq I\\ 1\leq |J|\leq |I|\end{subarray}}\binom{I}{J} \frac{(|J|-2)! (|I|-|J|-2)!}{(|I|-2)!}\\
&<& 144(2n+6) (|I|-2)! C^{|I|-4},
\end{eqnarray*}
By Remark \ref{rem:recursion-boundOfIndex}, we can  assume
\begin{equation*}
	i_c\leq \frac{|I|+2}{n+3}. 
\end{equation*}
Then 
\[
\frac{2|I|-4}{n-1}-2i_c\geq \frac{8(|I|-n-1)}{(n-1)(n+3)},
\]
and thus
\begin{eqnarray*}
 &&\big|\partial_{\tau^{a}}\partial_{\tau^b}\partial_{\tau^I}F(0)\big|\nn\\	
&<& 144(2n+6) (|I|-3)! C^{|I|-4}\cdot (n-1)(n+3) \frac{|I|-2}{8(|I|-n-1)}.
\end{eqnarray*}
When $|I|\geq 2n+2$, 
\[
\frac{|I|-2}{8(|I|-n-1)}<\frac{1}{4}.
\]
So when $|I|\geq 2n+2$ and 
\[
C>\big(36(2n+6)(n-1)(n+3)\big)^2,
\]
we have
\begin{equation}\label{eq-convergence-primitive-abcc}
	\big|\partial_{\tau^{a}}\partial_{\tau^b}\partial_{\tau^I}F(0)\big|<(|I|-3)! C^{|I|-\frac{7}{2}}.
 \end{equation}

By (\ref{eq-recursion-primitive-aabb-tau-simplified-even(2,2)}) and (\ref{eq-convergence-primitive-abcc}), when $i_a=|I|$,
\begin{eqnarray*}
&&(\frac{2|I|-4}{n-1})\partial_{\tau^{a}}^2\partial_{\tau^I}F(0)\nn\\
&\leq &\Big|(\frac{2|I|-4}{n-1}-2|I|)\partial_{\tau^{b}}^2\partial_{\tau^I}F(0)\Big|\\
&&+\Big|\frac{1}{4}\sum_{\begin{subarray}{c}0\leq J\leq I\\ 2\leq |J|\leq |I|\end{subarray}}
\binom{I}{J}\partial_{\tau^1}\partial_{\tau^{n-1}}\partial_{\tau^J}\partial_{\tau^e}F(0)\eta^{ef}
	\partial_{\tau^{f}}\partial_{\tau^{I-J}}\partial_{\tau^a}^2F(0) \nn\\
&&- \frac{1}{4}\sum_{\begin{subarray}{c}0\leq J\leq I\\ 1\leq |J|\leq |I|-1\end{subarray}}\binom{I}{J}\partial_{\tau^1}\partial_{\tau^{a}}\partial_{\tau^J}\partial_{\tau^e}F(0)\eta^{ef}
	\partial_{\tau^{f}}\partial_{\tau^{I-J}}\partial_{\tau^{n-1}}\partial_{\tau^a}F(0)\nn\\	
&&+\frac{1}{4}\sum_{\begin{subarray}{c}0\leq J\leq I\\ 2\leq |J|\leq |I|\end{subarray}}\binom{I}{J}\partial_{\tau^1}\partial_{\tau^{n-1}}\partial_{\tau^J}\partial_{\tau^e}F(0)\eta^{ef}
	\partial_{\tau^{f}}\partial_{\tau^{I-J}}\partial_{\tau^b}^2F(0) \nn\\
&&- \frac{1}{4}\sum_{\begin{subarray}{c}0\leq J\leq I\\ 1\leq |J|\leq |I|-1\end{subarray}}\binom{I}{J}\partial_{\tau^1}\partial_{\tau^{b}}\partial_{\tau^J}\partial_{\tau^e}F(0)\eta^{ef}
	\partial_{\tau^{f}}\partial_{\tau^{I-J}}\partial_{\tau^{n-1}}\partial_{\tau^b}F(0)\nn\\
&&- \sum_{\begin{subarray}{c}0\leq J\leq I\\ 2\leq |J|\leq |I|-2\end{subarray}}\binom{I}{J}\partial_{\tau^a}^2\partial_{\tau^J}\partial_{\tau^e}F(0)\eta^{ef}
	\partial_{\tau^{f}}\partial_{\tau^{I-J}}\partial_{\tau^b}^2F(0)\nn\\
&&+ \sum_{\begin{subarray}{c}0\leq J\leq I\\ 2\leq |J|\leq |I|-2\end{subarray}}\binom{I}{J}\partial_{\tau^a}\partial_{\tau^b}\partial_{\tau^J}\partial_{\tau^e}F(0)\eta^{ef}
	\partial_{\tau^{f}}\partial_{\tau^{I-J}}\partial_{\tau^a}\partial_{\tau^b}F(0)\Big|\\
&\leq & 2|I|\cdot(|I|-3)! C^{|I|-\frac{7}{2}}\\
&&	+ 6\times(2n+6)\times 4\times \sum_{\begin{subarray}{c}0\leq J\leq I\\ 1\leq |J|\leq |I|\end{subarray}}\binom{I}{J} (|J|-2)!C^{|J|-2} (|I|-|J|-2)!C^{|I|-|J|-2}\\
&\leq & 2|I|\cdot (|I|-3)! C^{|I|-\frac{7}{2}}\\
&&+ 24(2n+6) (|I|-2)! C^{|I|-4} \sum_{\begin{subarray}{c}0\leq J\leq I\\ 1\leq |J|\leq |I|\end{subarray}}\binom{I}{J} \frac{(|J|-2)! (|I|-|J|-2)!}{(|I|-2)!}\\
&\leq & (|I|-3)! C^{|I|-3}\big(\frac{2|I|}{C^{\frac{1}{2}}}+\frac{216(2n+6)(|I|-2)}{C}\big),
\end{eqnarray*}
so
\begin{eqnarray*}
&&\partial_{\tau^{a}}^2\partial_{\tau^I}F(0)
\leq (|I|-3)! C^{|I|-3}
\big(\frac{|I|(n-1)}{(|I|-2)C^{\frac{1}{2}}}+\frac{108(2n+6)(n-1)}{C}\big)<(|I|-3)! C^{|I|-3}
\end{eqnarray*}
when 
\begin{equation}\label{eq-convergence-estimateC}
\frac{|I|(n-1)}{(|I|-2)C^{\frac{1}{2}}}+\frac{108(2n+6)(n-1)}{C}<1.
\end{equation}
We choose $C$ such that (\ref{eq-convergence-02}) holds for $|I|<2n+1$, and such that 
\[
C>\max\{65,\sqrt{72(2n+6)},\big(36(2n+6)(n-1)(n+3)\big)^2
\}
\]
and (\ref{eq-convergence-estimateC}) holds. Then by Algorithm \ref{algorithm-correlator-even(2,2)} and the above estimates, (\ref{eq-convergence-02}) holds for all $I$.
\end{proof}

\begin{corollary}\label{cor-convergence}
There exists an open neighborhood of 0 in $\mathbb{C}^{2n+4}$ on which the generating function $F$ is an analytic function of $t^0,\dots,t^{2n+3}$.
\end{corollary}

\section{Semisimplicity of the Frobenius manifold}\label{sec:semisimplicity}
In this section, we still fix a smooth complete intersection $X$ of two quadrics in $\mathbb{P}^{n+2}$, where $n$ is even and $\geq 4$. By Corollary \ref{cor-convergence}, there exists an open neighborhood $U$ of $0\in \mathbb{C}^{2n+4}$, on which the generating function $F$ defines a Frobenius manifold $\mathcal{M}_X$. In this section we show
\begin{theorem}\label{thm-semisimplicity}
The Frobenius manifold $\mathcal{M}_X$ is (generically) tame semisimple.
\end{theorem}
By definition (\cite[\S 7.1]{Man99}), we need to show that at a general point of $U$, the multiplication by the Euler vector field $E$ has only simple eigenvalues. By \cite[Theorem 3.1]{Dub99}, this implies that $\mathcal{M}_X$ is (generically) semisimple. We work in the $\tau$-coordinates. We use Einstein's summation convention, where the range of the indices runs over $0,\dots,2n+3$.

\subsection{Quantum multiplication by the Euler vector field}
In $\tau$-coordinates, using (\ref{eq-EulerField-tau-Coordinates}), the big quantum multiplication by the Euler vector field $E$ is
\begin{eqnarray*}
 E\bqp \partial_{\tau^j}&=&\sum_{i=0}^{n}(1-i)\tau^i(\partial_{\tau^i}\partial_{\tau^j}\partial_{\tau^e}F)\eta^{ef}\partial_{\tau^f}+(4n-4) \tau^{n-1}\partial_{\tau^{j}}\\
&&	+(12n-12) \tau^n(\partial_{\tau^1}\partial_{\tau^j}\partial_{\tau^e}F)\eta^{ef}\partial_{\tau^f}
+\sum_{i=n+1}^{2n+3}(1-\frac{n}{2})\tau^{i}(\partial_{\tau^i}\partial_{\tau^j}\partial_{\tau^e}F)\eta^{ef}\partial_{\tau^f}\\
&&+(n-1)(\partial_{\tau^1}\partial_{\tau^j}\partial_{\tau^e}F)\eta^{ef}\partial_{\tau^f}.
\end{eqnarray*}
Denote by $\widetilde{\mathcal{E}}$ the matrix of the big quantum multiplication by $E$ in the basis $\partial_{\tau^0},\dots,\partial_{\tau^{2n+3}}$.
Then by (\ref{eq-etaInversePairing-even(2,2)}) we get
\begin{eqnarray*}
\widetilde{\mathcal{E}}_j^0&=&(4n-4) \tau^{n-1}\delta_{j}^0-4\sum_{i=0}^{n}(1-i)\tau^i \partial_{\tau^i}\partial_{\tau^j}\partial_{\tau^1}F
+\frac{1}{4}\sum_{i=0}^{n}(1-i)\tau^i \partial_{\tau^i}\partial_{\tau^j}\partial_{\tau^{n}}F\\
&&	-48(n-1) \tau^n\partial_{\tau^1}\partial_{\tau^j}\partial_{\tau^1}F
+3(n-1) \tau^n\partial_{\tau^1}\partial_{\tau^j}\partial_{\tau^n}F\\
&&-4\sum_{i=n+1}^{2n+3}(1-\frac{n}{2})\tau^{i}\partial_{\tau^i}\partial_{\tau^j}\partial_{\tau^1}F
+\frac{1}{4}\sum_{i=n+1}^{2n+3}(1-\frac{n}{2})\tau^{i}\partial_{\tau^i}\partial_{\tau^j}\partial_{\tau^{n}}F\\
&&-4(n-1)\partial_{\tau^1}\partial_{\tau^j}\partial_{\tau^1}F
+\frac{1}{4}(n-1)\partial_{\tau^1}\partial_{\tau^j}\partial_{\tau^n}F,
\end{eqnarray*}
\begin{eqnarray*}
\widetilde{\mathcal{E}}_j^1&=&(4n-4) \tau^{n-1}\delta_{j}^1-4\sum_{i=0}^{n}(1-i)\tau^i \partial_{\tau^i}\partial_{\tau^j}\partial_{\tau^0}F
+\frac{1}{4}\sum_{i=0}^{n}(1-i)\tau^i \partial_{\tau^i}\partial_{\tau^j}\partial_{\tau^{n-1}}F\\
&&	-48(n-1) \tau^n\partial_{\tau^1}\partial_{\tau^j}\partial_{\tau^0}F
+3(n-1) \tau^n\partial_{\tau^1}\partial_{\tau^j}\partial_{\tau^{n-1}}F\\
&&-4\sum_{i=n+1}^{2n+3}(1-\frac{n}{2})\tau^{i}\partial_{\tau^i}\partial_{\tau^j}\partial_{\tau^0}F
+\frac{1}{4}\sum_{i=n+1}^{2n+3}(1-\frac{n}{2})\tau^{i}\partial_{\tau^i}\partial_{\tau^j}\partial_{\tau^{n-1}}F\\
&&-4(n-1)\partial_{\tau^1}\partial_{\tau^j}\partial_{\tau^0}F
+\frac{1}{4}(n-1)\partial_{\tau^1}\partial_{\tau^j}\partial_{\tau^{n-1}}F,
\end{eqnarray*}
and for $2\leq k\leq n$,
\begin{eqnarray*}
\widetilde{\mathcal{E}}_j^k&=&(4n-4) \tau^{n-1}\delta_{j}^k
+\frac{1}{4}\sum_{i=0}^{n}(1-i)\tau^i \partial_{\tau^i}\partial_{\tau^j}\partial_{\tau^{n-k}}F\\
&&+3(n-1) \tau^n\partial_{\tau^1}\partial_{\tau^j}\partial_{\tau^{n-k}}F
+\frac{1}{4}\sum_{i=n+1}^{2n+3}(1-\frac{n}{2})\tau^{i}\partial_{\tau^i}\partial_{\tau^j}\partial_{\tau^{n-k}}F\\
&&+\frac{1}{4}(n-1)\partial_{\tau^1}\partial_{\tau^j}\partial_{\tau^{n-k}}F,
\end{eqnarray*}
and for $n+1\leq k\leq 2n+3$,
\begin{eqnarray*}
\widetilde{\mathcal{E}}_j^k&=&(4n-4) \tau^{n-1}\delta_{j}^k
+\sum_{i=0}^{n}(1-i)\tau^i \partial_{\tau^i}\partial_{\tau^j}\partial_{\tau^{k}}F\\
&&+12(n-1) \tau^n\partial_{\tau^1}\partial_{\tau^j}\partial_{\tau^{k}}F
+\sum_{i=n+1}^{2n+3}(1-\frac{n}{2})\tau^{i}\partial_{\tau^i}\partial_{\tau^j}\partial_{\tau^{k}}F\\
&&+(n-1)\partial_{\tau^1}\partial_{\tau^j}\partial_{\tau^{k}}F.
\end{eqnarray*}

\subsection{2nd order cutoff of the quantum multiplication by Euler vector field}
Taking the 2nd order  cutoff of the matrix $\widetilde{\mathcal{E}}$, and taking $\tau^i=0$ for  $0\leq i\leq n$, we denote the resulted matrix by $\mathcal{E}$, i.e.
\[
\mathcal{E}(\tau_{n+1},\dots,\tau_{2n+3})=\widetilde{\mathcal{E}}(0,\dots,0,\tau_{n+1},\dots,\tau_{2n+3})
+o(\tau^2).
\]
 Then by the formula of $\widetilde{\mathcal{E}}$ in the previous section, and the results of correlators of length 4 in Section \ref{sec:correlators-lengt-atMost4}, we get the following formula of $\mathcal{E}$.
\begin{eqnarray}\label{eq-cutoff-EulverVectorField-1}
\mathcal{E}_j^0&=&-4\sum_{i=n+1}^{2n+3}(1-\frac{n}{2})\tau^{i}\partial_{\tau^i}\partial_{\tau^j}\partial_{\tau^1}F
+\frac{1}{4}\sum_{i=n+1}^{2n+3}(1-\frac{n}{2})\tau^{i}\partial_{\tau^i}\partial_{\tau^j}\partial_{\tau^{n}}F\nn\\
&&-4(n-1)\partial_{\tau^1}\partial_{\tau^j}\partial_{\tau^1}F
+\frac{1}{4}(n-1)\partial_{\tau^1}\partial_{\tau^j}\partial_{\tau^n}F\nn\\
&\stackrel{\mod o(\mathbf{\tau}^2)}{=}& \begin{cases}
-4(1-\frac{n}{2})\cdot (-4)\sum_{i=n+1}^{2n+3}(\tau^i)^2
+\frac{1}{4}(1-\frac{n}{2})\cdot (-64)\sum_{i=n+1}^{2n+3}(\tau^i)^2&\\
-4(n-1)\cdot(-4)\cdot \frac{1}{2}\sum_{i=n+1}^{2n+3}(\tau^i)^2 &\\
+\frac{1}{4}(n-1)\cdot \big(2\times (-64)-12\times (-4)\big)\cdot \frac{1}{2}\sum_{i=n+1}^{2n+3}(\tau^i)^2
,& \mbox{if}\ j=n-1,\\
0 & \mbox{if otherwise},\\
\end{cases}\nn\\
&=& \begin{cases}
-2(n-1)\sum_{i=n+1}^{2n+3}(\tau^i)^2,& \mbox{if}\ j=n-1,\\
0 & \mbox{if otherwise}.
\end{cases}
\end{eqnarray}
\begin{eqnarray}\label{eq-cutoff-EulverVectorField-2}
\mathcal{E}_j^1&=&-4\sum_{i=n+1}^{2n+3}(1-\frac{n}{2})\tau^{i}\partial_{\tau^i}\partial_{\tau^j}\partial_{\tau^0}F
+\frac{1}{4}\sum_{i=n+1}^{2n+3}(1-\frac{n}{2})\tau^{i}\partial_{\tau^i}\partial_{\tau^j}\partial_{\tau^{n-1}}F\nn\\
&&-4(n-1)\partial_{\tau^1}\partial_{\tau^j}\partial_{\tau^0}F
+\frac{1}{4}(n-1)\partial_{\tau^1}\partial_{\tau^j}\partial_{\tau^{n-1}}F\nn\\
&\stackrel{\mod o(\mathbf{\tau}^2)}{=}&\begin{cases}
n-1, & \mbox{if}\ j=0,\\
\frac{1}{4}\sum_{i=n+1}^{2n+3}(1-\frac{n}{2})\tau^{i}\cdot(-4 \tau^i)
+\frac{1}{4}(n-1)\cdot(-4)\times \frac{1}{2}\sum_{i=n+1}^{2n+3}(\tau^i)^2,& \mbox{if}\ j=1,\\
-16(n-1)+\frac{1}{4}(n-1)\times 64,& \mbox{if}\ j=n-1,\\
\frac{1}{4}(1-\frac{n}{2})\cdot (-64)\sum_{i=n+1}^{2n+3}(\tau^i)^2 &\\
+\frac{1}{4}(n-1)\cdot \big(2\times (-64)-12\times (-4)\big)\cdot \frac{1}{2}\sum_{i=n+1}^{2n+3}(\tau^i)^2,& \mbox{if}\ j=n,\\
-4(1-\frac{n}{2})\tau^j+\frac{1}{4}(n-1)\cdot(-4 \tau^j),&\mbox{if}\ j\geq n+1
\end{cases}\nn\\
&=&\begin{cases}
n-1, & \mbox{if}\ j=0,\\
-\frac{1}{2}\sum_{i=n+1}^{2n+3}(\tau^{i})^2
,& \mbox{if}\ j=1,\\
0,& \mbox{if}\ j=n-1,\\
(-2n-6)\sum_{i=n+1}^{2n+3}(\tau^{i})^2,& \mbox{if}\ j=n,\\
(n-3)\tau^j,&\mbox{if}\ n+1\leq j\leq 2n+3.
\end{cases}
\end{eqnarray}
For $2\leq k\leq n-1$,
\begin{eqnarray}\label{eq-cutoff-EulverVectorField-3}
\mathcal{E}_j^k&=&\frac{1}{4}\sum_{i=n+1}^{2n+3}(1-\frac{n}{2})\tau^{i}\partial_{\tau^i}\partial_{\tau^j}\partial_{\tau^{n-k}}F
+\frac{1}{4}(n-1)\partial_{\tau^1}\partial_{\tau^j}\partial_{\tau^{n-k}}F\nn\\
&\stackrel{\mod o(\mathbf{\tau}^2)}{=}&\begin{cases}
n-1, & \mbox{if}\ j=k-1,\\
\frac{1}{4}\sum_{i=n+1}^{2n+3}(1-\frac{n}{2})\tau^{i}\cdot(-4 \tau^i)
+\frac{1}{4}(n-1)\cdot (-4)\cdot\frac{1}{2}\sum_{i=n+1}^{2n+3}(\tau^{i})^2,& \mbox{if}\ j=k,\\
16(n-1), & \mbox{if}\ (j,k)=(n,2),\\
0,\mbox{otherwise},
\end{cases}\nn\\
&=&\begin{cases}
n-1, & \mbox{if}\ j=k-1,\\
-\frac{1}{2}\sum_{i=n+1}^{2n+3}(\tau^{i})^2,& \mbox{if}\ j=k,\\
16(n-1), & \mbox{if}\ (j,k)=(n,2),\\
0,\mbox{otherwise}.
\end{cases}
\end{eqnarray}
and
\begin{eqnarray}\label{eq-cutoff-EulverVectorField-4}
\mathcal{E}_j^n&=&\frac{1}{4}\sum_{i=n+1}^{2n+3}(1-\frac{n}{2})\tau^{i}\partial_{\tau^i}\partial_{\tau^j}\partial_{\tau^{0}}F
+\frac{1}{4}(n-1)\partial_{\tau^1}\partial_{\tau^j}\partial_{\tau^{0}}F\nn\\
&\stackrel{\mod o(\mathbf{\tau}^2)}{=}&\begin{cases}
n-1,& \mbox{if}\ j=n-1,\\
0, & \mbox{if}\ 0\leq j\leq n-2\ \mbox{of}\ j=n,\\
\frac{2-n}{8}\tau^j,& \mbox{if}\ n+1\leq j\leq 2n+3.
\end{cases}
\end{eqnarray}
For $n+1\leq k\leq 2n+3$,
\begin{eqnarray}\label{eq-cutoff-EulverVectorField-5}
\mathcal{E}_j^k&=&\sum_{i=n+1}^{2n+3}(1-\frac{n}{2})\tau^{i}\partial_{\tau^i}\partial_{\tau^j}\partial_{\tau^{k}}F
+(n-1)\partial_{\tau^1}\partial_{\tau^j}\partial_{\tau^{k}}F\nn\\
&\stackrel{\mod o(\mathbf{\tau}^2)}{=}&\begin{cases}
(1-\frac{n}{2})\tau^k,& \mbox{if}\ j=0,\\
-4(n-1)\tau^k, & \mbox{if}\ j=n-1,\\
0,& \mbox{if}\ 1\leq j\leq n-2\ \mbox{or}\ j=n,\\
(2-n)\tau^j \tau^k+(n-1)\tau^j \tau^k, & \mbox{if}\ n+1\leq j\neq k\leq 2n+3,\\
(1-\frac{n}{2})\sum_{i=n+1}^{2n+3}(\tau^{i})^2+(n-1)\cdot\frac{1}{2}\sum_{i=n+1}^{2n+3}(\tau^{i})^2, & \mbox{if}\ j=k,
\end{cases}\nn\\
&=&\begin{cases}
(1-\frac{n}{2})\tau^k,& \mbox{if}\ j=0,\\
-4(n-1)\tau^k, & \mbox{if}\ j=n-1,\\
0,& \mbox{if}\ 0\leq j\leq n-2\ \mbox{or}\ j=n,\\
\tau^j \tau^k, & \mbox{if}\ n+1\leq j\neq k\leq 2n+3,\\
\frac{1}{2}\sum_{i=n+1}^{2n+3}(\tau^{i})^2, & \mbox{if}\ j=k.
\end{cases}
\end{eqnarray}

\subsection{The characteristic polynomial}

Let $\mathcal{P}_n(z,\tau^{n+1},\dots,\tau^{2n+3})$ be the characteristic polynomial of $\widetilde{\mathcal{E}}$ at $(0,\dots,0,\tau_{n+1},\tau_{n+2},\dots,\tau_{2n+3})$, and let $P_n(z,\tau^{n+1},\dots,\tau^{2n+3})$ be the characteristic polynomial of  $\mathcal{E}$ at $(\tau_{n+1},\tau_{n+2},\dots,\tau_{2n+3})$.
For brevity of expressions we define (recall (\ref{eq-invariantsOf-typeD-1}))
\[
s=2s_1=\sum_{i=n+1}^{2n+3}(\tau^{i})^2.
\]
Denote by $E_{i,j}$ the elementary matrix whose only nonzero entry is $1$ at the position $(i,j)$.
\begin{proposition}\label{prop-charPoly-EV-2ndOrder}
\begin{eqnarray}\label{eq-charPoly-EV-2ndOrder}
&&P_n(z,\tau^{n+1},\dots,\tau^{2n+3})\nn\\
&=&\bigg((n-1)^{n-1}\big(-\frac{(n-1)(n-2)^2 s^2}{4}+2(n-1)(n-4)sz
-4(n-5)z^2\big)\sum_{i=n+1}^{2n+3}\frac{(\tau^i)^2}{z-\frac{s}{2}+(\tau^i)^2}\nn\\
&&+4(n-1)^n s z-16(n-1)^{n-1}z^2+z^2(z+\frac{s}{2})^{n-1}
 -z^2(z+\frac{s}{2})^{n-1}\sum_{i=n+1}^{2n+3}\frac{(\tau^i)^2}{z-\frac{s}{2}+(\tau^i)^2})\bigg)\nn\\
&&\cdot \prod_{i=n+1}^{2n+3}\big(z-\frac{s}{2}+(\tau^i)^2\big).
\end{eqnarray}
\end{proposition}
\begin{proof}
By (\ref{eq-cutoff-EulverVectorField-1})-(\ref{eq-cutoff-EulverVectorField-5}), we have
{\footnotesize
\begin{equation*}
\mathcal{E}=\begin{pNiceArray}{ccccccc|ccc}
	0 & n-1 & 0 & 0 & \dots & 0 & 0 & (1-\frac{n}{2})\tau^{n+1} & \dots & (1-\frac{n}{2})\tau^{2n+3} \\
	0 & -\frac{s}{2} & n-1 & 0 & \dots & 0 & 0 & 0 & \dots & 0 \\
	0 & 0 & -\frac{s}{2} & n-1 & \dots & 0 & 0 &  0 & \dots & 0 \\
	\vdots & \vdots  & \vdots & \vdots& \vdots & \vdots& \vdots & 	 \vdots& \vdots& \vdots\\
	0 & 0 & \dots & 0 & \dots & n-1 & 0 & 	 0 & \dots& 0\\
	-2(n-1)s & 0 & \dots & \dots & \dots&  -\frac{s}{2} & n-1 & (4-4n)\tau^{n+1} & \dots & (4-4n)\tau^{2n+3}\\
	0 & (-2n-6)s & 16(n-1) & 0 & \dots& 0 & 0 &  0 & \dots & 0 \\
	\hline
	0 & (n-3)\tau^{n+1} & 0 & 0 & \dots & 0 & \frac{2-n}{8}\tau^{n+1} & \frac{s}{2} & \dots & \tau^{n+1}\tau^{2n+3}\\ 
	0 & (n-3)\tau^{n+2} & \dots & \dots & \dots & 0 & \frac{2-n}{8}\tau^{n+2} & \tau^{n+1}\tau^{n+2} & \dots & \tau^{n+2}\tau^{2n+2} \\
	\vdots & \vdots  & \vdots & \vdots& \vdots & \vdots& \vdots & 	\vdots& \vdots& \vdots\\
	0 & (n-3)\tau^{2n+3} & 0 & \dots & \dots & 0 & \frac{2-n}{8}\tau^{2n+3} & 	\tau^{n+1}\tau^{2n+3} & \dots & \frac{s}{2}
 	\end{pNiceArray}.
\end{equation*}
}
In the following we index the rows and columns by $0\leq i\leq 2n+3$. In the above we have blocked the matrix $\mathcal{E}$ by the first $n+1$ rows and first $n+1$ columns.
We perform several similarity transforms of $\mathcal{E}$.
Let 
\[
\mathcal{E}_1=\mathrm{Diag}(\underbrace{1,\dots,1}_{n+1},\tau^{n+1},\dots,\tau^{2n+3})\cdot
\mathcal{E}\cdot \mathrm{Diag}(\underbrace{1,\dots,1}_{n+1},\frac{1}{\tau^{n+1}},\dots,\frac{1}{\tau^{2n+3}}). 
\]
The effect is the $i$-th row $\times \tau^i$, the $i$-th column $\times \frac{1}{\tau^i}$, for $n+1\leq i\leq 2n+3$.
Then
{\footnotesize
\begin{equation*}
\mathcal{E}_1=\begin{pNiceArray}{ccccccc|ccc}
	0 & n-1 & 0 & 0 & \dots & 0 & 0 & 1-\frac{n}{2} & \dots & 1-\frac{n}{2} \\
	0 & -\frac{s}{2} & n-1 & 0 & \dots & 0 & 0 & 0 & \dots & 0 \\
	0 & 0 & -\frac{s}{2} & n-1 & \dots & 0 & 0 &  0 & \dots & 0 \\
	\vdots & \vdots  & \vdots & \vdots& \vdots & \vdots& \vdots & 	 \vdots& \vdots& \vdots\\
	0 & 0 & \dots & 0 & \dots & n-1 & 0 & 	 0 & \dots& 0\\
	-2(n-1)s & 0 & \dots & \dots & \dots&  -\frac{s}{2} & n-1 & 4-4n & \dots & 4-4n\\
	0 & (-2n-6)s & 16(n-1) & 0 & \dots& 0 & 0 &  0 & \dots & 0 \\
	\hline
	0 & (n-3)(\tau^{n+1})^2 & 0 & 0 & \dots & 0 & \frac{2-n}{8}(\tau^{n+1})^2 & \frac{s}{2} & \dots & (\tau^{n+1})^2\\ 
	0 & (n-3)(\tau^{n+2})^2 & \dots & \dots & \dots & 0 & \frac{2-n}{8}(\tau^{n+2})^2 & (\tau^{n+2})^2 & \dots & (\tau^{n+2})^2 \\
	\vdots & \vdots  & \vdots & \vdots& \vdots & \vdots& \vdots & 	\vdots& \vdots& \vdots\\
	0 & (n-3)(\tau^{2n+3})^2 & 0 & \dots & \dots & 0 & \frac{2-n}{8}(\tau^{2n+3})^2 & 	(\tau^{2n+3})^2 & \dots & \frac{s}{2}
 	\end{pNiceArray}.
\end{equation*}
}

Let 
\[
\mathcal{E}_2=(I_{2n+4}-\frac{8(n-3)}{n-2}E_{n,1})\cdot \mathcal{E}_1\cdot (I_{2n+4}+\frac{8(n-3)}{n-2}E_{n,1}).
\]
The effect is
\begin{eqnarray*}
&&\mbox{1st column}=>\mbox{1st column}+\frac{8(n-3)}{n-2}\times \mbox{$n$-th column},\\
&&\mbox{$n$-th row}=>\mbox{$n$-th row}-\frac{8(n-3)}{n-2}\times \mbox{$1$st row}.
\end{eqnarray*}
Then
{\footnotesize
\begin{equation*}
\mathcal{E}_2=\begin{pNiceArray}{ccccccc|ccc}
	0 & n-1 & 0 & 0 & \dots & 0 & 0 & 1-\frac{n}{2} & \dots & 1-\frac{n}{2} \\
	0 & -\frac{s}{2} & n-1 & 0 & \dots & 0 & 0 & 0 & \dots & 0 \\
	0 & 0 & -\frac{s}{2} & n-1 & \dots & 0 & 0 &  0 & \dots & 0 \\
	\vdots & \vdots  & \vdots & \vdots& \vdots & \vdots& \vdots & 	 \vdots& \vdots& \vdots\\
	0 & 0 & \dots & 0 & \dots & n-1 & 0 & 	 0 & \dots& 0\\
	-2(n-1)s & \frac{8(n-3)(n-1)}{n-2} & 0 & \dots & \dots&  -\frac{s}{2} & n-1 & 4-4n & \dots & 4-4n\\
	0 & \frac{-2n^2+2n}{n-2} s & \frac{8(n-1)^2}{n-2} & 0 & \dots& 0 & 0 &  0 & \dots & 0 \\
	\hline
	0 & 0 & 0 & 0 & \dots & 0 & \frac{2-n}{8}(\tau^{n+1})^2 & \frac{s}{2} & \dots & (\tau^{n+1})^2\\ 
	0 & 0 & \dots & \dots & \dots & 0 & \frac{2-n}{8}(\tau^{n+2})^2 & (\tau^{n+2})^2 & \dots & (\tau^{n+2})^2 \\
	\vdots & \vdots  & \vdots & \vdots& \vdots & \vdots& \vdots & 	\vdots& \vdots& \vdots\\
	0 & 0 & 0 & \dots & \dots & 0 & \frac{2-n}{8}(\tau^{2n+3})^2 & 	(\tau^{2n+3})^2 & \dots & \frac{s}{2}
 	\end{pNiceArray}.
\end{equation*}
}
For $i=n+1,\dots,2n+2$, we make the following transformation 
\[
U=> (I_{2n+4}+E_{2n+3,i})\cdot U\cdot (I_{2n+4}-E_{2n+3,i}). 
\]
The effect is, for  $i=n+1,\dots,2n+2$,
\begin{eqnarray*}
&&\mbox{$i$-th column}=>\mbox{$i$-th column}- \mbox{$(2n+3)$-th column},\\
&&\mbox{$(2n+3)$-th row}=> \mbox{$(2n+3)$-th row}+\mbox{$i$-th row}.
\end{eqnarray*}
Then let
\[
\mathcal{E}_4=(I_{2n+4}-\frac{8}{n-2}E_{n,2n+3})\cdot \mathcal{E}_4\cdot (I_{2n+4}+\frac{8}{n-2}E_{n,2n+3}).
\]
The effect is
\begin{eqnarray*}
&&\mbox{$(2n+3)$-th column}=>\mbox{$(2n+3)$-th column}+\frac{8}{n-2}\times \mbox{$n$-th column},\\
&&\mbox{$n$-th row}=>\mbox{$n$-th row}-\frac{8}{n-2}\times \mbox{$(2n+3)$-th row}.
\end{eqnarray*}
The matrices $\mathcal{E}_3$ and $\mathcal{E}_4$ are presented in the next page, and we block 
\[
\mathcal{E}_4=\begin{pmatrix}
A & B \\
C & D 
\end{pmatrix}
\]
as indicated there.  We compute the characteristic polynomial by the formula
\begin{eqnarray*}
&& P_n(z,\tau^{n+1},\dots,\tau^{2n+3})\\
&=& \det \begin{pmatrix}
z I-A & -B \\
-C & zI- D 
\end{pmatrix}
=\det \begin{pmatrix}
z I-A - B(zI-D)^{-1}C & 0 \\
-C & zI- D 
\end{pmatrix}.
\end{eqnarray*}

\begin{landscape}
\begin{eqnarray*}
\mathcal{E}_3=\begin{pNiceArray}{ccccccc|cccc}
	0 & n-1 & 0 & 0 & \dots & 0 & 0 & 0 & \dots & 0 & 1-\frac{n}{2} \\
	0 & -\frac{s}{2} & n-1 & 0 & \dots & 0 & 0 & 0 & \dots & 0 & 0 \\
	0 & 0 & -\frac{s}{2} & n-1 & \dots & 0 & 0 &  0 & \dots & 0  & 0 \\
	\vdots & \vdots  & \vdots & \vdots& \vdots & \vdots& \vdots & \vdots &	 \vdots& \vdots& \vdots\\
	0 & 0 & \dots & 0 & \dots & n-1 & 0 & 	 0 & \dots& 0 & 0\\
	-2(n-1)s & \frac{8(n-3)(n-1)}{n-2} & \dots & \dots & \dots&  -\frac{s}{2} & n-1 & 0 & \dots & 0 & 4-4n\\
	0 & \frac{-2n^2+2n}{n-2} s & \frac{8(n-1)^2}{n-2} & 0 & \dots& 0 & 0 &  0 & \dots & 0 & 0 \\
	\hline
	0 & 0 & 0 & 0 & \dots & 0 & \frac{2-n}{8}(\tau^{n+1})^2 & \frac{s}{2}-(\tau^{n+1})^2 & \dots & 0 &  (\tau^{n+1})^2\\ 
	0 & 0 & \dots & \dots & \dots & 0 & \frac{2-n}{8}(\tau^{n+2})^2 & 0 & \ddots & 0 & (\tau^{n+2})^2 \\
	\vdots & \vdots  & \vdots & \vdots& \vdots & \vdots& \vdots & 	\vdots&  & \vdots & \vdots\\
	0 & \dots & \dots & \dots & \dots & 0 & \frac{2-n}{8}(\tau^{2n+2})^2 & 0 & \dots & \frac{s}{2}-(\tau^{2n+2})^2 & (\tau^{2n+2})^2 \\
	0 & 0 & 0 & \dots & \dots & 0 & \frac{2-n}{8}s  & (\tau^{2n+3})^2-(\tau^{n+1})^2 & \dots & (\tau^{2n+3})^2-(\tau^{2n+2})^2 & \frac{3s}{2}-(\tau^{2n+3})^2 
 	\end{pNiceArray}.
\end{eqnarray*}

\begin{eqnarray*}
\mathcal{E}_4=\begin{pNiceArray}{ccccccc|cccc}
	0 & n-1 & 0 & 0 & \dots & 0 & 0 & 0 & \dots & 0 & 1-\frac{n}{2} \\
	0 & -\frac{s}{2} & n-1 & 0 & \dots & 0 & 0 & 0 & \dots & 0 & 0 \\
	0 & 0 & -\frac{s}{2} & n-1 & \dots & 0 & 0 &  0 & \dots & 0  & 0 \\
	\vdots & \vdots  & \vdots & \vdots& \vdots & \vdots& \vdots & \vdots &	 \vdots& \vdots& \vdots\\
	0 & 0 & \dots & 0 & \dots & n-1 & 0 & 	 0 & \dots& 0 & 0\\
	-2(n-1)s & \frac{8(n-3)(n-1)}{n-2} & \dots & \dots & \dots&  -\frac{s}{2} & n-1 & 0 & \dots & 0 & -\frac{4(n-1)(n-4)}{n-2}\\
	0 & \frac{-2n^2+2n}{n-2} s & \frac{8(n-1)^2}{n-2} & 0 & \dots& 0 & s &  -\frac{8\big((\tau^{2n+3})^2-(\tau^{n+1})^2\big)}{n-2} & \dots & -\frac{8\big((\tau^{2n+3})^2-(\tau^{2n+2})^2\big)}{n-2} & -\frac{8}{n-2}\big(\frac{s}{2}-(\tau^{2n+3})^2\big) \\
	\hline
	0 & 0 & 0 & 0 & \dots & 0 & \frac{2-n}{8}(\tau^{n+1})^2 & \frac{s}{2}-(\tau^{n+1})^2 & \dots & 0 &  0 \\ 
	0 & 0 & \dots & \dots & \dots & 0 & \frac{2-n}{8}(\tau^{n+2})^2 & 0 & \ddots & 0 & 0 \\
	\vdots & \vdots  & \vdots & \vdots& \vdots & \vdots& \vdots & 	\vdots&  & \vdots & \vdots\\
	0 & \dots & \dots & \dots & \dots & 0 & \frac{2-n}{8}(\tau^{2n+2})^2 & 0 & \dots & \frac{s}{2}-(\tau^{2n+2})^2 & 0 \\
	0 & 0 & 0 & \dots & \dots & 0 & \frac{2-n}{8}s  & (\tau^{2n+3})^2-(\tau^{n+1})^2 & \dots & (\tau^{2n+3})^2-(\tau^{2n+2})^2 & \frac{s}{2}-(\tau^{2n+3})^2 
 	\end{pNiceArray}.
\end{eqnarray*}
\end{landscape}
Since
\[
(zI-D)^{-1}=
\begin{pmatrix}
\frac{1}{z-\frac{s}{2}+(\tau^{n+1})^2} & 0 & \dots & 0 \\
0 & \ddots & & 0\\
0 & & \ddots & 0\\
\frac{(\tau^{2n+3})^2-(\tau^{n+1})^2}{(z-\frac{s}{2}+(\tau^{n+1})^2)(z-\frac{s}{2}+(\tau^{2n+3})^2)} & \dots & 
\frac{(\tau^{2n+3})^2-(\tau^{2n+2})^2}{(z-\frac{s}{2}+(\tau^{2n+2})^2)(z-\frac{s}{2}+(\tau^{2n+3})^2)} & \frac{1}{z-\frac{s}{2}+(\tau^{2n+3})^2}
\end{pmatrix}
\]
we get
\begin{gather*}
B(zI-D)^{-1}C=\\
\begin{pNiceArray}{c|c} 
\Block{6-1}<\large>{0_{{\scriptscriptstyle (n+1)\times n}}}  & 
{\scriptstyle
\frac{(n-2)^2}{16}\big(\sum_{i=n+1}^{2n+2}\frac{\big((\tau^{2n+3})^2-(\tau^i)^2\big)(\tau^i)^2}{\big(z-\frac{s}{2}+(\tau^{i})^2\big)\big(z-\frac{s}{2}+(\tau^{2n+3})^2\big)}+\frac{s}{z-\frac{s}{2}+(\tau^{2n+3})^2}\big) }\\
 \hspace*{1cm} & {\scriptstyle
0} \\ 
\hspace*{1cm} & \vdots \\
 \hspace*{1cm}  & {\scriptstyle
0} \\ 
 \hspace*{1cm}  & \vspace{1cm} {\scriptstyle
 \frac{(n-2)(n-4)}{2}\big(\sum_{i=n+1}^{2n+2}\frac{\big((\tau^{2n+3})^2-(\tau^i)^2\big)(\tau^i)^2}{\big(z-\frac{s}{2}+(\tau^{i})^2\big)\big(z-\frac{s}{2}+(\tau^{2n+3})^2\big)}+\frac{s}{z-\frac{s}{2}+(\tau^{2n+3})^2}\big) } \\ 
\hspace*{1cm}  &  {\scriptstyle
\sum_{i=n+1}^{2n+2}\frac{\big((\tau^{2n+3})^2-(\tau^i)^2\big)(\tau^i)^2}{z-\frac{s}{2}+(\tau^{i})^2}+\big(\frac{s}{2}-(\tau^{2n+3})^2\big) \big(\sum_{i=n+1}^{2n+2}\frac{\big((\tau^{2n+3})^2-(\tau^i)^2\big)(\tau^i)^2}{(z-\frac{s}{2}+(\tau^{i})^2)(z-\frac{s}{2}+(\tau^{2n+3})^2)}+\frac{s}{z-\frac{s}{2}+(\tau^{2n+3})^2}\big) }
\end{pNiceArray},
\end{gather*}
where $0_{(n+1)\times n}$ stands for a 0-matrix of size $(n+1)\times n$.
Using 
\[
\sum_{i=n+1}^{2n+2}\frac{\big((\tau^{2n+3})^2-(\tau^i)^2\big)(\tau^i)^2}{\big(z-\frac{s}{2}+(\tau^{i})^2\big)\big(z-\frac{s}{2}+(\tau^{2n+3})^2\big)}+\frac{s}{z-\frac{s}{2}+(\tau^{2n+3})^2}
=\sum_{i=n+1}^{2n+3}\frac{(\tau^i)^2}{z-\frac{s}{2}+(\tau^i)^2},
\]
and
\begin{eqnarray*}
&&\sum_{i=n+1}^{2n+2}\frac{\big((\tau^{2n+3})^2-(\tau^i)^2\big)(\tau^i)^2}{z-\frac{s}{2}+(\tau^{i})^2}\\
&&+\big(\frac{s}{2}-(\tau^{2n+3})^2\big) \big(\sum_{i=n+1}^{2n+2}\frac{\big((\tau^{2n+3})^2-(\tau^i)^2\big)(\tau^i)^2}{(z-\frac{s}{2}+(\tau^{i})^2)(z-\frac{s}{2}+(\tau^{2n+3})^2)}+\frac{s}{z-\frac{s}{2}+(\tau^{2n+3})^2}\big)\\
&=& z\sum_{i=n+1}^{2n+3}\frac{(\tau^i)^2}{z-\frac{s}{2}+(\tau^i)^2} -s,
\end{eqnarray*}
we get
\begin{eqnarray*}
B(zI-D)^{-1}C=\begin{pNiceArray}{c|c} 
\Block{6-1}<\large>{0_{{\scriptscriptstyle
(n+1)\times n}}}  & 
\frac{(n-2)^2}{16}\sum_{i=n+1}^{2n+3}\frac{(\tau^i)^2}{z-\frac{s}{2}+(\tau^i)^2} \\
 \hspace*{1cm} & 0 \\ 
\hspace*{1cm} & \vdots \\
 \hspace*{1cm}  & 0 \\ 
 \hspace*{1cm}  & \frac{(n-2)(n-4)}{2}\sum_{i=n+1}^{2n+3}\frac{(\tau^i)^2}{z-\frac{s}{2}+(\tau^i)^2}  \\ 
\hspace*{1cm}  &  z\sum_{i=n+1}^{2n+3}\frac{(\tau^i)^2}{z-\frac{s}{2}+(\tau^i)^2} -s 
\end{pNiceArray},
\end{eqnarray*}
and thus
\begin{gather*}
zI-A-B(zI-D)^{-1}C=\\
\begin{pNiceArray}{ccccccc}
	z & 1-n & 0 & 0 & \dots & 0 & -\frac{(n-2)^2}{16}\sum_{i=n+1}^{2n+3}\frac{(\tau^i)^2}{z-\frac{s}{2}+(\tau^i)^2}  \\
	0 & z+\frac{s}{2} & 1-n & 0 & \dots & 0 & 0 \\
	0 & 0 & z+\frac{s}{2} & 1-n & \dots & 0 & 0  \\
	\vdots & \vdots  & \vdots & \vdots& \vdots & \vdots& \vdots \\
	0 & 0 & \dots & 0 & \dots & 1-n & 0 \\
	2(n-1)s & -\frac{8(n-3)(n-1)}{n-2} & \dots & \dots & \dots&  z+\frac{s}{2} & 1-n-\frac{(n-2)(n-4)}{2}\sum_{i=n+1}^{2n+3}\frac{(\tau^i)^2}{z-\frac{s}{2}+(\tau^i)^2} \\
	0 & \frac{2n^2-2n}{n-2} s & -\frac{8(n-1)^2}{n-2} & 0 & \dots& 0 & z\big(1-\sum_{i=n+1}^{2n+3}\frac{(\tau^i)^2}{z-\frac{s}{2}+(\tau^i)^2}\big) 
 	\end{pNiceArray}.
\end{gather*}
We compute the determinant of $zI-A-B(zI-D)^{-1}C$ by expanding the determinant according to the last two columns and then the first column.
So the characteristic polynomial is
\begin{eqnarray*}
&&\det \big(zI-A-B(zI-D)^{-1}C\big)\\
&=&\Big(\frac{(n-2)^2}{16}\sum_{i=n+1}^{2n+3}\frac{(\tau^i)^2}{z-\frac{s}{2}+(\tau^i)^2}
\cdot(2n-2)s\\
&& +z\big(1-n-\frac{(n-1)(n-4)}{2}\sum_{i=n+1}^{2n+3}\frac{(\tau^i)^2}{z-\frac{s}{2}+(\tau^i)^2}
\big)\Big)\cdot (1-n)\\
&&\cdot \big((-1)^{n-3}\frac{2n^2-2n}{n-2}s\cdot  (1-n)^{n-3}
+(-1)^{n-2}(-\frac{8(n-1)^2}{n-2})\cdot (1-n)^{n-4}(z+\frac{s}{2})\big)\\
&&+z(1-\sum_{i=n+1}^{2n+3}\frac{(\tau^i)^2}{z-\frac{s}{2}+(\tau^i)^2})\big(
(-1)^{n-1}(2n-2)s\cdot (1-n)^{n-1}\\
&&+(-1)^n(-\frac{8(n-3)(n-1)}{n-2})\cdot (1-n)^{n-2}z+z(z+\frac{s}{2})^{n-1}
\big).
\end{eqnarray*}
We simplify its expression as follows, where the change in each step is indicated in \textcolor{blue}{blue}. Recall that $n$ is even.
\begin{eqnarray*}
&&\det \big(zI-A-B(zI-D)^{-1}C\big)\\
&=&\Big(-\frac{(n-2)^2 \textcolor{blue}{s}}{\textcolor{blue}{8}}\sum_{i=n+1}^{2n+3}\frac{(\tau^i)^2}{z-\frac{s}{2}+(\tau^i)^2} 
+z\big(1+\frac{n-4}{2}\sum_{i=n+1}^{2n+3}\frac{(\tau^i)^2}{z-\frac{s}{2}+(\tau^i)^2}
\big)\Big)\cdot \textcolor{blue}{(1-n)^2}\\
&&\cdot \big(\frac{2n(n-1)^{\textcolor{blue}{n-2}}}{n-2}s
-\frac{8(n-1)^{\textcolor{blue}{n-2}}}{n-2})(z+\frac{s}{2})\big)\\
&&+z(1-\sum_{i=n+1}^{2n+3}\frac{(\tau^i)^2}{z-\frac{s}{2}+(\tau^i)^2})\big(
2(n-1)^{\textcolor{blue}{n}} s
-\frac{8(n-3)(n-1)^{\textcolor{blue}{n-1}}}{n-2}z+z(z+\frac{s}{2})^{n-1}
\big)\\
&=&(n-1)^2\Big(-\frac{(n-2)^2 s}{8}\sum_{i=n+1}^{2n+3}\frac{(\tau^i)^2}{z-\frac{s}{2}+(\tau^i)^2}
+z\big(1+\frac{n-4}{2}\sum_{i=n+1}^{2n+3}\frac{(\tau^i)^2}{z-\frac{s}{2}+(\tau^i)^2}
\big)\Big)\\
&&\cdot \big(\textcolor{blue}{2(n-1)^{n-2}s}-\frac{8(n-1)^{n-2}}{n-2}z\big)\\
&&+z(1-\sum_{i=n+1}^{2n+3}\frac{(\tau^i)^2}{z-\frac{s}{2}+(\tau^i)^2})\big(
2(n-1)^n s
-\frac{8(n-3)(n-1)^{n-1}}{n-2}z+z(z+\frac{s}{2})^{n-1}
\big)\\
&=&(n-1)^{\textcolor{blue}{n}}\Big(-\frac{(n-2)^2 s}{8}\sum_{i=n+1}^{2n+3}\frac{(\tau^i)^2}{z-\frac{s}{2}+(\tau^i)^2}
+z\big(1+\frac{n-4}{2}\sum_{i=n+1}^{2n+3}\frac{(\tau^i)^2}{z-\frac{s}{2}+(\tau^i)^2}
\big)\Big)\\
&&\cdot \big(2s-\frac{8}{n-2}z\big)\\
&&+z(1-\sum_{i=n+1}^{2n+3}\frac{(\tau^i)^2}{z-\frac{s}{2}+(\tau^i)^2})\big(
2(n-1)^n s
-\frac{\textcolor{blue}{8(n-1)}(n-1)^{n-1}}{n-2}z
\big)\\
&& +z(1-\sum_{i=n+1}^{2n+3}\frac{(\tau^i)^2}{z-\frac{s}{2}+(\tau^i)^2})\big(
\frac{\textcolor{blue}{16}(n-1)^{n-1}}{n-2}z+z(z+\frac{s}{2})^{n-1}
\big)\\
&=&(n-1)^n\Big(-\frac{(n-2)^2 s}{8}\sum_{i=n+1}^{2n+3}\frac{(\tau^i)^2}{z-\frac{s}{2}+(\tau^i)^2}
+z\big(1+\frac{n-4}{2}\sum_{i=n+1}^{2n+3}\frac{(\tau^i)^2}{z-\frac{s}{2}+(\tau^i)^2}
\big)\Big)\\
&&\cdot \big(2s-\frac{8}{n-2}z\big)\\
&&+\textcolor{blue}{(n-1)^n} z(1-\sum_{i=n+1}^{2n+3}\frac{(\tau^i)^2}{z-\frac{s}{2}+(\tau^i)^2})\big(
2 s -\frac{8}{n-2}z\big)\\
&& +z^2(1-\sum_{i=n+1}^{2n+3}\frac{(\tau^i)^2}{z-\frac{s}{2}+(\tau^i)^2})\big(
\frac{16(n-1)^{n-1}}{n-2}+(z+\frac{s}{2})^{n-1}
\big)\\
&=&(n-1)^n\big(2s-\frac{8}{n-2}z\big)\Big(-\frac{(n-2)^2 s}{8}\sum_{i=n+1}^{2n+3}\frac{(\tau^i)^2}{z-\frac{s}{2}+(\tau^i)^2}
+\textcolor{blue}{2z+\frac{n-6}{2}z}\sum_{i=n+1}^{2n+3}\frac{(\tau^i)^2}{z-\frac{s}{2}+(\tau^i)^2}\Big)\\
&& +z^2(1-\sum_{i=n+1}^{2n+3}\frac{(\tau^i)^2}{z-\frac{s}{2}+(\tau^i)^2})\big(
\frac{16(n-1)^{n-1}}{n-2}+(z+\frac{s}{2})^{n-1}
\big)\\
&=&(n-1)^n\big(2s-\frac{8}{n-2}z\big)\textcolor{blue}{\Big(}-\frac{(n-2)^2 s}{8}
+\frac{n-6}{2}z\textcolor{blue}{\Big)}\sum_{i=n+1}^{2n+3}\frac{(\tau^i)^2}{z-\frac{s}{2}+(\tau^i)^2}\\
&&+(n-1)^n\big(2s-\frac{8}{n-2}z\big)2z\\
&& +z^2(1-\sum_{i=n+1}^{2n+3}\frac{(\tau^i)^2}{z-\frac{s}{2}+(\tau^i)^2})\big(
\frac{16(n-1)^{n-1}}{n-2}+(z+\frac{s}{2})^{n-1}
\big)\\
&=&(n-1)^n\big(s-\frac{\textcolor{blue}{4}}{n-2}z\big)\Big(-\frac{(n-2)^2 s}{\textcolor{blue}{4}}
+(n-6)z\Big)\sum_{i=n+1}^{2n+3}\frac{(\tau^i)^2}{z-\frac{s}{2}+(\tau^i)^2}\\
&&+\textcolor{blue}{4}(n-1)^n\big(s-\frac{4}{n-2}z\big)z+\textcolor{blue}{z^2\big(
\frac{16(n-1)^{n-1}}{n-2}+(z+\frac{s}{2})^{n-1}\big)}
\\
&& -z^2\sum_{i=n+1}^{2n+3}\frac{(\tau^i)^2}{z-\frac{s}{2}+(\tau^i)^2})\big(
\frac{16(n-1)^{n-1}}{n-2}+(z+\frac{s}{2})^{n-1}
\big)\\
&=&(n-1)^n\textcolor{blue}{\big(-\frac{(n-2)^2 s^2}{4}+(2n-8)sz-
\frac{4(n-6)}{n-2}z^2\big)}\sum_{i=n+1}^{2n+3}\frac{(\tau^i)^2}{z-\frac{s}{2}+(\tau^i)^2}\\
&&+4(n-1)^n s z\textcolor{blue}{-16(n-1)^{n-1}z^2}+z^2(z+\frac{s}{2})^{n-1}
\\
&& -z^2\sum_{i=n+1}^{2n+3}\frac{(\tau^i)^2}{z-\frac{s}{2}+(\tau^i)^2})\big(
\frac{16(n-1)^{n-1}}{n-2}+(z+\frac{s}{2})^{n-1}
\big)\\
&=&(n-1)^{\textcolor{blue}{n-1}}\big(-\frac{(n-1)(n-2)^2 s^2}{4}+2(n-1)(n-4)sz\\
&&\textcolor{blue}{-\frac{16}{n-2}z^2}-\frac{4(n-1)(n-6)}{n-2}z^2\big)\sum_{i=n+1}^{2n+3}\frac{(\tau^i)^2}{z-\frac{s}{2}+(\tau^i)^2}\\
&&+4(n-1)^n s z-16(n-1)^{n-1}z^2+z^2(z+\frac{s}{2})^{n-1}
-z^2\textcolor{blue}{(z+\frac{s}{2})^{n-1}}\sum_{i=n+1}^{2n+3}\frac{(\tau^i)^2}{z-\frac{s}{2}+(\tau^i)^2}\\
&=&(n-1)^{n-1}\big(-\frac{(n-1)(n-2)^2 s^2}{4}+2(n-1)(n-4)sz
-4\textcolor{blue}{(n-5)}z^2\big)\sum_{i=n+1}^{2n+3}\frac{(\tau^i)^2}{z-\frac{s}{2}+(\tau^i)^2}\\
&&+4(n-1)^n s z-16(n-1)^{n-1}z^2+z^2(z+\frac{s}{2})^{n-1}
 -z^2(z+\frac{s}{2})^{n-1}\sum_{i=n+1}^{2n+3}\frac{(\tau^i)^2}{z-\frac{s}{2}+(\tau^i)^2}.
\end{eqnarray*}
Hence we obtain (\ref{eq-charPoly-EV-2ndOrder}).
\end{proof}

\subsection{Proof of semisimplicity}

\begin{proposition}\label{prop-charPoly-EV-2ndOrder-simpleRoots}
For general $(\tau^{n+1},\dots,\tau^{2n+3})\in \mathbb{C}^{n+3}$, the polynomial in $z$
\begin{eqnarray}\label{eq-charPoly-EV-2ndOrder-1}
&&\mathcal{P}_n(z,\tau^{n+1},\dots,\tau^{2n+3})=
P_n(z,\tau^{n+1},\dots,\tau^{2n+3})+o(\tau^2)
\nn\\
&=&\bigg((n-1)^{n-1}\big(-\frac{(n-1)(n-2)^2 s^2}{4}+2(n-1)(n-4)sz
-4(n-5)z^2\big)\sum_{i=n+1}^{2n+3}\frac{(\tau^i)^2}{z-\frac{s}{2}+(\tau^i)^2}\nn\\
&&+4(n-1)^n s z-16(n-1)^{n-1}z^2+z^2(z+\frac{s}{2})^{n-1}
 -z^2(z+\frac{s}{2})^{n-1}\sum_{i=n+1}^{2n+3}\frac{(\tau^i)^2}{z-\frac{s}{2}+(\tau^i)^2})\bigg)\nn\\
&&\cdot \prod_{i=n+1}^{2n+3}\big(z-\frac{s}{2}+(\tau^i)^2\big)+o(\tau^2)
\end{eqnarray}
has only simple roots.
\end{proposition}
\begin{proof}
We have
\begin{eqnarray*}
P_n(z,0,\dots,0)=z^{n+5}(-16(n-1)^{n-1}+z^{n-1}).
\end{eqnarray*}
It has $n-1$ nonzero simple roots and a root of multiplicity $n+5$ at 0. We are going to find a (germ of) line $L$ on the $\tau$-space starting at 0, parametrized by $\theta$, and find the branches of solutions on this ray. Let the line $L$ be, say $\tau^i=\alpha_i \theta$, where $\alpha_i\in \mathbb{C}$, for $n+1\leq i\leq 2n+3$ are temporarily not specified.  By the theory of Puiseux expansion, near the origin of $L$, there are solutions in the form 
\begin{equation*}
	z=\sum_{i=1}^{\infty}a_{\frac{i}{M}}\theta^{\frac{i}{M}}
\end{equation*}
where $M$ is a natural number. The index $\frac{i}{M}$ of the least possible nonzero coefficient $a_{\frac{i}{M}}$ can be read out from the expression (\ref{eq-charPoly-EV-2ndOrder}). It is also  given by the \emph{Newton polygon}; indeed, regarding the restriction of $P_n$ on $L$ as a polynomial $P_n(\theta,z)$ the pair of indices that gives the first slope of the Newton polygon is $(0,n+5)$ and $(2,n+4)$. So the slope is $2$, and we get the expansion 
\[
z=a_2 \theta^2+o(\theta^{2}).
\]
To show that $\mathcal{P}_n$ has $n+5$ distinct branches of solutions of $z$ at 0, we need only to show that the equation 
\begin{eqnarray*}
	&&\big(-\frac{(n-1)(n-2)^2 s}{4}+2(n-1)(n-4)sz-4(n-5)z^2\big)\sum_{i=0}^{n+2}\frac{\tau^2_i}{z-\frac{s}{2}+\tau^2_i}\nn\\
&&+4(n-1) s z -16z^2=0
\end{eqnarray*}
for $a_2$ has $n+5$ simple roots, for generic choices of $\tau^{n+1},\dots,\tau^{2n+3}$. This follows  from the following Proposition \ref{prop-secondOrderCoeff-distinctRoots}.
\end{proof}

\begin{lemma}\label{lem-simpleRoots}
When $n\geq 3$ and $a\in \mathbb{Q}$, the equation $z^n-az+1=0$ has only simple roots.
\end{lemma}
\begin{proof}
$z^n-az+1$ and $nz^{n-1}-a$ have common factors only if $a=\frac{n}{n-1}\cdot (n-1)^{\frac{1}{n}}$, which is not rational when $n\geq 3$.
\end{proof}

\begin{proposition}\label{prop-secondOrderCoeff-distinctRoots}
Let $n\geq 4$ be an even natural number.
Let $\varsigma=\varsigma(y_0,\dots,y_{n+2})=y_0+\dots+y_{n+2}$. Then for general $(y_0,\dots,y_{n+2})\in \mathbb{C}^{n+3}$,  the equation 
\begin{eqnarray}\label{eq-secondOrderCoeff-distinctRoots}
&&	\Big(\big(-\frac{(n-1)(n-2)^2 \varsigma^2}{4}+2(n-1)(n-4)\varsigma z
-4(n-5)z^2\big)\sum_{i=0}^{n+2}\frac{y_i}{z-\frac{\varsigma}{2}+y_i}\nn\\
&&+4(n-1) \varsigma z-16z^2\Big)\cdot \prod_{i=0}^{n+2}\big(z-\frac{\varsigma}{2}+y_i\big)=0
\end{eqnarray}
for $z$ has $n+5$ simple roots.
\end{proposition}
\begin{proof}
Let $\zeta_{n+3}=e^{2\pi \sqrt{-1}}$ be a $(n+3)$-th root of unity, and set 
\begin{equation*}
	y_i=\zeta_{n+3}^i+y.
\end{equation*}
It suffices to show that for general $y\in \mathbb{C}$, (\ref{eq-secondOrderCoeff-distinctRoots}) has only simple roots.
Since
\begin{equation*}
	\prod_{i=0}^{n+2}(z+\zeta_{n+3}^i)=z^{n+3}+1,
\end{equation*}
and 
\begin{equation*}
	\sum_{i=0}^{n+2}\frac{1}{z+\zeta_{n+3}^i}=\frac{n+3}{z^{n+3}+1},
\end{equation*}
we have
\begin{equation*}
	\prod_{i=0}^{n+2}(z-\frac{s}{2}+y_i)=(z-\frac{s}{2}+y)^{n+3}+1,
\end{equation*}
and
\begin{equation*}
	\sum_{i=0}^{n+2}\frac{y_i}{z-\frac{s}{2}+y_i}
	=n+3-(z-\frac{s}{2})\sum_{i=0}^{n+2}\frac{1}{z-\frac{s}{2}+y_i}
	=(n+3)(1-\frac{z-\frac{s}{2}}{(z-\frac{s}{2}+y)^{n+3}+1}).
\end{equation*}
So (\ref{eq-secondOrderCoeff-distinctRoots}) reads
\begin{eqnarray}\label{eq-prop-secondOrderCoeff-distinctRoots-1}
	&&\big(-\frac{(n-1)(n-2)^2(n+3)^2y^2}{4}+2(n-1)(n-4)(n+3)yz
-4(n-5)z^2\big)\nn\\
&&\cdot(n+3)((z-\frac{n+1}{2}y)^{n+3}+1-z+\frac{n+3}{2}y) \nn\\
&&+\big(4(n-1)(n+3)y z-16z^2\big)\big((z-\frac{n+1}{2}y)^{n+3}+1\big)
=0.
\end{eqnarray}
When $y=0$, this equation is
\begin{equation*}
	4z^2\big((- n^2+2 n+11) z^{n+3}+( n^2-2 n-15)z+(- n^2 +2 n +11) \big)=0,
\end{equation*}
i.e.
\begin{equation*}
	4(- n^2+2 n+11)z^2\big( z^{n+3}-\frac{n^2-2 n-15}{n^2-2 n-11}z+1 \big)=0.
\end{equation*}
By Lemma \ref{lem-simpleRoots}, the second factor has only simple roots.
We study the branches $z=z(y)$ with $z(0)=0$. They can be written in the form of Puiseux series
\[
z=a_{\frac{1}{2}}y^{\frac{1}{2}}+a_1 y+O(y^{\frac{3}{2}}).
\]
Expanding (\ref{eq-prop-secondOrderCoeff-distinctRoots-1}), the equation for $a_{\frac{1}{2}}$ is
\[
0=-4(n-5)(n+3)a_{\frac{1}{2}}^2-16a_{\frac{1}{2}}^2=0
\]
so $a_{\frac{1}{2}}=0$. One can also argue by using Newton polygon.
Then the equation for $a_1$ is
\begin{eqnarray*}
	&&(n+3)\big(-\frac{(n-1)(n-2)^2(n+3)^2}{4}+2(n-1)(n-4)(n+3)a_1
-4(n-5)a_1^2\big)\\
&&+\big(4(n-1)(n+3)a_1-16a_1^2\big)=0.
\end{eqnarray*}
This is a quadratic equation with a discriminant equal to
\[
64 (-1 + n) (2 + n) (3 + n)^2.
\]
So we have two distinct branches with $z(0)=0$.
\end{proof}

\begin{proof}[Proof of Theorem \ref{thm-semisimplicity}]
By Proposition \ref{prop-charPoly-EV-2ndOrder} and \ref{prop-charPoly-EV-2ndOrder-simpleRoots}, the Euler field $E$ has pairwise distinct eigenvalues on a general point in an open neighborhood of 0. 
\end{proof}


\section{Enumerative geometry on even dimensional intersections of two quadrics}\label{sec:EnumerativeGeometry-Even(2,2)}
In this section we study the special correlator 
\begin{equation}\label{eq-specialCorrelator}
	\langle \epsilon_{1},\dots,\epsilon_{n+3}\rangle^X_{0,n+3,\frac{n}{2}}
\end{equation}
by relating its value to enumerative geometry of $X$. Then we prove the $n=4$ case of Conjecture \ref{conj-unknownCorrelator-Even(2,2)}. 
We begin by recalling some facts about the even dimensional intersections of two quadrics from \cite{Rei72}. Let $\lambda_0,\dots,\lambda_{n+2}\in\mathbb{C}$ be pairwise distinct. Let 
\begin{equation}\label{eq-definingEquationOfQi}
	\varphi_1(Y_0,\dots,Y_{n+2})=\sum_{i=0}^{n+2}Y_i^2,\
	\varphi_2(Y_0,\dots,Y_{n+2})=\sum_{i=0}^{n+2} \lambda_i Y_i^2,
\end{equation}
and $X=\{\varphi_1=\varphi_2=0\}\subset \mathbb{P}^{n+2}$. By \cite[Prop. 2.1]{Rei72},  every smooth complete intersections of two quadrics can be obtained in this way by choosing appropriate coordinates on $\mathbb{P}^{n+2}$.
\begin{lemma}
For $0\leq k\leq \frac{n}{2}$, let $P_k$ be the point in $\mathbb{P}^{n+2}$  whose $i$-th homogeneous coordinate is
\begin{equation*}
\frac{\lambda_i^k}{\sqrt{\prod_{\begin{subarray}{c}0\leq j\leq n+2\\ j\neq i\end{subarray}}
(\lambda_i- \lambda_j)}}. 
\end{equation*}
Then
\begin{equation}\label{eq-vanish-varphiPk}
	\varphi_1(P_k)=\varphi_2(P_k)=0,
\end{equation}
and $P_0,\dots,P_{\frac{n}{2}}$ span an $\frac{n}{2}$-plane contained in $X$.  Here we  choose a root $\sqrt{\lambda_i- \lambda_j}$ uniformly in the expressions of $P_0,\dots,P_k$, for every pair ${i,j}$, for $0\leq i\neq j\leq n+2$.
\end{lemma}
\begin{proof}
Directly check (\ref{eq-vanish-varphiPk}) using 
\[
\frac{\lambda^l}{\prod_{\begin{subarray}{c}0\leq j\leq n+2\\ j\neq i\end{subarray}}(\lambda_i- \lambda_j)}=0
\]
for $0\leq l\leq n+1$.
\end{proof}

We denote this $\frac{n}{2}$-plane by $S$. Then $S$ is defined by the following linear equations:
\begin{equation*}
\sum_{i=0}^{n+2}\frac{\lambda_i^k}{\sqrt{\prod_{\begin{subarray}{c}0\leq j\leq n+2\\ j\neq i\end{subarray}}
(\lambda_i- \lambda_j)}}Y_i=0,\ \mbox{for}\ 0\leq k\leq \frac{n}{2}+1.
\end{equation*}
Make a change of coordinates
\begin{equation}\label{eq-Wcoordinates}
	W_i=\frac{Y_i}{\sqrt{\prod_{\begin{subarray}{c}0\leq j\leq n+2\\ j\neq i\end{subarray}}
(\lambda_i- \lambda_j)}}.
\end{equation}
Then $S$ is defined by
\begin{equation}\label{eq-definingEquations-WCoordinates}
	\sum_{i=0}^{n+2}\lambda_i^k W_i=0,\ \mbox{for}\ 0\leq k\leq \frac{n}{2}+1.
\end{equation}

Recall from \cite[Lemma 3.10]{Rei72}:
\begin{lemma}\label{lem-intersectionNumber-middelDimension}
\begin{enumerate}
	\item[(i)] Let $T$ be a $\frac{n}{2}$-plane contained in $X$. Then $[T]\cdot \sfh_{n/2}=1$.
	\item[(ii)] Let $T_1$ and $T_2$ be $\frac{n}{2}$-planes contained in $X$.
	Suppose $\dim(T_1\cap T_2)=r$, then 
	\begin{equation*}
			[T_1]\cdot [T_2]=(-1)^{r}(\lfloor\frac{r}{2}\rfloor+1).
	\end{equation*}
\end{enumerate}
\end{lemma}

\begin{notation}
For integers $a\leq b$ let $[a,b]$ be the set of integers $c$ satisfying $a\leq c\leq b$. 
For a subset $I\subset [0,n+2]$, let $S_I$ be the $\frac{n}{2}$-plane obtained by reversing the sign of the $i$-th homogeneous coordinate of the points on $S$ for all $i\in I$. Denote the complement of $I$ by $C(I)$. Then $S_I=S_{C(I)}$. 
\end{notation}
By \cite[Theorem 3.8]{Rei72}, every $\frac{n}{2}$-planes of $X$ is equal to $S_I$ for some $I\subset [0,n+2]$.
In particular, let $S_i=S_{\{i\}}$. Let $\varsigma_i=[S_i]$, the homology class of $S_i$. Let $\varsigma=[S]$. By \cite[Lemma 3.13]{Rei72},
\[
\sfh_{\frac{n}{2}},\varsigma_0,\varsigma_1,\dots,\varsigma_{n+2}
\]
is a basis of $H^{\frac{n}{2}}(X;\mathbb{Q})$, and 
\begin{equation*}
	\varsigma=\frac{\frac{n}{2}+1}{n+1}\sfh_{n/2}-\frac{1}{n+1}\sum_{i=0}^{n+2}\varsigma_i.
\end{equation*}
For $I\subset [0,n+2]$, let $a(I)=(x_0,\dots,x_{n+2})\in \mathbb{Z}^{n+3}$ such that 
\[
 x_i=\begin{cases}
 1,& \mbox{if}\ i\not\in I,\\
 -1,& \mbox{if}\ i\in I.
 \end{cases}
 \] 
 For $I$ and $J \subset [0,n+2]$, suppose $a(I)=(x_0,\dots,x_{n+2})$ and $a(J)=(y_0,\dots,y_{n+2})$. Define
 \[
 a(I,J):=(x_0 y_0,\dots,x_{n+2}y_{n+2}).
 \]
 \begin{lemma}\label{lem-intersectionDimension}
Let $m(I,J)$ be the number of $-1$ in the $(n+3)$-tuple $a(I,J)$. Then
\begin{equation*}
	\dim S_I\cap S_J=\begin{cases}
	\frac{n}{2}-m(I,J),& \mbox{if}\ m(I,J)\leq \frac{n}{2},\\
	-1,& \mbox{if}\ m(I,J)=\frac{n}{2}+1\ \mbox{or}\ \frac{n}{2}+2,\\
	m(I,J)-\frac{n}{2}-3,& \mbox{if}\ m(I,J)\geq \frac{n}{2}+3,
	\end{cases}	
\end{equation*}
where $\dim S_I\cap S_J=-1$ means that the intersection is empty.
 \end{lemma}
\begin{proof}
One can show this using directly the defining equations (\ref{eq-definingEquations-WCoordinates}) of $S$, or from \cite[Theorem 3.8]{Rei72}.
\end{proof} 
By Lemma \ref{lem-intersectionNumber-middelDimension} and Lemma \ref{lem-intersectionDimension}, 
\begin{equation}\label{eq-intersectionNumber-Si-Sj}
	\varsigma_i\cdot \varsigma_j=\begin{cases}
	(-1)^{\frac{n}{2}-2}(\lfloor \frac{\frac{n}{2}-2}{2}\rfloor+1), & \mbox{if}\ i\neq j,\\
	(-1)^{\frac{n}{2}}(\lfloor \frac{n}{4}\rfloor+1),& \mbox{if}\ i=j.
	\end{cases}
\end{equation}

\subsection{An explicit orthonormal basis}\label{sec:explictD-Lattice}
For $1\leq i\leq n+3$, we define
 \begin{equation*}
	\varepsilon_i=\varsigma_{i-1}-\frac{1}{n+1}\sum_{i=0}^{n+2}\varsigma_{i}+\frac{1}{2(n+1)}\sfh_{n/2}.
\end{equation*}
\begin{lemma}\label{lem-expressingEpsilonbyMiddleDimPlanes}
\begin{equation*}
	\varsigma_{i-1}=\varepsilon_i-\frac{1}{2}\sum_{i=1}^{n+3}\varepsilon_i+\frac{\sfh_{n/2}}{4},\ \mbox{for}\ 1\leq i\leq n+3,
\end{equation*}
\begin{equation*}
	\varsigma=\frac{1}{4}\sfh_{n/2}+\frac{1}{2}\sum_{i=1}^{n+3}\varepsilon_i.
\end{equation*}
\end{lemma}
\begin{lemma}\label{lem-intersectionNumber-epsilon}
For $1\leq i\leq n+3$,
\begin{equation}\label{eq-intersectionNumber-epsilon-h}
\varepsilon_i\cdot \sfh_{n/2}=0,
\end{equation}
\begin{equation}\label{eq-intersectionNumber-epsiloni-self}
	\varepsilon_i \cdot \varepsilon_i=(-1)^{\frac{n}{2}},
\end{equation}
for $i\neq j$,
\begin{equation}\label{eq-intersectionNumber-epsiloni-epsilonj}
	\varepsilon_i \cdot \varepsilon_j=0.
\end{equation}
\end{lemma}
\begin{proof}
(\ref{eq-intersectionNumber-epsilon-h}) follows from Lemma \ref{lem-intersectionNumber-middelDimension} (i).
For (\ref{eq-intersectionNumber-epsiloni-self}) and (\ref{eq-intersectionNumber-epsiloni-epsilonj}) we use Lemma \ref{lem-intersectionNumber-middelDimension} (i) and  (\ref{eq-intersectionNumber-Si-Sj}): for $1\leq i\neq j\leq n+3$,
\begin{eqnarray*}
&& \varepsilon_i \cdot \varepsilon_j\\
&=& \big(1-\frac{2(n+2)}{n+1}+\frac{(n+3)(n+2)}{(n+1)^2}\big)\cdot (-1)^{\frac{n}{2}-2}(\lfloor \frac{\frac{n}{2}-2}{2}\rfloor+1)\\
&&+ \big(-\frac{2}{n+1}+\frac{n+3}{(n+1)^2}\big)\cdot (-1)^{\frac{n}{2}}(\lfloor \frac{n}{4}\rfloor+1)\\
&&+ 2(1-\frac{n+3}{n+1})\frac{1}{2(n+1)}+\frac{1}{4(n+1)^2}\cdot 4\\
&=& \begin{cases}
\frac{n+3}{(n+1)^2}\cdot \frac{n}{4}-\frac{n-1}{(n+1)^2}\cdot (\frac{n}{4}+1)-\frac{1}{(n+1)^2},& \mbox{if}\ n\equiv 0 \mod 4,\\
-\frac{n+3}{(n+1)^2}\cdot \frac{n-2}{4}+\frac{n-1}{(n+1)^2}\cdot (\frac{n-2}{4}+1)-\frac{1}{(n+1)^2},& \mbox{if}\ n\equiv 2 \mod 4,\\
\end{cases}\\
&=& 0,
\end{eqnarray*}
and
\begin{eqnarray*}
&& \varepsilon_i \cdot \varepsilon_i\\
&=& \big(-\frac{2(n+2)}{n+1}+\frac{(n+3)(n+2)}{(n+1)^2}\big)\cdot (-1)^{\frac{n}{2}-2}(\lfloor \frac{\frac{n}{2}-2}{2}\rfloor+1)\\
&&+ \big(1-\frac{2}{n+1}+\frac{n+3}{(n+1)^2}\big)\cdot (-1)^{\frac{n}{2}}(\lfloor \frac{n}{4}\rfloor+1)\\
&&+ 2(1-\frac{n+3}{n+1})\frac{1}{2(n+1)}+\frac{1}{4(n+1)^2}\cdot 4\\
&=& \begin{cases}
\frac{(-n+1)(n+2)}{(n+1)^2}\cdot \frac{n}{4}+\frac{n^2+n+2}{(n+1)^2}\cdot (\frac{n}{4}+1)-\frac{1}{(n+1)^2},& \mbox{if}\ n\equiv 0 \mod 4,\\
-\frac{(-n+1)(n+2)}{(n+1)^2}\cdot \frac{n-2}{4}-\frac{n^2+n+2}{(n+1)^2}\cdot (\frac{n-2}{4}+1)-\frac{1}{(n+1)^2},& \mbox{if}\ n\equiv 2 \mod 4,\\
\end{cases}\\
&=& (-1)^{\frac{n}{2}}.
\end{eqnarray*}
\end{proof}
We define 
\begin{equation}\label{eq-roots-D-usingExplicitEpsilonClass}
	\begin{cases}
	\alpha_i=\varepsilon_{i}-\varepsilon_{i+1}\ \mbox{for}\ 1\leq i\leq n+2,\\
	\alpha_{n+3}=\varepsilon_{n+2}+\varepsilon_{n+3}.
	\end{cases}
\end{equation} 
Then from Lemma \ref{lem-expressingEpsilonbyMiddleDimPlanes} we have
\begin{equation}\label{eq-roots-D-usingExplicitMiddleDimPlanes}
	\begin{cases}
	\alpha_i=\varsigma_{i-1}-\varsigma_{i}\ \mbox{for}\ 1\leq i\leq n+2,\\
	\alpha_{n+3}=\varsigma_{n+1}+\varsigma_{n+2}+ 2\varsigma-\sfh_{n/2}.
	\end{cases}
\end{equation} 
So $\alpha_i\in H^*(X;\mathbb{Z})$. Using Lemma (\ref{lem-intersectionNumber-epsilon}), and comparing (\ref{eq-roots-D-usingExplicitEpsilonClass}) and (\ref{eq-roots-D}),  one sees  that the group generated  by the reflections with respect to $\alpha_i$'s is the Weyl group $D_{n+3}$. By the Picard-Lefschetz formula, one can take the class $\alpha_i$  in Section \ref{sec:monodromy-lattice} to be $\alpha_i$ defined here. This justifies the notations. Moreover, because of (\ref{eq-intersectionNumber-epsiloni-self}), we define an  orthonormal basis $\epsilon_i$ of $H^*_{\mathrm{prim}}(X)$ exactly as (\ref{eq-normalizedOrthonormalBasis}).

\subsection{A potentially enumerative correlator}\label{sec:enumerativeCorrelator}
We begin with giving a working definition of \emph{enumerative correlators}.
\begin{definition}\label{def-enumerativeCorrelator}
Let $n\geq 4$ be an even natural number, and $X$ be an $n$-dimensional smooth complete intersection of two quadrics in $\mathbb{P}^{n+2}$.  Denote the $i$-th projection from $X^{n+3}$ to $X$ by $q_i$. Consider the product of the evaluation morphisms
\begin{equation*}
	\mathrm{ev}_1\times\cdots \mathrm{ev}_{n+3}: \Mbar_{0,n+3}(X,\frac{n}{2})\rightarrow X^{n+3}.
\end{equation*}
Let $I_1,\dots, I_{n+3}\subset [0,n+2]$. We say that the correlator 
\begin{equation}\label{eq-enumerativeCorrelator-0}
	\langle \varsigma_{I_1},\dots,\varsigma_{I_{n+3}}\rangle
\end{equation}
is \emph{enumerative} if there exists an irreducible component $M$ of $\Mbar_{0,n+3}(X,\frac{n}{2})$ satisfying the following:
\begin{enumerate}
	\item[(i)] $\dim M$ equals the expected dimension.
	\item[(ii)] The cycles  $(\mathrm{ev}_1\times\cdots \mathrm{ev}_{n+3})(M)$ and $q_1^{-1} S_{I_1}$,\dots, $q_{n+3}^{-1}S_{I_{n+3}}$ intersect \emph{properly}, i.e. the dimension of their (scheme theoretic) intersection is 0.
	\item[(iii)] Each irreducible component of $\Mbar_{0,n+3}(X,\frac{n}{2})$ other than $M$ has empty intersection with $q_1^{-1} S_{I_1}$,\dots, $q_{n+3}^{-1} S_{I_n+3}$. 
\end{enumerate}
\end{definition}

This definition is of course not standard in any sense. It just facilitates the following presentation.  

Our strategy to compute the special correlator by  the enumerative geometry of $X$ consists of three steps:
\begin{enumerate}
 	\item Select $I_1,\dots, I_{n+3}\subset [0,n+2]$, such that the correlator (\ref{eq-enumerativeCorrelator-0}) is enumerative.
	\item Express (\ref{eq-enumerativeCorrelator-0}) in terms of the special correlator, using Lemma \ref{lem-expressingEpsilonbyMiddleDimPlanes}.
	\item Solve the corresponding enumerative problem by counting curves. More precisely, compute the intersection multiplicities of the intersection 
	\begin{equation*}
		(\mathrm{ev}_1\times\cdots \mathrm{ev}_{n+3})_*[M]\cap q_1^{*} [S_{I_1}]\cap\dots\cap q_{n+3}^{*} [S_{I_{n+3}}]
	\end{equation*}
	 in condition (ii) above. 
 \end{enumerate} 

 \begin{example}\label{exp-nonEnumerative-correlator}
 The correlator 
 \begin{equation*}
 	\langle \varsigma_{0},\dots,\varsigma_{n+3}\rangle
 \end{equation*}
 should not be enumerative in general. For example, let $n=4$. Then the intersection $S\cap S_i$ is a line, for  $0\leq i\leq 6$.  The moduli space of conics on $S$ passing through the seven lines has a positive dimension. So there are infinitely many conics passing through $S_0,\dots,S_6$. Then the conditions (ii) and (iii) in Definition \ref{def-enumerativeCorrelator} cannot be true simultaneously.
 \end{example}

Stimulated by Example \ref{exp-nonEnumerative-correlator}, we propose a possibly enumerative correlator.
Let $S_{[i,i+k-1]}$ be the $\frac{n}{2}$-plane $S_{i,i+1,\dots,i+k-1}$ in $X$, where we understand the indices $i$ in the subscript in mod $n+3$ sense. For example, when $n=4$, $S_{[6,6+1]}=S_{6,0}$, $S_{[5,5+3]}=S_{5,6,0,1}$.
\begin{lemma}\label{lem-uniqueS}
$S$ is the only $\frac{n}{2}$-plane in $X$ that has non-empty intersections with each of $S_{[i,i+\frac{n}{2}-1]}$, for $0\leq i\leq n+2$. Moreover $S$ meets $S_{[i,i+\frac{n}{2}-1]}$ at exactly one point.
\end{lemma}
\begin{proof}
The second statement follows from Lemma \ref{lem-intersectionDimension}.

We show the first statement.
Let $I\subset [0,n+2]$, and suppose $a(I)=(x_0,\dots,x_{n+2})$. Let $J_i=[i,i+\frac{n}{2}-1]$ in the mod $n+3$ sense. By Lemma \ref{lem-intersectionDimension}, the first statement is equivalent to that when $a(I)\neq (1,1,\dots,1)$ or $(-1,-1,\dots,-1)$, there exists $i$ such that $a(I,J_i)=\frac{n}{2}+1$ or $\frac{n}{2}+2$, or  equivalently the sum of all components of $a(I,J_i)$ is equal to $\pm 1$. Put
\[
p(x_0,\dots,x_{n+2})=\sum_{i=0}^{n+2}\big(a(I,J_i)^2-1\big).
\]
We regard $p(x_0,\dots,x_{n+2})$ as a polynomial of indeterminates $x_0,\dots,x_{n+2}$.
Then we are left to show that when $(x_0,\dots,x_{n+2})\in \{1,-1\}^{n+3}$ and $(x_0,\dots,x_{n+2})\neq (1,1,\dots,1)$ or $(-1,-1,\dots,-1)$, 
\begin{equation*}
	p(x_0,\dots,x_{n+2})=0.
\end{equation*}
Put 
\[
y_i=\sum_{j\in J_i}x_j.
\]
Then
\begin{equation*}
	a(I,J_i)=\sum_{i=0}^{n+2}x_i-2y_i.
\end{equation*}
So $p(x_0,\dots,x_{n+2})$ is manifestly symmetric in $y_i$. It follows that $p(x_0,\dots,x_{n+2})$ is  symmetric in $x_i$, rather than only cyclicly symmetric. So it suffices to show the statement for all $I\subset[0,n+2]$ of the  form $I=[0,k]$ for some $0\leq k<n+2$. Let 
\[
i=\begin{cases}\frac{k}{2}+1,& \mbox{if}\ 2|k,\\
\frac{k+3}{2},& \mbox{if}\ 2\neq k.
\end{cases}
\]
Then one checks that $a(I,J_i)=\pm 1$, and thus $p(x_0,\dots,x_{n+2})=0$. 
\end{proof}

\begin{lemma}
\begin{equation}\label{eq-sigmaInterval-epsilon}
	\varsigma_{[i-1,i-1+\frac{n}{2}-1]}=\frac{1}{4}\sfh_{n/2}+(-1)^{\frac{n}{2}}\big(\frac{1}{2}\sum_{i=1}^{n+3}\varepsilon_i-\sum_{j=0}^{\frac{n}{2}-1}\varepsilon_{i+j}\big)
\end{equation}
\end{lemma}
\begin{proof}
By \cite[Lemma 3.12]{Rei72},
\begin{equation*}
	\varsigma_{[i-1,i-1+\frac{n}{2}-1]}=
	 (-1)^{\frac{n}{2}}\big(\lfloor \frac{n}{4}\rfloor \sfh_{n/2}-(\frac{n}{2}-1)\varsigma-\sum_{j=0}^{\frac{n}{2}-1}\varsigma_{i+j}\big).
\end{equation*}
Then using Lemma \ref{lem-expressingEpsilonbyMiddleDimPlanes} we compute
\begin{eqnarray*}
&& \varsigma_{[i-1,i-1+\frac{n}{2}-1]}\\
&=&(-1)^{\frac{n}{2}}\big(\lfloor \frac{n}{4}\rfloor \sfh_{n/2}-(\frac{n}{2}-1)(\frac{1}{4}\sfh_{n/2}+\frac{1}{2}\sum_{i=1}^{n+3}\varepsilon_i)-\sum_{j=0}^{\frac{n}{2}-1}( \varepsilon_{i+j}-\frac{1}{2}\sum_{i=1}^{n+3}\varepsilon_i+\frac{\sfh_{n/2}}{4})\big)\\
&=&(-1)^{\frac{n}{2}}\big((\lfloor \frac{n}{4}\rfloor-\frac{n}{4}+\frac{1}{4} )\sfh_{n/2}+\frac{1}{2}\sum_{i=1}^{n+3}\varepsilon_i-\sum_{j=0}^{\frac{n}{2}-1}\varepsilon_{i+j}\big)\\
&=& \frac{1}{4}\sfh_{n/2}+(-1)^{\frac{n}{2}}\big(\frac{1}{2}\sum_{i=1}^{n+3}\varepsilon_i-\sum_{j=0}^{\frac{n}{2}-1}\varepsilon_{i+j}\big).
\end{eqnarray*}
\end{proof}

\begin{definition}
Let $n\geq 4$ be an even integer. 
Let $X$ be an $n$-dimensional smooth complete intersection of two quadrics in $\mathbb{P}^{n+2}$. We define
\begin{equation}\label{eq-possiblyEnumerativeCorrelator}
	f(n):=\langle   \varsigma_{[0,\frac{n}{2}-1]},\dots,\varsigma_{[n+2,n+2+\frac{n}{2}-1]}\rangle_{0,n+3,\frac{n}{2}}.
\end{equation}

In Section \ref{sec:specialCorrelator-Dim4}, we will show that for general $4$-dimensional smooth complete intersection of two quadrics, $f(n)$ is an enumerative correlator. In following we relate the value of $f(n)$ to the special correlator.

\begin{lemma}\label{lem-correlatorsOfLength7-4dim}
Let $X$ be a 4-dim smooth complete intersection of two quadrics in $\mathbb{P}^6$. Then
\begin{eqnarray*}
&\langle\sfh_2,\sfh_2,\sfh_2,\sfh_2,\sfh_2,\sfh_2,\sfh_2\rangle=46656,\
\langle \sfh_2,\sfh_2,\sfh_2,\sfh_2,\sfh_2,\varepsilon_1,\varepsilon_1\rangle=-624,\
\langle \sfh_2,\sfh_2,\sfh_2,\varepsilon_1,\varepsilon_1,\varepsilon_1,\varepsilon_1\rangle=36,\\
&\langle \sfh_2,\sfh_2,\sfh_2,\varepsilon_1,\varepsilon_1,\varepsilon_2,\varepsilon_2\rangle=4,\
\langle \sfh_2,\varepsilon_1,\varepsilon_1,\varepsilon_1,\varepsilon_1,\varepsilon_1,\varepsilon_1\rangle=-7,\
\langle \sfh_2,\varepsilon_1,\varepsilon_1,\varepsilon_1,\varepsilon_1,\varepsilon_2,\varepsilon_2\rangle=1,\\
&\langle \sfh_2,\varepsilon_1,\varepsilon_1,\varepsilon_2,\varepsilon_2,\varepsilon_3,\varepsilon_3\rangle=1.
\end{eqnarray*}
\end{lemma}
\begin{proof}
This is a computation via our package \texttt{QuantumCohomologyFanoCompleteIntersection} in Macaulay2 based on the algorithm in the proof of Theorem \ref{thm-reconstruction-even(2,2)}; 
see the examples after Algorithm \ref{algorithm-correlator-even(2,2)}.
\end{proof}

\begin{lemma}\label{lem-relate-specialCorrelatorToEnumerativeCorrelator}
Let $X$ be a 4-dim smooth complete intersection of two quadrics in $\mathbb{P}^6$. Then $\langle \varepsilon_1,\varepsilon_2,\varepsilon_3,\varepsilon_4,\varepsilon_5,\varepsilon_6,\varepsilon_7\rangle=\frac{1}{2}$ if only if 
\begin{equation}\label{eq-enumerativeCorrelator-Dim4}
	 \langle \varsigma_{01},\varsigma_{12},\varsigma_{23},\varsigma_{34},\varsigma_{45},\varsigma_{56},\varsigma_{60}\rangle_{0,7,2}=1.
\end{equation}
\end{lemma}
\begin{proof}
We compute by (\ref{eq-sigmaInterval-epsilon}) and Lemma \ref{lem-correlatorsOfLength7-4dim},
\begin{eqnarray}\label{eq-enumerativeCorrelator-inTermsOfSpecialCorr-dim4}
&& \langle \varsigma_{01},\varsigma_{12},\varsigma_{23},\varsigma_{34},\varsigma_{45},\varsigma_{56},\varsigma_{60}\rangle_{0,7,2}\nn\\
&=& \frac{1}{16384}\langle \sfh_2,\dots,\sfh_2\rangle_{0,7,2}+\frac{7}{4096} \langle \sfh_2,\sfh_2,\sfh_2,\sfh_2,\sfh_2,\varepsilon_1,\varepsilon_1\rangle\nn\\
&& -\frac{35}{1024}\langle \sfh_2,\sfh_2,\sfh_2,\varepsilon_1,\varepsilon_1,\varepsilon_1,\varepsilon_1\rangle
+\frac{77}{512}\langle \sfh_2,\sfh_2,\sfh_2,\varepsilon_1,\varepsilon_1,\varepsilon_2,\varepsilon_2\rangle\nn\\
&&+\frac{21}{256}\langle \sfh_2,\varepsilon_1,\varepsilon_1,\varepsilon_1,\varepsilon_1,\varepsilon_1,\varepsilon_1\rangle
-\frac{63}{128}\langle \sfh_2,\varepsilon_1,\varepsilon_1,\varepsilon_1,\varepsilon_1,\varepsilon_2,\varepsilon_2\rangle\nn\\
&&+\frac{77}{128}\langle \sfh_2,\varepsilon_1,\varepsilon_1,\varepsilon_2,\varepsilon_2,\varepsilon_3,\varepsilon_3\rangle
+\frac{5}{8}\langle \varepsilon_1,\varepsilon_2,\varepsilon_3,\varepsilon_4,\varepsilon_5,\varepsilon_6,\varepsilon_7\rangle\nn\\
&=& \frac{1}{16384}\times 46656+\frac{7}{4096} \times(-624) -\frac{35}{1024}\times 36
+\frac{77}{512}\times 4\nn\\
&&+\frac{21}{256}\times (-7)
-\frac{63}{128}\times 1
+\frac{77}{128}\times 1
+\frac{5}{8}\langle \varepsilon_1,\varepsilon_2,\varepsilon_3,\varepsilon_4,\varepsilon_5,\varepsilon_6,\varepsilon_7\rangle\nn\\
 &=& \frac{11}{16}+\frac{5}{8}\langle \varepsilon_1,\varepsilon_2,\varepsilon_3,\varepsilon_4,\varepsilon_5,\varepsilon_6,\varepsilon_7\rangle.
\end{eqnarray}
\end{proof}
We have defined functions in the package \texttt{QuantumCohomologyFanoCompleteIntersection} to perform similar computations for any even $n$. We present the cases $n=6$ and $8$ in the following.
\end{definition}

\begin{example}\label{example-f(6)}
\[
f(6)=\frac{19303}{16}+\frac{39}{8}\langle \varepsilon_1,\dots,\varepsilon_9\rangle.
\]
So $\langle \varepsilon_1,\dots,\varepsilon_9\rangle=-\frac{1}{2}$ iff $f(6)=1204$. Note that by the integrality of genus 0 Gromov-Witten invariants of semipositive symplectic manifolds, $f(n)$ is an integer. This is why we conjecture $\langle \varepsilon_1,\dots,\varepsilon_9\rangle=-\frac{1}{2}$, not $\frac{1}{2}$.
\end{example}

\begin{example}\label{example-f(8)}
\[
f(8)=\frac{6441821}{4}+\frac{135}{2}\langle \varepsilon_1,\dots,\varepsilon_{11}\rangle.
\]
So $\langle \varepsilon_1,\dots,\varepsilon_{11}\rangle=\frac{1}{2}$ iff $f(8)=1610489$.
\end{example}

\begin{example}\label{example-f(10)}
\[
f(10)=\frac{16232959575}{4}+\frac{2467}{2}\langle \varepsilon_1,\dots,\varepsilon_{13}\rangle.
\]
So $\langle \varepsilon_1,\dots,\varepsilon_{13}\rangle=-\frac{1}{2}$ iff $f(10)=4058239277$.
\end{example}

\begin{remark}\label{rem:compute-specialCorrelator-fromHigherGenusInv}
From (\ref{eq-enumerativeCorrelator-inTermsOfSpecialCorr-dim4}), Example \ref{example-f(6)}-\ref{example-f(10)}, we see another possible approach to computing the special correlator. As a consequence of  Theorem \ref{thm-semisimplicity}, we can compute Gromov-Witten invariants of $X$ in all genera using Givental's formalism (\cite{Giv01},\cite{Tel12}), in terms of functions of the special correlator $\langle \epsilon_{1},\dots,\epsilon_{n+3}\rangle_{0,n+3,\frac{n}{2}}$ or equivalently $\langle \varepsilon_{1},\dots,\varepsilon_{n+3}\rangle_{0,n+3,\frac{n}{2}}$. In genus $g=n+3$  there exist Gromov-Witten invariants with only ambient insertions that depend on $\langle \epsilon_{1},\dots,\epsilon_{n+3}\rangle_{0,n+3,\frac{n}{2}}$. If one can compute such Gromov-Witten invariants with only ambient insertions, one can determine the special correlator, but only up to signs for the same reason as in Lemma \ref{lem-specialCorrelator-freedomOfSign}. In dimensions 4, 6, 8, and 10, the integrality in genus 0 can then  determine the sign of the special correlator. 
We expect that this is the case for all even $n\geq 4$.
\end{remark}

\subsection{Counting conics: a preparation}\label{sec:countingConics-preparation}

From this section on, we will focus on the dimension 4 case, and prove (\ref{eq-enumerativeCorrelator-Dim4}). Then we need to count the number of conics on $X$ passing through $S_{01},\dots,S\cap S_{56},S\cap S_{60}$; we will see in Theorem \ref{thm-countingConics-4dim-general} a precise argument relating such counting to  $f(n)$.

\begin{proposition}\label{prop-noConicOnS}
For general choices of $\lambda_0,\dots,\lambda_6$, there is no conic on $S$ passing through the 7 points $S\cap S_{01},\dots,S\cap S_{56},S\cap S_{60}$.
\end{proposition}
\begin{proof}
Make a change of coordinates
\begin{equation}\label{eq-Zcoordinates}
	Z_i=\sqrt{\prod_{\begin{subarray}{c}0\leq j\leq n+2\\ j\neq i\end{subarray}}
(\lambda_i- \lambda_j)}\cdot Y_i.
\end{equation}
Then 
$S$ is spanned by
\[
(1,1,1,1,1,1,1),\ (\lambda_1,\lambda_2,\lambda_3,\lambda_4,\lambda_5,\lambda_6,\lambda_7),\
 (\lambda_1^2,\lambda_2^2,\lambda_3^2,\lambda_4^2,\lambda_5^2,\lambda_6^2,\lambda_7^2),
\]
and $S_{01}$ is spanned by
\[
(-1,-1,1,1,1,1,1),\ (-\lambda_1,-\lambda_2,\lambda_3,\lambda_4,\lambda_5,\lambda_6,\lambda_7),\
 (-\lambda_1^2,-\lambda_2^2,\lambda_3^2,\lambda_4^2,\lambda_5^2,\lambda_6^2,\lambda_7^2),
\]
and similar for other $S_{i,i+1}$.
Let $X,Y,Z$ be homogeneous coordinates on $S$ such that $[X:Y:Z]$ is identified with
\[
(1,1,1,1,1,1,1)\cdot X+ (\lambda_1,\lambda_2,\lambda_3,\lambda_4,\lambda_5,\lambda_6,\lambda_7)\cdot Y+
 (\lambda_1^2,\lambda_2^2,\lambda_3^2,\lambda_4^2,\lambda_5^2,\lambda_6^2,\lambda_7^2)\cdot Z
\]
in the $Z_i$ coordinates.
Then 
\begin{equation*}
	S\cap S_{ij}=[\lambda_i \lambda_j,-\lambda_i- \lambda_j,1].
\end{equation*}
Let $A_{i,j}=(\lambda_i \lambda_j,-\lambda_i- \lambda_j)$. There exists a conic on $S$ passing through the 7 points  with coordinates $A_{0,1},A_{1,2},\dots,A_{6,0}$ if and only if the matrix
\begin{equation}\label{eq-matrix-conicOnS}
	\left(
\begin{array}{cccccc}
 \lambda_{0}^2 \lambda_{1}^2 & (-\lambda_{0}-\lambda_{1})^2 & 1 & \lambda_{0} \lambda_{1} (-\lambda_{0}-\lambda_{1}) & \lambda_{0} \lambda_{1} & -\lambda_{0}-\lambda_{1} \\
 \lambda_{1}^2 \lambda_{2}^2 & (-\lambda_{1}-\lambda_{2})^2 & 1 & \lambda_{1} \lambda_{2} (-\lambda_{1}-\lambda_{2}) & \lambda_{1} \lambda_{2} & -\lambda_{1}-\lambda_{2} \\
 \lambda_{2}^2 \lambda_{3}^2 & (-\lambda_{2}-\lambda_{3})^2 & 1 & \lambda_{2} \lambda_{3} (-\lambda_{2}-\lambda_{3}) & \lambda_{2} \lambda_{3} & -\lambda_{2}-\lambda_{3} \\
 \lambda_{3}^2 \lambda_{4}^2 & (-\lambda_{3}-\lambda_{4})^2 & 1 & \lambda_{3} \lambda_{4} (-\lambda_{3}-\lambda_{4}) & \lambda_{3} \lambda_{4} & -\lambda_{3}-\lambda_{4} \\
 \lambda_{4}^2 \lambda_{5}^2 & (-\lambda_{4}-\lambda_{5})^2 & 1 & \lambda_{4} \lambda_{5} (-\lambda_{4}-\lambda_{5}) & \lambda_{4} \lambda_{5} & -\lambda_{4}-\lambda_{5} \\
 \lambda_{5}^2 \lambda_{6}^2 & (-\lambda_{5}-\lambda_{6})^2 & 1 & \lambda_{5} \lambda_{6} (-\lambda_{5}-\lambda_{6}) & \lambda_{5} \lambda_{6} & -\lambda_{5}-\lambda_{6} \\
 \lambda_{0}^2 \lambda_{6}^2 & (-\lambda_{0}-\lambda_{6})^2 & 1 & \lambda_{0} \lambda_{6} (-\lambda_{0}-\lambda_{6}) & \lambda_{0} \lambda_{6} & -\lambda_{0}-\lambda_{6} \\
\end{array}
\right)
\end{equation}
has rank $<6$. One can take for example $(\lambda_0,\dots,\lambda_6)=(1,2,3,4,5,6,7)$ and check that the matrix (\ref{eq-matrix-conicOnS}) has rank 6. Then the set of parameters $(\lambda_0,\dots,\lambda_6)$ satisfying that  rank of (\ref{eq-matrix-conicOnS}) $=6$ is Zariski dense in $\mathbb{C}^7$.
\end{proof}

Every conic in a projective space lies on a plane. When a conic is not a double line, it spans a unique plane. To find conics on $X$ passing through the planes $S_{01},\dots,S_{60}$, we will first find all the planes in $\mathbb{P}^6$ that meets $S_{01},\dots,S_{60}$. Furthermore, by Lemma \ref{lem-uniqueS} and Proposition \ref{prop-noConicOnS}, we need to find planes $\Sigma\not\subset X$ that meets $S_{01},\dots,S_{60}$.

We work in the $W$-coordinates (\ref{eq-Wcoordinates}). Then $S_{j,j+1}$ is defined by 
\begin{equation*}
	N_j \cdot (W_0,\dots, W_6)^{\mathrm{T}}=0
\end{equation*}
where (according to (\ref{eq-definingEquations-WCoordinates}))
\begin{equation*}
	N_j=\begin{pmatrix}
	1 &         \dots & 1 		    & -1         & -1              & 1             &  \dots& 1 \\
			\vspace{0.1cm}
	\lambda_0 & \dots & \lambda_{j-1} & -\lambda_j & - \lambda_{j+1} & \lambda_{j+2} & \dots & \lambda_{6}\\
		\vspace{0.1cm}
	\lambda_0^2 & \dots & \lambda_{j-1}^2 & -\lambda_j^2 & - \lambda_{j+1}^2 & \lambda_{j+2}^2& \dots & \lambda_{6}^2\\
	\lambda_0^3 & \dots & \lambda_{j-1}^3 & -\lambda_j^3 & - \lambda_{j+1}^3 & \lambda_{j+2}^3 & \dots & \lambda_{6}^3
	\end{pmatrix}
\end{equation*}
where the subscripts of $\lambda_i$ are understood in the mod $n+3$ sense. Let $B=(b_{i,j})_{0\leq i\leq 6,0\leq j\leq 3}$. The plane defined by $B.\{W_0,\dots,W_6\}^{\mathrm{T}}=0$ meets $S_{01},\dots,S_{60}$ if and only if for $0\leq i\leq 6$
the matrices
\begin{equation}\label{eq-NjoverB}
	\begin{pmatrix}
	N_j\\
	B
	\end{pmatrix}
\end{equation}
have ranks $\leq 6$, i.e
\begin{equation}\label{eq-condition-vanishingOfMinors}
	\mbox{all $7\times 7$ minors of (\ref{eq-NjoverB}) vanish.}
\end{equation}
 This gives a system of polynomial equations for the entries $b_{i,j}$. It turns out to be an effective  approach (indeed the unique one that I can carry out) to find $\Sigma$  by using its Plücker coordinates on the Grassmannian parametrizing the planes in $\mathbb{P}^6$. For each subset $I\subset[0,6]$ of cardinality 3, denote by $p_I$ the determinant of the submatrix of $B$ formed by $j$-th columns for $j\in C(I)=[0,6]\backslash I$. For each $j\in [0,6]$, the condition (\ref{eq-condition-vanishingOfMinors}) yields 4 linear equations for the Plücker coordinates $p_I$'s. 
 Due to lack of space, we present only a small part of the equations when $j=0$:

{\footnotesize
\begin{eqnarray*}
&&\begin{pmatrix}
\lambda_0 \lambda_1 \lambda_2 (\lambda_0-\lambda_1)(\lambda_1-\lambda_2)(\lambda_2-\lambda_0) & -\lambda_0 \lambda_1 \lambda_3 (\lambda_0-\lambda_1)(\lambda_1-\lambda_3)(\lambda_3-\lambda_0) &    \dots  \\
(\lambda_0 \lambda_1+\lambda_0 \lambda_2+\lambda_1 \lambda_2)(\lambda_0-\lambda_1)(\lambda_1-\lambda_2)(\lambda_2-\lambda_0) & 
-(\lambda_0 \lambda_1+\lambda_0 \lambda_3+\lambda_1 \lambda_3)(\lambda_0-\lambda_1)(\lambda_1-\lambda_3)(\lambda_3-\lambda_0) & \dots \\
(\lambda_0+\lambda_1+\lambda_2) (\lambda_0-\lambda_1)(\lambda_1-\lambda_2)(\lambda_2-\lambda_0) &  -(\lambda_0+\lambda_1+\lambda_3) (\lambda_0-\lambda_1)(\lambda_1-\lambda_3)(\lambda_3-\lambda_0) & \dots \\
(\lambda_0-\lambda_1)(\lambda_1-\lambda_2)(\lambda_2-\lambda_0) & -(\lambda_0-\lambda_1)(\lambda_1-\lambda_3)(\lambda_3-\lambda_0) &  \dots \\
\end{pmatrix}\nn\\
&&. (p_{012}, p_{013}, p_{023}, p_{123}, p_{014}, p_{024}, p_{124}, p_{034}, p_{134}, p_{234}, p_{015}, 
p_{025}, p_{125}, p_{035}, p_{135}, p_{235}, p_{045},\nn\\
&& p_{145}, p_{245}, p_{345}, p_{016}, p_{026}, 
p_{126}, p_{036}, p_{136}, p_{236}, p_{046}, p_{146}, p_{246}, p_{346}, p_{056}, p_{156}, p_{256}, 
p_{356}, p_{456})^{\mathrm{T}}=0.
\end{eqnarray*}
}
\begin{notation}
We denote such a system of 4 equations arising from (\ref{eq-condition-recursion-even(2,2)}) by $\mathrm{EC}_j$. For example, the exhibited above is $\mathrm{EC}_0$.  
\end{notation}

A priori, the system $\bigcup_{j=0}^6 \mathrm{EC}_j$ of  28 linear equations is an implication of, but not equivalent to, (\ref{eq-condition-vanishingOfMinors}). But as we will see in Section \ref{sec:conic-(1,2,3,4,5,6,7)}, the planes it defines  do satisfy (\ref{eq-condition-vanishingOfMinors}).

Perfectionists should be willing to solve this good-looking system of 28 linear equations in 35 unknowns, with $\lambda_0,\dots,\lambda_6$ as indeterminates.   I did not manage to achieve this.
 In the next section we choose a specific 7-tuple of parameters $\lambda_i$ and solve the corresponding curve counting problem. This is sufficient for computing the Gromov-Witten invariant (\ref{eq-enumerativeCorrelator-Dim4}).

\subsection{A case study}\label{sec:conic-(1,2,3,4,5,6,7)}
In this section we take $(\lambda_0,\dots,\lambda_6)=(1,2,3,4,5,6,7)$. We write a Macaulay2 package named \texttt{ConicsOn4DimIntersectionOfTwoQuadrics}. Within this package, running
\begin{equation*}
	\mbox{\texttt{syslinearEqsPluckerCoord}}\ \{1,2,3,4,5,6,7\}
\end{equation*}
 gives the whole system $\bigcup_{j=0}^6 \mathrm{EC}_j$.
  Solving  this  system of linear equations  of the Plücker coordinates, we get
\begin{equation}\label{eq-solveSyslinearEqsPluckerCoord}
\begin{cases}
{p}_{0,3,4} = \frac{40}{17}\,{p}_{0,1,2}-5\,{p}_{0,1,3}+\frac{82}{17}\,{p}_{0,2,3}-\frac{175}{102}\,{p}_{1,2,3}-\frac{75}{34}\,{p}_{0,1,4}+\frac{259}{408}\,{p}_{0,2,4}+\frac{55}{102}\,{p}_{1,2,4},\\ 
{p}_{1,3,4} = \frac{972}{17}\,{p}_{0,1,2}-129\,{p}_{0,1,3}+\frac{1792}{17}\,{p}_{0,2,3}-\frac{1544}{51}\,{p}_{1,2,3}+\frac{470}{17}\,{p}_{0,1,4}-\frac{3629}{102}\,{p}_{0,2,4}+\frac{698}{51}\,{p}_{1,2,4},\\
  {p}_{2,3,4} =  \frac{3658}{17}\,{p}_{0,1,2}-486\,{p}_{0,1,3}+\frac{6603}{17}\,{p}_{0,2,3}-\frac{1835}{17}\,{p}_{1,2,3}+\frac{2085}{17}\,{p}_{0,1,4}-\frac{9259}{68}\,{p}_{0,2,4}+\frac{771}{17}\,{p}_{1,2,4},\\ {p}_{0,1,5} =  \frac{126}{17}\,{p}_{0,1,2}-\frac{171}{10}\,{p}_{0,1,3}+\frac{1164}{85}\,{p}_{0,2,3}-\frac{317}{85}\,{p}_{1,2,3}+\frac{474}{85}\,{p}_{0,1,4}-\frac{883}{170}\,{p}_{0,2,4}+\frac{138}{85}\,{p}_{1,2,4},\\ {p}_{0,2,5} = \frac{105}{17}\,{p}_{0,1,2}-\frac{57}{4}\,{p}_{0,1,3}+\frac{194}{17}\,{p}_{0,2,3}-\frac{167}{51}\,{p}_{1,2,3}+\frac{79}{17}\,{p}_{0,1,4}-\frac{883}{204}\,{p}_{0,2,4}+\frac{86}{51}\,{p}_{1,2,4},\\ {p}_{1,2,5} =  -\frac{18}{17}\,{p}_{0,1,2}+\frac{9}{2}\,{p}_{0,1,3}-\frac{76}{17}\,{p}_{0,2,3}+\frac{11}{17}\,{p}_{1,2,3}-\frac{66}{17}\,{p}_{0,1,4}+\frac{135}{34}\,{p}_{0,2,4}-\frac{2}{17}\,{p}_{1,2,4},\\ {p}_{0,3,5} =  -\frac{4}{85}\,{p}_{0,1,2}+\frac{3}{5}\,{p}_{0,1,3}+\frac{36}{85}\,{p}_{0,2,3}-\frac{28}{51}\,{p}_{1,2,3}-\frac{46}{17}\,{p}_{0,1,4}+\frac{379}{510}\,{p}_{0,2,4}+\frac{146}{255}\,{p}_{1,2,4},\\ {p}_{1,3,5} = \frac{368}{17}\,{p}_{0,1,2}-45\,{p}_{0,1,3}+\frac{632}{17}\,{p}_{0,2,3}-\frac{395}{34}\,{p}_{1,2,3}+\frac{75}{34}\,{p}_{0,1,4}-\frac{993}{136}\,{p}_{0,2,4}+\frac{163}{34}\,{p}_{1,2,4},\\ {p}_{2,3,5} = \frac{1236}{17}\,{p}_{0,1,2}-162\,{p}_{0,1,3}+\frac{2204}{17}\,{p}_{0,2,3}-\frac{625}{17}\,{p}_{1,2,3}+\frac{690}{17}\,{p}_{0,1,4}-\frac{1467}{34}\,{p}_{0,2,4}+\frac{262}{17}\,{p}_{1,2,4},\\ {p}_{0,4,5} =  \frac{4}{85}\,{p}_{0,1,2}+\frac{9}{10}\,{p}_{0,1,3}-\frac{36}{85}\,{p}_{0,2,3}-\frac{2}{17}\,{p}_{1,2,3}-\frac{39}{17}\,{p}_{0,1,4}+\frac{191}{170}\,{p}_{0,2,4}+\frac{8}{85}\,{p}_{1,2,4},\\ {p}_{1,4,5} =  -\frac{90}{17}\,{p}_{0,1,2}+\frac{27}{2}\,{p}_{0,1,3}-\frac{108}{17}\,{p}_{0,2,3}+\frac{21}{17}\,{p}_{1,2,3}-\frac{126}{17}\,{p}_{0,1,4}+\frac{165}{34}\,{p}_{0,2,4}-\frac{10}{17}\,{p}_{1,2,4},\\ {p}_{2,4,5} = \frac{234}{17}\,{p}_{0,1,2}-\frac{117}{4}\,{p}_{0,1,3}+\frac{614}{17}\,{p}_{0,2,3}-\frac{160}{17}\,{p}_{1,2,3}+\frac{195}{17}\,{p}_{0,1,4}-\frac{705}{68}\,{p}_{0,2,4}+\frac{43}{17}\,{p}_{1,2,4},\\ {p}_{3,4,5} =  \frac{360}{17}\,{p}_{0,1,2}-45\,{p}_{0,1,3}+\frac{908}{17}\,{p}_{0,2,3}-\frac{220}{17}\,{p}_{1,2,3}+\frac{300}{17}\,{p}_{0,1,4}-\frac{262}{17}\,{p}_{0,2,4}+\frac{40}{17}\,{p}_{1,2,4},\\ {p}_{0,1,6} =  \frac{56}{17}\,{p}_{0,1,2}-\frac{38}{5}\,{p}_{0,1,3}+\frac{1552}{255}\,{p}_{0,2,3}-\frac{1268}{765}\,{p}_{1,2,3}+\frac{632}{255}\,{p}_{0,1,4}-\frac{1834}{765}\,{p}_{0,2,4}+\frac{569}{765}\,{p}_{1,2,4},\\ {p}_{0,2,6} = \frac{160}{51}\,{p}_{0,1,2}-\frac{43}{6}\,{p}_{0,1,3}+\frac{899}{153}\,{p}_{0,2,3}-\frac{3061}{1836}\,{p}_{1,2,3}+\frac{298}{153}\,{p}_{0,1,4}-\frac{3971}{1836}\,{p}_{0,2,4}+\frac{364}{459}\,{p}_{1,2,4},\\ {p}_{1,2,6} =  \frac{76}{85}\,{p}_{0,1,2}-\frac{9}{10}\,{p}_{0,1,3}+\frac{81}{85}\,{p}_{0,2,3}-\frac{33}{68}\,{p}_{1,2,3}-\frac{27}{17}\,{p}_{0,1,4}+\frac{237}{340}\,{p}_{0,2,4}+\frac{16}{85}\,{p}_{1,2,4},\\ {p}_{0,3,6} = \frac{74}{51}\,{p}_{0,1,2}-\frac{10}{3}\,{p}_{0,1,3}+\frac{433}{153}\,{p}_{0,2,3}-\frac{407}{459}\,{p}_{1,2,3}+\frac{151}{153}\,{p}_{0,1,4}-\frac{1915}{1836}\,{p}_{0,2,4}+\frac{227}{459}\,{p}_{1,2,4},\\ {p}_{1,3,6} =  \frac{468}{85}\,{p}_{0,1,2}-\frac{56}{5}\,{p}_{0,1,3}+\frac{2324}{255}\,{p}_{0,2,3}-\frac{440}{153}\,{p}_{1,2,3}+\frac{98}{51}\,{p}_{0,1,4}-\frac{3013}{1530}\,{p}_{0,2,4}+\frac{182}{153}\,{p}_{1,2,4},\\ {p}_{2,3,6} = \frac{344}{51}\,{p}_{0,1,2}-\frac{43}{3}\,{p}_{0,1,3}+\frac{1762}{153}\,{p}_{0,2,3}-\frac{3223}{918}\,{p}_{1,2,3}+\frac{445}{306}\,{p}_{0,1,4}-\frac{5381}{3672}\,{p}_{0,2,4}+\frac{943}{918}\,{p}_{1,2,4},\\ {p}_{0,4,6} =  \frac{64}{17}\,{p}_{0,1,2}-\frac{17}{2}\,{p}_{0,1,3}+\frac{121}{17}\,{p}_{0,2,3}-\frac{407}{204}\,{p}_{1,2,3}+\frac{42}{17}\,{p}_{0,1,4}-\frac{517}{204}\,{p}_{0,2,4}+\frac{44}{51}\,{p}_{1,2,4},\\ {p}_{1,4,6} = \frac{111}{17}\,{p}_{0,1,2}-15\,{p}_{0,1,3}+\frac{242}{17}\,{p}_{0,2,3}-\frac{407}{102}\,{p}_{1,2,3}+\frac{84}{17}\,{p}_{0,1,4}-\frac{517}{102}\,{p}_{0,2,4}+\frac{88}{51}\,{p}_{1,2,4},\\ {p}_{2,4,6}\ =\ \frac{24}{17}\,{p}_{0,1,2}-3\,{p}_{0,1,3}+\frac{124}{17}\,{p}_{0,2,3}-\frac{95}{51}\,{p}_{1,2,3}+\frac{20}{17}\,{p}_{0,1,4}-\frac{49}{51}\,{p}_{0,2,4}+\frac{8}{51}\,{p}_{1,2,4},\\ {p}_{3,4,6} = \frac{96}{17}\,{p}_{0,1,2}-12\,{p}_{0,1,3}+\frac{241}{17}\,{p}_{0,2,3}-\frac{176}{51}\,{p}_{1,2,3}+\frac{80}{17}\,{p}_{0,1,4}-\frac{196}{51}\,{p}_{0,2,4}+\frac{32}{51}\,{p}_{1,2,4},\\ {p}_{0,5,6} = \frac{158}{51}\,{p}_{0,1,2}-\frac{431}{60}\,{p}_{0,1,3}+\frac{4646}{765}\,{p}_{0,2,3}-\frac{3784}{2295}\,{p}_{1,2,3}+\frac{1771}{765}\,{p}_{0,1,4}-\frac{21803}{9180}\,{p}_{0,2,4}+\frac{1657}{2295}\,{p}_{1,2,4},\\ {p}_{1,5,6} = \frac{558}{85}\,{p}_{0,1,2}-\frac{153}{10}\,{p}_{0,1,3}+\frac{1132}{85}\,{p}_{0,2,3}-\frac{307}{85}\,{p}_{1,2,3}+\frac{414}{85}\,{p}_{0,1,4}-\frac{853}{170}\,{p}_{0,2,4}+\frac{26}{17}\,{p}_{1,2,4},\\ {p}_{2,5,6} =  -\frac{48}{17}\,{p}_{0,1,2}+6\,{p}_{0,1,3}-\frac{61}{17}\,{p}_{0,2,3}+\frac{199}{204}\,{p}_{1,2,3}-\frac{23}{17}\,{p}_{0,1,4}+\frac{239}{204}\,{p}_{0,2,4}-\frac{16}{51}\,{p}_{1,2,4},\\ {p}_{3,5,6} =  -\frac{184}{51}\,{p}_{0,1,2}+\frac{26}{3}\,{p}_{0,1,3}-\frac{1016}{153}\,{p}_{0,2,3}+\frac{808}{459}\,{p}_{1,2,3}-\frac{460}{153}\,{p}_{0,1,4}+\frac{923}{459}\,{p}_{0,2,4}-\frac{184}{459}\,{p}_{1,2,4},\\ {p}_{4,5,6} = -\frac{121}{17}\,{p}_{0,1,2}+\frac{69}{4}\,{p}_{0,1,3}-\frac{220}{17}\,{p}_{0,2,3}+\frac{353}{102}\,{p}_{1,2,3}-\frac{115}{17}\,{p}_{0,1,4}+\frac{821}{204}\,{p}_{0,2,4}-\frac{46}{51}\,{p}_{1,2,4}.
\end{cases}
\end{equation}
Pluging this solution into the Plücker equations, the Plücker ideal is transformed into the ideal (we have taken a reduced Gröbner basis as generators)
{\footnotesize
\begin{eqnarray*}
	&&(9\,{p}_{0,2,4}{p}_{1,2,4}-4\,{p}_{1,2,4}^{2},45\,{p}_{1,2,3}{p}_{1,2,4}+1755\,{p}_{0,1,4}{p}_{1,2,4
       }-254\,{p}_{1,2,4}^{2},180\,{p}_{0,2,3}{p}_{1,2,4}+1755\,{p}_{0,1,4}{p}_{1,2,4}-274\,{p}_{1,2,4}^{2},\nn\\
      && 135\,{p}_{0,1,3}{p}_{1
       ,2,4}-330\,{p}_{0,1,4}{p}_{1,2,4}+32\,{p}_{1,2,4}^{2},270\,{p}_{0,1,2}{p}_{1,2,4}-885\,{p}_{0,1,4}{p}_{1,2,4}+112\,{p}_{1,2
       ,4}^{2},81\,{p}_{0,2,4}^{2}-16\,{p}_{1,2,4}^{2},\nn\\
      && 9\,{p}_{0,1,4}{p}_{0,2,4}-4\,{p}_{0,1,4}{p}_{1,2,4},405\,{p}_{1,2,3}{p}_{0,
       2,4}+7020\,{p}_{0,1,4}{p}_{1,2,4}-1016\,{p}_{1,2,4}^{2},\nn\\
       && 405\,{p}_{0,2,3}{p}_{0,2,4}+1755\,{p}_{0,1,4}{p}_{1,2,4}-274\,{p}_{
       1,2,4}^{2}, 1215\,{p}_{0,1,3}{p}_{0,2,4}-1320\,{p}_{0,1,4}{p}_{1,2,4}+128\,{p}_{1,2,4}^{2},\nn\\
      &&  1215\,{p}_{0,1,2}{p}_{0,2,4}-1770
       \,{p}_{0,1,4}{p}_{1,2,4}+224\,{p}_{1,2,4}^{2},101475\,{p}_{0,1,4}^{2}-28320\,{p}_{0,1,4}{p}_{1,2,4}+1972\,{p}_{1,2,4}^{2},\nn\\
      && 101475\,{p}_{1,2,3}{p}_{0,1,4}+531710\,{p}_{0,1,4}{p}_{1,2,4}-76908\,{p}_{1,2,4}^{2},202950\,{p}_{0,2,3}{p}_{0,1,4}+243305
       \,{p}_{0,1,4}{p}_{1,2,4}-38454\,{p}_{1,2,4}^{2},\nn\\
      && 83025\,{p}_{0,1,3}{p}_{0,1,4}-36960\,{p}_{0,1,4}{p}_{1,2,4}+3944\,{p}_{1,2,
       4}^{2}, 913275\,{p}_{0,1,2}{p}_{0,1,4}-456600\,{p}_{0,1,4}{p}_{1,2,4}+58174\,{p}_{1,2,4}^{2},\nn\\
      && 913275\,{p}_{1,2,3}^{2}+
       14412060\,{p}_{0,1,4}{p}_{1,2,4}-2102008\,{p}_{1,2,4}^{2},913275\,{p}_{0,2,3}{p}_{1,2,3}+7560540\,{p}_{0,1,4}{p}_{1,2,4}-
       1098272\,{p}_{1,2,4}^{2},\nn\\
      && 9963\,{p}_{0,1,3}{p}_{1,2,3}+35508\,{p}_{0,1,4}{p}_{1,2,4}-5128\,{p}_{1,2,4}^{2},2739825\,{p}_{0,1
       ,2}{p}_{1,2,3}+2732055\,{p}_{0,1,4}{p}_{1,2,4}-391334\,{p}_{1,2,4}^{2},\\
      && 3653100\,{p}_{0,2,3}^{2}+11518065\,{p}_{0,1,4}{p}_{1
       ,2,4}-1716142\,{p}_{1,2,4}^{2},49815\,{p}_{0,1,3}{p}_{0,2,3}+30855\,{p}_{0,1,4}{p}_{1,2,4}-5098\,{p}_{1,2,4}^{2},\nn\\
      && 10959300
       \,{p}_{0,1,2}{p}_{0,2,3}-1259295\,{p}_{0,1,4}{p}_{1,2,4}+113786\,{p}_{1,2,4}^{2},83025\,{p}_{0,1,3}^{2}-42240\,{p}_{0,1,4}{
       p}_{1,2,4}+4976\,{p}_{1,2,4}^{2},\nn\\
       && 27675\,{p}_{0,1,2}{p}_{0,1,3}-12320\,{p}_{0,1,4}{p}_{1,2,4}+1588\,{p}_{1,2,4}^{2},913275
       \,{p}_{0,1,2}^{2}-254880\,{p}_{0,1,4}{p}_{1,2,4}+33533\,{p}_{1,2,4}^{2}).
\end{eqnarray*}
}
Its primary decomposition is
\begin{equation}\label{eq-primaryDecom-transformedPluckerIdeal}
\mathrm{ideal}(\bigcup_{j=0}^6 \mathrm{EC}_j)+\mbox{Plücker\ ideal}=	\mathfrak{p}_1 \cap \mathfrak{p}_2 \cap \mathfrak{p}_3,
\end{equation}
where 
\begin{eqnarray*}
\mathfrak{p}_1&=&(9\,{p}_{0,2,4}-4\,{p}_{1,2,4},15\,{p}_{0,1,4}-2\,{p}_{1,2,4},9\,{p}_{1,2,3}-4\,{p}_{1,2,4},9
       \,{p}_{0,2,3}-2\,{p}_{1,2,4},\nn\\
       && 45\,{p}_{0,1,3}-4\,{p}_{1,2,4},45\,{p}_{0,1,2}-{p}_{1,2,4}),\nn\\
 \mathfrak{p}_2&=&(9\,{p}_{0,2,4}-4\,{p}_{1,2,4},6765\,{p}_{0,1,4}-986\,{p}_{1,2,4},20295\,{p}_{1,2,3}+808\,{p}_{1,2,4},20295\,{p}_{0,2,3}-2053\,{p}_{1,2,4},\nn\\
 		&&369\,{p}_{0,1,3}-44\,{p}_{1,2,4}, 20295\,{p}_{0,1,2}-1277\,{p}_{1,2,4}),\nn\\
\mathfrak{p}_3&=& ({p}_{1,2,4},{p}_{0,1,4},{p}_{0,2,4}^{2},{p}_{1,2,3}{p}_{0,2,4},{p}_{0,2,3}{p}_{0,2,4},{p}_{0,1,3}{p}_{0,2,4},{p}_{0,1,2}{p}_{0,2,4},{
       p}_{1,2,3}^{2},{p}_{0,2,3}{p}_{1,2,3},\nn\\
       &&{p}_{0,1,3}{p}_{1,2,3},{p}_{0,1,2}{p}_{1,2,3},{p}_{0,2,3}^{2},{p}_{0,1,3}{p}_{0,2,3},{p}_{0,1,2}{p}_{0,2,3},{p}_{0,1,3}^{2},{p}_{0,1,2}{p}_{0,1,3},{p}_{0,1,2}^{2}).
\end{eqnarray*}
We thus get exactly two nonzero non-proportional solutions:
\begin{equation}\label{eq-intersectingPlane-example-2}
\begin{array}{|c|c|c|c|c|c|c|c|c|c|c|c|}
\hline
p_{012} &  p_{013} &  p_{023} &  p_{123} & p_{014} &  p_{024} &  p_{124} &  p_{034} &  p_{134} & p_{234} &  p_{015} & p_{025}\\
\hline
1 & 4 &  10 &  20 &  6 &  20 &  45 &  20 &  60 &  50 &  4 &  15 \\
\hline
p_{125} &  p_{035} & p_{135} &  p_{235} & p_{045} & p_{145} & p_{245} &  p_{345} & p_{016} &  p_{026} & p_{126} &  p_{036} \\
\hline
36 &  20 &  64 & 60 & 10 & 36 & 45 & 20 &  1 &  4 &  10 & 6 \\
\hline 
 p_{136} &  p_{236} & p_{046} & p_{146} &  p_{246} & p_{346} & p_{056} & p_{156} & p_{256} & p_{356} &  p_{456} &\\
\hline
 20 & 20 & 4 & 15 & 20 & 10 & 1 & 4 & 6 & 4 & 1 & \\
 \hline
\end{array}\
\end{equation}
and
\begin{equation}\label{eq-intersectingPlane-example-1}
\begin{array}{|c|c|c|c|c|c|c|c|c|c|c|c|}
\hline
p_{012} &  p_{013} &  p_{023} &  p_{123} & p_{014} &  p_{024} &  p_{124} &  p_{034} &  p_{134} & p_{234} &  p_{015} & p_{025}\\
\hline
1277 &  2420 &  2053 & -808 &  2958 & 9020 & 20295 &  12338 &  40332 &  38335 &  1804 & 8403\\
\hline
p_{125} &  p_{035} & p_{135} &  p_{235} & p_{045} & p_{145} & p_{245} &  p_{345} & p_{016} &  p_{026} & p_{126} &  p_{036} \\
\hline
21780 & 13024 & 42416 & 40332 & 6722 & 21780 & 20295 &  -808 & 451 & 2420 &  6722 & 3861 \\
\hline 
 p_{136} &  p_{236} & p_{046} & p_{146} &  p_{246} & p_{346} & p_{056} & p_{156} & p_{256} & p_{356} &  p_{456} &\\
\hline
13024 & 12338 & 2420 & 8403 & 9020 & 2053 &  451 &  1804 & 2958 & 2420 & 1277 & \\
\hline
\end{array}.
\end{equation}
(The equality or symmetry of certain coordinates arises from the special choice of $\lambda_i$; they are not true in general.)

 Both solutions have nonzero $p_{012}$. So the corresponding 2-plane is the solution (column) space of the matrix
\begin{equation*}
\begin{pmatrix}
\vspace{0.1cm}
	\frac{p_{123}}{p_{012}} & \frac{p_{023}}{p_{012}} & \frac{p_{013}}{p_{012}} & 1 & 0 & 0 & 0 \\
\vspace{0.1cm}	
	-\frac{p_{124}}{p_{012}} & -\frac{p_{024}}{p_{012}} & - \frac{p_{014}}{p_{012}} & 0 & 1 & 0 & 0 \\
\vspace{0.1cm}	
	\frac{p_{125}}{p_{012}} & \frac{p_{025}}{p_{012}} & \frac{p_{015}}{p_{012}} & 0 & 0 & 1 & 0 \\
	-\frac{p_{126}}{p_{012}} & -\frac{p_{026}}{p_{012}} & -\frac{p_{016}}{p_{012}} & 0 & 0 & 0 & 1
\end{pmatrix}.
\end{equation*}
Then one can check that the solution (\ref{eq-intersectingPlane-example-2}) gives the plane $S$. As we have seen in the proof of Proposition \ref{prop-noConicOnS}, the intersections of $S$ with $S_{i,i+1}$'s do not lie on a conic.  
The plane $\Sigma$ given by the solution (\ref{eq-intersectingPlane-example-1}) is defined by 
\begin{equation}\label{eq-definingMatrixOfSigma-matrix}
	\left(
\begin{array}{ccccccc}
\vspace{0.1cm}
 -\frac{808}{1277} & \frac{2053}{1277} & \frac{2420}{1277} & 1 & 0 & 0 & 0 \\
 \vspace{0.1cm}
 -\frac{20295}{1277} & -\frac{9020}{1277} & -\frac{2958}{1277} & 0 & 1 & 0 & 0 \\
 \vspace{0.1cm}
 \frac{21780}{1277} & \frac{8403}{1277} & \frac{1804}{1277} & 0 & 0 & 1 & 0 \\
 -\frac{6722}{1277} & -\frac{2420}{1277} & -\frac{451}{1277} & 0 & 0 & 0 & 1 \\
\end{array}
\right).
\end{equation}
It has the following intersecting points with $S_{i,i+1}$'s:
\begin{equation}\label{eq-intersectingPoints-example}
\begin{split}
&\Sigma\cap S_{01}=[1,-4,-\frac{90}{7},\frac{220}{7},-\frac{295}{7},\frac{192}{7},-\frac{48}{7}],\\
&\Sigma\cap S_{12}=[1,0,\frac{7}{2},-6,24,-22,\frac{13}{2}],\\
&\Sigma\cap S_{23}=[1,-\frac{36}{25},\frac{63}{50},\frac{14}{25},\frac{216}{25},-\frac{234}{25},\frac{149}{50}],\\
&\Sigma\cap S_{34}=[1,-\frac{468}{149},\frac{432}{149},\frac{28}{149},\frac{63}{149},-\frac{72}{149},\frac{50}{149}],\\
&\Sigma\cap S_{45}=[1,-\frac{44}{13},\frac{48}{13},-\frac{12}{13},\frac{7}{13},0,\frac{2}{13}],\\
&\Sigma\cap S_{56}=[1,-4,\frac{295}{48},-\frac{55}{12},\frac{15}{8},\frac{7}{12},-\frac{7}{48}],\\
&\Sigma\cap S_{60}=[1,-\frac{108}{7},\frac{495}{7},-\frac{760}{7},\frac{495}{7},-\frac{108}{7},1].
\end{split}
\end{equation}
We use the first three coordinates $W_0,W_1,W_2$ as coordinates on $\Sigma$. Then in $\Sigma$ one can check that  the first 5 points of (\ref{eq-intersectingPoints-example}) determine a unique conic 
\begin{equation}\label{eq-conic-(1,2,3,4,5,6,7)}
	C=\{W_0^2+\frac{231982 }{286839}W_0 W_1-\frac{69410 }{286839}W_0 W_2+\frac{924289 }{4589424}W_1^2-\frac{68035}{1721034} W_1 W_2-\frac{512}{40977} W_2^2=0\}
\end{equation}
and one checks that it passes through also the latter two points. We are left to show that $C$ lies in $X$. The intersection number $[C]$ and $[Q_i]$ is 4, but $\#(C\cap Q_i)\geq 7$, so $\dim (C\cap Q_i)\geq 1$, i.e. $C\subset Q_i$, for $i=1,2$. Thus $C\subset X$.\\

Let us find an explicit parametrization $\rho:\mathbb{P}^1\rightarrow C$. First, there is a point $(X,Y,Z)=(2,0,7)$ on $C$. Considering the intersection of $C$ with the line $\{W_2=7,W_1=t(W_0-2)\}$, we get a rational parametrization
\begin{equation}\label{eq-rationalParametrizationOfC-1}
\begin{split}
&W_0=w_0(t)=5545734 t^2+3809960 t-4214784,\\
 &W_1=w_1(t)=-18460312 t^2-31751328 t,\\
 & W_2=w_2(t)=19410069 t^2+77945952 t+96377904.
\end{split}
\end{equation}
Then using (\ref{eq-definingMatrixOfSigma-matrix}) we get
\begin{equation}\label{eq-rationalParametrizationOfC-2}
\begin{split}
&W_3=w_3(t)= -3596236 t^2-94256288 t-185309376,\\
&W_4=w_4(t)= 2704496 t^2+16828728 t+156262176,\\
&W_5=w_5(t)= -532380 t^2+33837888 t-64266048,\\
&W_6=w_6(t)= 1063751 t^2-12587344 t+11851728.
\end{split}
\end{equation}
By computing resultants one finds that $w_i(t)$ are pairwise coprime.
The defining equation of $Q_1$ is
\[
60 W_0^2-10 W_1^2+4 W_2^2-3 W_3^2+4 W_4^2-10 W_5^2+60 W_6^2.
\]
The defining equation of $Q_1$ is
\[
15 W_0^2-5 W_1^2+3 W_2^2-3 W_3^2+5 W_4^2-15 W_5^2+105 W_6^2.
\]
One can use also this parametrization to check that $C\subset X$.\\

We study the splitting type of $\Omega_X^1|_C$. We follow  \cite[Exercise 3.8]{Deb15}.
There is an exact sequence
\begin{equation*}
	0\rightarrow N_{C/X}\rightarrow N_{C/\mathbb{P}^{6}}\rightarrow N_{X/\mathbb{P}^{6}}\big|_C\rightarrow 0
\end{equation*}
and
\begin{equation*}
	N_{C/\mathbb{P}^{4+2}}\cong \mathcal{O}_{C}(4)\oplus \mathcal{O}_C(2)^{\oplus 4},\
	N_{X/\mathbb{P}^{4+2}}\big|_C\cong \mathcal{O}_C(4)^{\oplus 2}.
\end{equation*}
Twisting by $\mathcal{O}_C(-1)$ yields
\begin{equation*}
	0\rightarrow N_{C/X}(-1)\rightarrow N_{C/\mathbb{P}^{6}}(-1)
 	\rightarrow N_{X/\mathbb{P}^{6}}\big|_C(-1) \rightarrow 0.
\end{equation*}
Then $H^1(C, N_{C/X}(-1))=0$ iff
\begin{equation}\label{eq-mapExpectedSurj}
	H^0\big(C,N_{C/\mathbb{P}^{6}}(-1)\big)
 	\rightarrow H^0\big(C,N_{X/\mathbb{P}^{6}}\big|_C(-1)\big) 
\end{equation}
is surjective. Set
 \begin{equation*}
 	g_0(W_0,\dots,W_6)=W_0,\ g_1(W_0,\dots,W_6)=W_1,\ g_2(W_0,\dots,W_6)=W_2,
 \end{equation*}
 \begin{eqnarray*}	 
	g_3(W_0,\dots,W_6)&=& -\frac{808}{1277} W_0 + \frac{2053}{1277}W_1 + \frac{2420}{1277} W_2 + W_3, \\
	g_4(W_0,\dots,W_6)&=&  -\frac{20295}{1277}W_0 -\frac{9020}{1277}W_1 -\frac{2958}{1277}W_2 + W_4, \\
	g_5(W_0,\dots,W_6)&=&  \frac{21780}{1277}W_0+ \frac{8403}{1277}W_1+ \frac{1804}{1277}W_2+ W_5, \\
	g_6(W_0,\dots,W_6)&=&  -\frac{6722}{1277}W_0 -\frac{2420}{1277}W_1 -\frac{451}{1277}W_2+W_6.
\end{eqnarray*}
Recall that $\varphi_i$ is the defining equation (\ref{eq-definingEquationOfQi}) of $Q_i$. Then 
\begin{eqnarray*}
\frac{\partial}{\partial g_3}=\sum_{i=0}^{6}\frac{\partial W_i}{\partial g_3}\frac{\partial}{\partial W_i}=
\frac{\partial}{\partial W_3}
&\stackrel{(d\varphi_1,d\varphi_2)}{\mapsto} &  (\frac{\partial \varphi_1}{\partial W_3}\frac{\partial}{\partial \varphi_1},
\frac{\partial \varphi_2}{\partial W_3}\frac{\partial}{\partial \varphi_2})
= (2W_3\frac{\partial}{\partial \varphi_1},8W_3\frac{\partial}{\partial \varphi_2}).
\end{eqnarray*}
Let $[t:u]$ be the homogeneous coordinate on $C\cong \mathbb{P}^1$ such that $\frac{t}{u}$ is identified to the affine coordinate $t$ in the parametrization (\ref{eq-rationalParametrizationOfC-1}). Moreover let
\[
w_i(t,u)=u^2w_i(\frac{t}{u}),\ \mbox{for}\ 0\leq i\leq 6.
\]
 Then under the map (\ref{eq-mapExpectedSurj})
\[
t\frac{\partial}{\partial g_3}\mapsto \big(2t w_3(t,u)\frac{\partial}{\partial \varphi_1},8t w_3(t,u)\frac{\partial}{\partial \varphi_2}\big),\
u\frac{\partial}{\partial g_3}\mapsto \big(2u w_3(t,u)\frac{\partial}{\partial \varphi_1},8u w_3(t,u)\frac{\partial}{\partial \varphi_2}\big).
\]
Similarly,
\[
t\frac{\partial}{\partial g_4}\mapsto \big(2t w_4(t,u)\frac{\partial}{\partial \varphi_1},10t w_4(t,u)\frac{\partial}{\partial \varphi_2}\big),\
u\frac{\partial}{\partial g_4}\mapsto \big(2u w_4(t,u)\frac{\partial}{\partial \varphi_1},10u w_4(t,u)\frac{\partial}{\partial \varphi_2}\big),
\]
\[
t\frac{\partial}{\partial g_5}\mapsto \big(2t w_5(t,u)\frac{\partial}{\partial \varphi_1},12t w_5(t,u)\frac{\partial}{\partial \varphi_2}\big),\
u\frac{\partial}{\partial g_5}\mapsto \big(2u w_5(t,u)\frac{\partial}{\partial \varphi_1},12u w_5(t,u)\frac{\partial}{\partial \varphi_2}\big),
\]
\[
t\frac{\partial}{\partial g_6}\mapsto \big(2t w_6(t,u)\frac{\partial}{\partial \varphi_1},14t w_6(t,u)\frac{\partial}{\partial \varphi_2}\big),\
u\frac{\partial}{\partial g_6}\mapsto \big(2u w_6(t,u)\frac{\partial}{\partial \varphi_1},14u w_6(t,u)\frac{\partial}{\partial \varphi_2}\big).
\]
Then from (\ref{eq-rationalParametrizationOfC-1}) and (\ref{eq-rationalParametrizationOfC-2}) we get
{\footnotesize
\begin{eqnarray*}
\begin{pmatrix}
\vspace{0.1cm}
t\frac{\partial}{\partial g_3} \\
\vspace{0.1cm}
 u\frac{\partial}{\partial g_3} \\
 \vspace{0.1cm}
t\frac{\partial}{\partial g_4} \\
\vspace{0.1cm}
 u\frac{\partial}{\partial g_4} \\
 \vspace{0.1cm}
t\frac{\partial}{\partial g_5} \\
\vspace{0.1cm}
u\frac{\partial}{\partial g_5} \\
\vspace{0.1cm}
t\frac{\partial}{\partial g_6} \\
 u\frac{\partial}{\partial g_6} 
\end{pmatrix}\stackrel{(\ref{eq-mapExpectedSurj})}{\mapsto}
\left(
\begin{array}{cc}
 -7192472 t^3-188512576 t^2 u-370618752 t u^2, -28769888 t^3-754050304 t^2 u-1482475008 t u^2 \\
 -7192472 t^2 u-188512576 t u^2-370618752 u^3, -28769888 t^2 u-754050304 t u^2-1482475008 u^3 \\
 5408992 t^3+33657456 t^2 u+312524352 t u^2, 27044960 t^3+168287280 t^2 u+1562621760 t u^2 \\
 5408992 t^2 u+33657456 t u^2+312524352 u^3, 27044960 t^2 u+168287280 t u^2+1562621760 u^3 \\
 -1064760 t^3+67675776 t^2 u-128532096 t u^2, -6388560 t^3+406054656 t^2 u-771192576 t u^2 \\
 -1064760 t^2 u+67675776 t u^2-128532096 u^3, -6388560 t^2 u+406054656 t u^2-771192576 u^3 \\
 2127502 t^3-25174688 t^2 u+23703456 t u^2, 14892514 t^3-176222816 t^2 u+165924192 t u^2 \\
 2127502 t^2 u-25174688 t u^2+23703456 u^3, 14892514 t^2 u-176222816 t u^2+165924192 u^3 \\
\end{array}
\right).
\end{eqnarray*}
}
One checks that the matrix of coefficients
{\footnotesize
\[
\left(
\begin{array}{cccccccc}
 0 & -370618752 & -188512576 & -7192472 & 0 & -1482475008 & -754050304 & -28769888 \\
 -370618752 & -188512576 & -7192472 & 0 & -1482475008 & -754050304 & -28769888 & 0 \\
 0 & 312524352 & 33657456 & 5408992 & 0 & 1562621760 & 168287280 & 27044960 \\
 312524352 & 33657456 & 5408992 & 0 & 1562621760 & 168287280 & 27044960 & 0 \\
 0 & -128532096 & 67675776 & -1064760 & 0 & -771192576 & 406054656 & -6388560 \\
 -128532096 & 67675776 & -1064760 & 0 & -771192576 & 406054656 & -6388560 & 0 \\
 0 & 23703456 & -25174688 & 2127502 & 0 & 165924192 & -176222816 & 14892514 \\
 23703456 & -25174688 & 2127502 & 0 & 165924192 & -176222816 & 14892514 & 0 \\
\end{array}
\right)
\]}
is nonsingular. So  (\ref{eq-mapExpectedSurj}) is surjective and thus $H^1(C, N_{C/X}(-1))=0$, i.e. $C$ is a \emph{free curve} in $X$. It follows that $X$ is covered by conics.\\

As we said below Definition \ref{def-enumerativeCorrelator}, to compute the correlator (\ref{eq-possiblyEnumerativeCorrelator}) enumeratively, we need not to only count curves, but also count certain multiplicities. As we will see in the next section, this forces us to  solve the system 
\begin{equation}\label{eq-sysOfLinearAndPluckerEquations}
\bigcup_{j=0}^6 \mathrm{EC}_j+\mbox{Plücker equations}
\end{equation}
 over the ring of  dual numbers $\mathbb{C}[\varepsilon]/(\varepsilon^2)$. Let $\tilde{p}_{i,j,k}=p_{i,j,k}+x_{i,j,k} \varepsilon$. We need only find the solutions with $p_{i,j,k}$ equal to those given by (\ref{eq-intersectingPlane-example-1}). So we replace the variables in the system (\ref{eq-sysOfLinearAndPluckerEquations}) by $\tilde{p}_{i,j,k}$'s, and plug the values of $p_{i,j,k}$ from (\ref{eq-intersectingPlane-example-1}), and  obtain
\begin{eqnarray}\label{eq-solveSysOfLinearAndPluckerEquations-withEpsilon}
\vspace{0.2cm}
&&\underbrace{\mathrm{ideal}(\bigcup_{j=0}^6 \mathrm{EC}_j)+\mbox{Plücker\ ideal}}_{\mbox{with variables replaced by}\ \tilde{p}_{i,j,k}}
+\sum_{i,j,k}\mathrm{ideal}\big(p_{i,j,k}-\mbox{the value of } p_{i,j,k}\mbox{ in (\ref{eq-intersectingPlane-example-1})}\big)\nn\\
	&=& (9x_{0,2,4}\varepsilon-4x_{1,2,4}\varepsilon,\ 6765x_{0,1,4}\varepsilon-986x_{1,2,4}\varepsilon,\ 20295x_{1,2,3}\varepsilon+808x_{1,2,4}\varepsilon,\nn\\
      &&  20295x_{0,2,3}\varepsilon-2053x_{1,2,4}\varepsilon,\ 369x_{0,1,3}\varepsilon-44x_{1,2,4}\varepsilon,\ 20295x_{0,1,2}\varepsilon-1277x_{1,2,4}\varepsilon)
\end{eqnarray}
Compare (\ref{eq-solveSysOfLinearAndPluckerEquations-withEpsilon}) to $\mathfrak{p}_1$ in (\ref{eq-primaryDecom-transformedPluckerIdeal}), we find that the solution of (\ref{eq-sysOfLinearAndPluckerEquations}) that deforms the solution (\ref{eq-intersectingPlane-example-1}) are of the form 
\begin{equation*}
	\tilde{p}_{i,j,k}=u p_{i,j,k}
\end{equation*}
where $u$ is a unit of $\mathbb{C}[\varepsilon]/(\varepsilon^2)$, and $p_{i,j,k}$ is given by (\ref{eq-intersectingPlane-example-1}). One can also understand this result from the multiplicity 1 of the primes in the decomposition (\ref{eq-primaryDecom-transformedPluckerIdeal}).
\\

We summarize the computations in this section as follows.
\begin{theorem}\label{thm-conic-(1,2,3,4,5,6,7)}
Let $X$ be the 4 dimensional smooth complete intersection of two quadrics, given by (\ref{eq-definingEquationOfQi}) with 
$(\lambda_0,\dots,\lambda_6)=(1,2,3,4,5,6,7)$. Then
\begin{enumerate}
	\item[(i)] There exists a unique conic $C$ in $X$ that meets $S_{i,i+1}$ for $i\in [0,6]$.
	\item[(ii)] The conic $C$ is a free curve in $X$.
	\item[(iii)] In the ring of dual numbers $\mathbb{C}[\varepsilon]/(\varepsilon^2)$, up to a common multiple, the system (\ref{eq-sysOfLinearAndPluckerEquations}) has only one solution whose reduction modulo $(\varepsilon)$ equals (\ref{eq-intersectingPlane-example-1}).
\end{enumerate}
\end{theorem}

\subsection{The special correlator in dimension 4}\label{sec:specialCorrelator-Dim4}
In this section we show how to deduce (\ref{eq-enumerativeCorrelator-Dim4}) from the computations in Section \ref{sec:conic-(1,2,3,4,5,6,7)}. The following lemma is slightly weaker than Theorem \ref{thm-conic-(1,2,3,4,5,6,7)} (ii)  in the case studied in Section \ref{sec:conic-(1,2,3,4,5,6,7)}.
\begin{lemma}\label{lem-weakFreeness}
Let $X$ be a 4-dimensional (2,2)-complete intersection.
Suppose that $C$ is a conic in $X$ that meets each of $S_{i,i+1}$ for $0\leq i\leq 6$. Then 
$\Mbar_{0,0}(X,2)$ is smooth at the point  $[C\hookrightarrow X]$ and its dimension at this point is the expected dimension 7.
\end{lemma}
\begin{proof}
It suffices to show $H^1(C,N_{C/X})=0$.
Let $\Sigma$ be the plane spanned by $C$. Since $N_{C/\Sigma}\cong \mathcal{O}_C(4)$ and $N_{\Sigma/\mathbb{P}^{6}}\big|_C\cong \mathcal{O}_C(2)^{\oplus 4}$, from the exact sequence
\begin{equation*}
	0\rightarrow N_{C/\Sigma}\rightarrow N_{C/\mathbb{P}^{6}}\rightarrow N_{\Sigma/\mathbb{P}^{6}}\big|_C\rightarrow 0,
\end{equation*}
it follows that there is a natural splitting 
\begin{equation*}
	N_{C/\mathbb{P}^{6}}=N_{C/\Sigma}\oplus N_{\Sigma/\mathbb{P}^{6}}\big|_C.
\end{equation*}
Thus 
 $N_{C/\mathbb{P}^{6}}\cong \mathcal{O}_{C}(4)\oplus \mathcal{O}_C(2)^{\oplus 4}$. 
Moreover, one has $N_{X/\mathbb{P}^{6}}|_C\cong \mathcal{O}_C(4)^{\oplus 2}$. Suppose 
$ N_{C/X}=\bigoplus_{i=1}^{3}\mathcal{O}_C(b_i)$. Then $b_1+b_2+b_3=4$. 
From the exact sequence
\begin{equation*}
	0\rightarrow N_{C/X}\rightarrow N_{C/\mathbb{P}^{6}}\rightarrow N_{X/\mathbb{P}^{6}}\big|_C\rightarrow 0
\end{equation*}
one finds that the only possibility that $H^1(C,N_{C/X})\neq 0$ is $\{b_1,b_2,b_3\}=\{4,2,-2\}$.
Note that the summand $\mathcal{O}_C(4)$ in $N_{C/\mathbb{P}^6}$ is identified with $N_{C/\Sigma}$. It follows that if $b_1=4$, $X$ is tangent to $\Sigma$ everywhere along $C$; more precisely, $T_X \cap T_{\Sigma}\neq 0$ at every point of $C$. 
Write $X=Q_1\cap Q_2$, where $Q_1$ and $Q_2$ are smooth quadric hypersurfaces in $\mathbb{P}^6$.
As we have seen, $\Sigma\not\subset X$. So without loss of generality we can assume $\Sigma\not\subset Q_1$. So if $b_1=4$, then $Q_1$ meets $\Sigma$ properly, and  the intersection multiplicity of $Q_1$ and $\Sigma$ is at least $2$ at the generic point of $C$. This contradicts that, in $\mathrm{CH}(\mathbb{P}^6)$,
\[
[\Sigma]\cap [Q_1]=[C].
\]
So $\{b_1,b_2,b_3\}\neq \{4,2,-2\}$, and $H^1(C,N_{C/X})=0$.
\end{proof}

\begin{theorem}\label{thm-unknownCorrelator-sij-Even(2,2)-4dim}
For any  $4$-dimensional complete intersections of two quadrics in $\mathbb{P}^{6}$, 
\begin{equation}\label{eq-unknownCorrelator-sij-Even(2,2)-4dim}
	\langle \varsigma_{0,1},\dots,\varsigma_{6,0}\rangle_{0,7,2}=1.
\end{equation}
\end{theorem}
\begin{proof}
It suffices to show (\ref{eq-unknownCorrelator-sij-Even(2,2)-4dim}) for a single $4$-dimensional complete intersection of two quadrics in $\mathbb{P}^{6}$. We take the one in Section \ref{sec:conic-(1,2,3,4,5,6,7)}. By the final computation in Section \ref{sec:conic-(1,2,3,4,5,6,7)}, or by Lemma \ref{lem-weakFreeness}, $\Mbar_{0,0}(X,2)$ is smooth at the point representing $C\hookrightarrow X$. Let $p_i=C\cap S_{i,i+1}$. Then $\Mbar_{0,7}(X,2)$ is smooth at the point representing the stable map $(C,{p_0,\dots,p_6})\hookrightarrow X$. Let $M$ be the irreducible component of $\Mbar_{0,7}(X,2)$ that containing $(C,{p_0,\dots,p_6})$. Then $M$ has the expected dimension, and 
\[
[\Mbar_{0,7}(X,2)]=[M]+\sum_i [M_i]
\]
where $M_i$ are cycles lying on other components that does not containing $(C,{p_0,\dots,p_6})$. 
Denote the $i$-th evaluation map by $\mathrm{ev}_i: \Mbar_{0,7}(X,2)\rightarrow X$, for $i=0,\dots,6$. Denote the $i$-th projection $X^7\rightarrow X$ by $q_i$.

 By the definition of Gromov-Witten invariants, one has
\begin{equation*}
	\langle \varsigma_{0,1},\dots,\varsigma_{6,0}\rangle_{0,7,2}=
	(\mathrm{ev}_0\times\cdots\times \mathrm{ev}_6)_*[\Mbar_{0,7}(X,2)]\cap \bigcup_{i=0}^6 q_i^* \varsigma_{i,i+1}.
\end{equation*}
The computational results in Section \ref{sec:conic-(1,2,3,4,5,6,7)} show, that the cycle $M$ has a proper intersection with the cycle $\bigcap_{i=0}^6 q_i^* S_{i,i+1}$, at exactly one point, i.e. $(p_0,\dots,p_6)\in X^7$, and the other cycles $M_i$ do not meet $\bigcap_{i=0}^6 q_i^* S_{i,i+1}$. To show (\ref{eq-unknownCorrelator-sij-Even(2,2)-4dim}), we are left to show that the intersection multiplicity of $M$ and $\bigcap_{i=0}^6 q_i^* S_{i,i+1}$ is 1. By \cite[Proposition 7.1(a)]{Ful98}, it suffices to show that the scheme theoretic intersection 
\begin{equation*}
\Spec(A):=
 	 \bigcap_{i=0}^6 (\mathrm{ev}_0\times\cdots\times \mathrm{ev}_6)^*q_i^* S_{i,i+1}=\bigcap_{i=0}^6 \mathrm{ev}_i^* S_{i,i+1}
\end{equation*}
has length 1 at $(C,p_0,\dots,p_6)$. Suppose, on the contrary, $\mathrm{length}(A)\geq 2$. Then there is a closed immersion $\Spec\ \mathbb{C}[\epsilon]/(\epsilon^2)\hookrightarrow \Spec(A)$. The composition $\iota:\Spec\ \mathbb{C}[\epsilon]\hookrightarrow \Mbar_{0,7}(X,2)$ induces a \emph{nontrivial} family of stable maps over  $\Spec\ \mathbb{C}[\epsilon]/(\epsilon^2)$, with the $i$-th  marked point lying in $S_{i,i+1}$, for $i=0,\dots,6$. More precisely, in the graph,
\begin{equation*}
	\xymatrix{
		\mathcal{C} \ar[r]^{f} \ar[d]_{\pi} & X \\
		\Spec \mathbb{C}[\epsilon]/(\epsilon^2) \ar@/_/[u]_{\sigma_i} & 
	}
\end{equation*}
the composition $f\circ \sigma_i$ is a closed immersion into $S_{i,i+1}$. We are going to show that such a nontrivial family does not exist. Indeed, suppose given such a family. Then
\[
H^0\big(\mathcal{C},\mathcal{O}_{\mathbb{P}^{n+2}_{\mathbb{C}[\epsilon]/(\epsilon^2)}}(1)\big)
\]
is a rank 2 free module over $\mathbb{C}[\epsilon]/(\epsilon^2)$. 
So it induces a solution of the system (\ref{eq-sysOfLinearAndPluckerEquations}) over $\mathbb{C}[\epsilon]/(\epsilon^2)$. By Theorem \ref{thm-conic-(1,2,3,4,5,6,7)} (iii) such a solution  is unique and is no other than the trivial deformation of the plane $\Sigma$. 
Then we are left to show that the system of equations for   the curve $C$ on $\Sigma$ to pass through the points $\Sigma\cap S_{i,i+1}$ for $i=0,\dots,4$ also has only trivial deformations over $\mathbb{C}[\varepsilon]/(\varepsilon^2)$. One can see this  by solving equations in the ring of dual numbers as the last part of Section \ref{sec:conic-(1,2,3,4,5,6,7)}. 
Here we  deduce it as an enumerative consequence of the Gromov-Witten invariant of $\mathbb{P}^2$. The moduli stack $\Mbar_{0,5}(\mathbb{P}^2,2)$ is smooth and irreducible and has a dense open smooth subscheme, and we know that 5 general points on $\mathbb{P}^2$ determine a unique conic. It follows that
\begin{equation}\label{eq-GWInvariant-P2}
\langle pt,\dots,pt\rangle_{0,5,2}^{\mathbb{P}^2}=1.
\end{equation}
Then given 5 points on $\mathbb{P}^2$, denote  by (E5) the system of equations of coefficients of conics $C$ for it to pass through these  5 points. Suppose 
\begin{equation*}
  	\mbox{(E5)  has only finitely many solutions}.
  \end{equation*} 
 Then (E5) has exactly one solution in $\mathbb{C}$, and has exactly one solution  also in $\mathbb{C}[\varepsilon]/(\varepsilon^2)$ which is a trivial deformation of its  solution in $\mathbb{C}$; this is a consequence of (\ref{eq-GWInvariant-P2}) by   the reasoning  as above. As we said around (\ref{eq-conic-(1,2,3,4,5,6,7)}), the condition (E5) holds for the first 5 points of (\ref{eq-intersectingPoints-example}). So we are done. 
\end{proof}

\begin{corollary}\label{cor-unknownCorrelator-Even(2,2)-4dim}
Let $\varepsilon_1,\dots,\varepsilon_{n+3}$ be the basis of $H^{n}_{\mathrm{prim}}(X)$ defined as in Section \ref{sec:explictD-Lattice}, and  $\epsilon_1,\dots,\epsilon_{n+3}$ be defined as in (\ref{eq-normalizedOrthonormalBasis}). 
Then for any  $4$-dimensional complete intersections of two quadrics in $\mathbb{P}^{6}$, 
\begin{equation}\label{eq-unknownCorrelator-Even(2,2)-4dim}
	\langle \epsilon_{1},\dots,\epsilon_{7}\rangle_{0,7,2}=\frac{1}{2}.
\end{equation}
\end{corollary}
\begin{proof}
This follows from  Lemma \ref{lem-relate-specialCorrelatorToEnumerativeCorrelator} and Theorem \ref{thm-unknownCorrelator-sij-Even(2,2)-4dim}.
\end{proof}

\begin{theorem}\label{thm-countingConics-4dim-general}
For general 4-dimensional  smooth complete intersections $X$ of two quadrics in $\mathbb{P}^6$, there exists exactly one smooth conic that meets each of the 2-planes $S_{i,i+1}$ in $X$ for $0\leq i\leq 6$.
\end{theorem}
\begin{proof}
Let $\mathscr{X}\rightarrow \mathscr{T}$ be the family of all 4-dimensional  smooth complete intersections of two quadrics in $\mathbb{P}^6$. We replace $\mathscr{T}$ by a finite Galois cover, such that the base change of this family to the cover has a family of 2-planes $\mathcal{S}_{i,i+1}$ for $0\leq i\leq 6$. We still denote the base changed family by $\mathscr{X}\rightarrow \mathscr{T}$. The moduli stacks of stable maps of degree 2 to the fibers of $\mathscr{X}\rightarrow \mathscr{T}$ form a family, which we denote by  $\overline{\mathscr{M}}_{0,k}(\mathscr{X}/\mathscr{T},2)$, over $\mathscr{T}$. Consider the fiber product
\begin{equation*}
	\mathscr{Z}:=\bigcap_{i=0}^6 \mathrm{ev}_i^* \mathcal{S}_{i,i+1}.
\end{equation*}
By Theorem \ref{thm-unknownCorrelator-sij-Even(2,2)-4dim}, there exists a closed point $0\in \mathscr{T}$ such that the fiber of $\mathscr{Z}$ over $0$ consists of exactly one point. Then by Chevalley's semi-continuity theorem, there exists an open neighborhood $U$ of $0$ in $\mathscr{T}$ such that $\mathscr{Z}_{U}$ is finite over $U$. 
By Lemma \ref{lem-weakFreeness}, 
 if a point in $\mathscr{Z}_v$ represents a smooth conic in $\mathscr{X}_v$,  it contributes at least 1 to the Gromov-Witten invariant (\ref{eq-unknownCorrelator-sij-Even(2,2)-4dim}). Since this invariant equals 1, there exists at most one point in $\mathscr{Z}_v$ representing a smooth conic in $\mathscr{X}_v$.

We are going to show that, there exists a non-empty open subset $V_0$ of $V$ such that for all $v\in V_0$,   there is a point in $\mathscr{Z}_v$ representing a smooth conic in $\mathscr{X}_v$ for $v\in V$. It suffices to show that, shrinking $V$ if necessary, there is no point in $\mathscr{Z}_v$ that represents a singular conic (i.e. a union of two lines) or a double cover of a line in $X$. First we note that this does not happen for the example $(\lambda_0,\dots,\lambda_6)=(1,2,3,4,5,6,7)$. Indeed, our computation in Section \ref{sec:conic-(1,2,3,4,5,6,7)} in finding conics covers the cases of  singular conics or double lines, and the result shows that no singular conics or double lines is meeting $S_{i,i+1}$ for $0\leq i\leq 6$. Then it suffices to note that, the existence of singular conics or double lines meeting $S_{i,i+1}$ for $0\leq i\leq 6$ is equivalent to the existence of solutions to a certain system of algebraic equations of $\lambda_0,\dots,\lambda_6$. So the locus of such existence is a closed subscheme of $V$. Hence we are done.
\end{proof}

For general $X$ that satisfies the conclusion of Theorem \ref{thm-countingConics-4dim-general}, the conic $C$ depends on the choice of a $2$-plane $S$ in $X$, and also on the choice of the order of the coordinates $Y_i$. Nevertheless, such  conics $C$ are rather distinctive. Therefore we wish for explicit descriptions of $C$ in terms of $(\lambda_0,\dots,\lambda_6)$. As we commented at the end of Section \ref{sec:enumerativeCorrelator}, it  seems not plausible to do this by solving the system of  equations (\ref{eq-sysOfLinearAndPluckerEquations}) in terms of $\lambda_0,\dots,\lambda_6$.
On the other hand, it is natural to expect that the family of smooth conic $C$ degenerates on $X$ which is not so general to satisfy the conclusion of Theorem \ref{thm-countingConics-4dim-general}. We propose thus the following problems.

\begin{problem}\label{problem-explicit-C}
Describe explicitly the conic $C$ for general $\lambda_0,\dots,\lambda_6$.
\end{problem}

\begin{question}\label{ques-degeneration-C}
Is the statement of Theorem \ref{thm-countingConics-4dim-general} true in an appropriate sense (e.g. allowing singular conics or double lines), for \emph{all} 4-dimensional smooth complete intersections of two quadrics in $\mathbb{P}^6$?
\end{question}

Clearly a solution to Problem \ref{problem-explicit-C} will be helpful to Question \ref{ques-degeneration-C}. 
For Problem \ref{problem-explicit-C}, we have only a partial conjectural statement.
\begin{conjecture}
In the $W$-coordinates (\ref{eq-Wcoordinates}) we define a quadric $Q\subset \mathbb{P}^{n+2}$ as follows. Define
\begin{eqnarray*}
h(\lambda_0,\dots,\lambda_6)&:=&\lambda_0^2 \lambda_1 \lambda_3-\lambda_0^2 \lambda_1 \lambda_5-\lambda_0^2 \lambda_2 \lambda_3+\lambda_0^2 \lambda_2 \lambda_6+\lambda_0^2 \lambda_4 \lambda_5-\lambda_0^2 \lambda_4 \lambda_6+\lambda_0 \lambda_1 \lambda_2 \lambda_5\\
&&-\lambda_0 \lambda_1 \lambda_2 \lambda_6-\lambda_0 \lambda_1 \lambda_3 \lambda_4-\lambda_0 \lambda_1 \lambda_3 \lambda_6+\lambda_0 \lambda_1 \lambda_4 \lambda_6+\lambda_0 \lambda_1 \lambda_5 \lambda_6+\lambda_0 \lambda_2 \lambda_3 \lambda_4\\
&&+\lambda_0 \lambda_2 \lambda_3 \lambda_5-\lambda_0 \lambda_2 \lambda_4 \lambda_5-\lambda_0 \lambda_2 \lambda_5 \lambda_6-\lambda_0 \lambda_3 \lambda_4 \lambda_5+\lambda_0 \lambda_3 \lambda_4 \lambda_6-\lambda_1 \lambda_2 \lambda_3 \lambda_5\\
&&+\lambda_1 \lambda_2 \lambda_3 \lambda_6+\lambda_1 \lambda_3 \lambda_4 \lambda_5-\lambda_1 \lambda_4 \lambda_5 \lambda_6-\lambda_2 \lambda_3 \lambda_4 \lambda_6+\lambda_2 \lambda_4 \lambda_5 \lambda_6,
\end{eqnarray*}
and
\[
\mu_i(\lambda_0,\dots,\lambda_6):=h(\lambda_i,\lambda_{i+1},\lambda_{i+2},\lambda_{i+3},\lambda_{i+4},\lambda_{i+5},\lambda_{i+6})\cdot \prod_{j=i+1}^{i+6}(\lambda_i-\lambda_j)
\]
for $0\leq i\leq 6$, where the subscripts are understood in the mod 7 sense. We define a quadric hypersurface $Q$  by
\begin{equation*}
	\sum_{i=0}^6 \mu_i(\lambda_0,\dots,\lambda_6) W_i^2=0.
\end{equation*}
Then the $2$-plane $S_C$ spanned by the conic $C$ is contained in $Q$.
\end{conjecture}

\begin{appendices}
\section{An algorithm}\label{sec:algorithm}
In this Appendix, we describe an algorithm based on the proof of Theorem \ref{thm-reconstruction-even(2,2)} to compute the primary genus 0 Gromov-Witten invariants, with the special correlator (\ref{eq-specialLength(n+3)Invariant-even(2,2)}) as an unknown. 
First, we rewrite (\ref{eq-recursion-primitive-aabb-even(2,2)}) as
\begin{eqnarray}\label{eq-recursion-primitive-aabb-tau-simplified-even(2,2)}
(\frac{2|I|-4}{n-1}-2i_b)\partial_{\tau^{a}}^2\partial_{\tau^I}F(0)
=-(\frac{2|I|-4}{n-1}-2i_a)\partial_{\tau^{b}}^2\partial_{\tau^I}F(0)
+\mbox{RHS of  (\ref{eq-recursion-primitive-aabb-even(2,2)})}.
\end{eqnarray}

\begin{algorithm}\label{algorithm-correlator-even(2,2)}
For $I=(i_0,\dots,i_{2n+3})\in \mathbb{Z}_{\geq 0}^{2n+4}$, we compute the correlator $\partial_{\tau^I}F(0)$ recursively as follows.
\begin{enumerate}
	\item  Define $\beta:=\frac{\sum_{k=0}^n k i_k+\frac{n}{2}\sum_{k=n+1}^{2n+3}i_k-(n-3+|I|)}{n-1}$. If $\beta\not\in \mathbb{Z}$,  return 0. In the following steps we assume $\beta\in \mathbb{Z}$. This excludes for example the case that $\sum_{k=0}^{n}i_k=0$ and $\sum_{k=n+1}^{2n+3}i_k=2$.  (By (\ref{eq-Dim})).
	\item If $\sum_{k=n+1}^{2n+3}i_k=0$, returns the ambient correlator $\partial_{\tau^I}F(0)$.
	\item If $i_0>0$ and $\sum_{k=0}^{2n+3}i_k>3$ then return 0. (By (\ref{eq-FCA}))
	\item If $i_1>0$ and $\sum_{k=0}^{2n+3}i_k>3$ then apply (\ref{eq-recursion-EulerVecField-even(2,2)}).
	\item If $\sum_{k=2}^{n}i_k>0$ and $\sum_{k=0}^{2n+3}i_k\geq 2$ then apply 
	(\ref{eq-recursion-ambient-simplified-even(2,2)}).
	\item If $\sum_{k=n+1}^{2n+3}i_k=1$,  return 0. (By Theorem \ref{thm-monodromy-evenDim(2,2)})
	\item If $\sum_{k=n+1}^{2n+3}i_k=3$, return 0. (By Theorem \ref{thm-monodromy-evenDim(2,2)})
	\item If $\sum_{k=0}^{n}i_k=0$ and $\sum_{k=n+1}^{2n+3}i_k=4$, then if the nonzero components of $(i_{n+1},\dots,i_{2n+3})$ is a 4, or two 2's, then return 1, otherwise return 0. (By Theorem \ref{thm-monodromy-evenDim(2,2)} and  Theorem \ref{thm-reconstruction-even(2,2)})
	\item If $\sum_{k=0}^{n}i_k=0$ and $i_k=1$ for all $n+1\leq k\leq 2n+3$, return an indeterminate \texttt{x}, which stands for the unknown special correlator (\ref{eq-specialLength(n+3)Invariant-even(2,2)}).
	\item If $\sum_{k=0}^{n}i_k=0$, $\sum_{k=n+1}^{2n+3}i_k>4$ and there is only one nonzero component $i_a$ in $(i_{n+1},\dots,i_{2n+3})$, we rearrange $I$ such that $I=(0,\dots,0,i_{2n+3})$, and take $b=n+1$ so that $i_b=0$. Then apply (\ref{eq-recursion-primitive-aabb-tau-simplified-even(2,2)}) to $a=2n+3$, $b=n+1$, and $I=I-2e_a$, and thus return the RHS of (\ref{eq-recursion-primitive-aabb-tau-simplified-even(2,2)}) divided by a nonzero number.
	\item We define a function \texttt{indexTriple}. The input is $(i_{n+1},\dots,i_{2n+3})$. Then we check that neither of the conditions
	\begin{enumerate}
	 	\item[(i)] there is only one nonzero component in $(i_{n+1},\dots,i_{2n+3})$
	 	\item[(ii)] $i_k=1$ for all $n+1\leq k\leq 2n+3$
	 \end{enumerate}
	is true. When this is the case,   the output is the indices $a,b,c$ as given in Remark \ref{rem:recursion-boundOfIndex}. This produces in a definite way the indices $a,b,c$ in the proof of Theorem \ref{thm-reconstruction-even(2,2)}.
	\item If $\sum_{k=0}^{n}i_k=0$, $\sum_{k=n+1}^{2n+3}i_k>4$, and the neither of the conditions (i) and (ii) in the last step is true, we apply \texttt{indexTriple} to $(i_{n+1},\dots,i_{2n+3})$ to get $(a,b,c)$. Then we apply (\ref{eq-recursion-primitive-abcc-even(2,2)}) to $(a,b,c)$ and $I=I-e_a-e_b$, and thus return the RHS of (\ref{eq-recursion-primitive-abcc-even(2,2)}) divided by a nonzero number. (As we have  shown in the proof of Theorem \ref{thm-reconstruction-even(2,2)} that  the coefficient on the LHS of (\ref{eq-recursion-primitive-abcc-even(2,2)}) is nonzero in this case.)
\end{enumerate}

\end{algorithm}
We implemented Algorithm \ref{algorithm-correlator-even(2,2)} by the command 
\[
\mbox{\texttt{correlatorEvenIntersecTwoQuadricsRecursionInTauCoord}}
\]
 in the Macaulay2 package \texttt{QuantumCohomologyFanoCompleteIntersection}. The input of this command is a list $\{n,I\}$, and the output is $\partial_{\tau^I}F(0)$. For example, running
\[
\mbox{\texttt{correlatorEvenIntersecTwoQuadricsRecursionInTauCoord} $\{4,\{0,0,7,0,0,0,0,0,0,0,0,0\}\}$}
\]
returns 46656, so $\langle\sfh_2,\sfh_2,\sfh_2,\sfh_2,\sfh_2,\sfh_2,\sfh_2\rangle=46656$. Running
\[
\mbox{\texttt{correlatorEvenIntersecTwoQuadricsRecursionInTauCoord} $\{4,\{0,0,5,0,0,2,0,0,0,0,0,0\}\}$}
\]
returns $-624$, so $\langle \sfh_2,\sfh_2,\sfh_2,\sfh_2,\sfh_2,\epsilon_1,\epsilon_1\rangle=-624$. Similarly one can get the correlators in Lemma \ref{lem-correlatorsOfLength7-4dim}. The special correlator (\ref{eq-specialLength(n+3)Invariant-even(2,2)}) is denoted by \texttt{x} in this package. For example running
\begin{multline*}
\texttt{correlatorEvenIntersecTwoQuadricsRecursionInTauCoord}\\ \{6,\{0,0,0,0,0,0,0,0,2,2,2,2,2,2,2,0,0\}\}
\end{multline*}
returns $8\mbox{\texttt{x}}^2-2$. One can check Conjecture \ref{conj-unknownCorrelator-Even(2,2)-quadraticEquation} in this way.

\end{appendices}

\textsc{School of Mathematics, Sun Yat-sen University, Guangzhou 510275, P.R. China}

 \emph{E-mail address:}  huxw06@gmail.com 

\begin{thebibliography}{000}

\bibitem[ABPZ21]{ABPZ21}
Argüz, H.; Bousseau, P.; Pandharipande, R.; Zvonkine, D. Gromov--Witten Theory of Complete Intersections. arXiv:2109.13323.

\bibitem[Bea95]{Bea95}
Beauville, Arnaud.
Quantum cohomology of complete intersections. Mat. Fiz. Anal. Geom. 2 (1995), no. 3-4, 384–398.

\bibitem[Beh97]{Beh97}
Behrend, K. Gromov-Witten invariants in algebraic geometry. Invent. Math. 127 (1997), no. 3, 601–617.

\bibitem[BF97]{BF97}
Behrend, K.; Fantechi, B. 
The intrinsic normal cone. 
Invent. Math. 128 (1997), no. 1, 45–88.

\bibitem[Deb15]{Deb15}
Debarre, Olivier. Rational curves on hypersurfaces. Cabo 1 (2015): 12. 

\bibitem[Del73]{Del73}Deligne, P. Le th\'{e}or\`{e}me de Noether, Groupes de Monodromie en G\'{e}om\'{e}trie Alg\'{e}brique, Lecture
Notes in Math., Vol. 340, Springer-Verlag, Berlin and New York, 1973, pp. 328-340.

\bibitem[Dub98]{Dub98}
Dubrovin, Boris. Geometry and analytic theory of Frobenius manifolds, in Proceedings of the International Congress of Mathematicians, Vol. II (Berlin, 1998), 315–326.

\bibitem[Dub99]{Dub99}
Dubrovin, Boris.
Painlevé transcendents in two-dimensional topological field theory.  The Painlevé property, 287–412,
CRM Ser. Math. Phys., Springer, New York, 1999.

\bibitem[Ful98]{Ful98}
Fulton, William. Intersection theory. Second edition. Ergebnisse der Mathematik und ihrer Grenzgebiete. 3. Folge. A Series of Modern Surveys in Mathematics, 2. Springer-Verlag, Berlin, 1998. 

\bibitem[Get98]{Get98}
Getzler, E.
Topological recursion relations in genus 2. Integrable systems and algebraic geometry (Kobe/Kyoto, 1997), 73–106, World Sci. Publ., River Edge, NJ, 1998.


\bibitem[Giv96]{Giv96}
Givental, Alexander B. Equivariant Gromov-Witten invariants. Internat. Math. Res. Notices 1996, no. 13, 613–663.

\bibitem[Giv01]{Giv01}
Givental, Alexander B.
Semisimple Frobenius structures at higher genus.
Internat. Math. Res. Notices 2001, no. 23, 1265–1286.


\bibitem[Hu15]{Hu15}
Hu, Xiaowen. Big quantum cohomology of Fano complete intersections. arXiv:1501.03683v4.

\bibitem[Hum90]{Hum90}
Humphreys, James E. Reflection groups and Coxeter groups. Cambridge Studies in Advanced Mathematics, 29. Cambridge University Press, Cambridge, 1990. 

\bibitem[Jor1880]{Jor1880}
Jordan, Camille. Mémoire sur l'equivalence des formes, J. EC. Pol. XL VIII (1880),
112-150.

\bibitem[Kit93]{Kit93}
Kitaoka, Yoshiyuki.
Arithmetic of quadratic forms. 
Cambridge Tracts in Mathematics, 106. Cambridge University Press, Cambridge, 1993.

\bibitem[KM83]{KM83}
Kobayashi, Zenji; Morita, Jun.
Automorphisms of certain root lattices.
Tsukuba J. Math. 7 (1983), no. 2, 323–336.

\bibitem[KM94]{KM94}
 Kontsevich, M.; Manin, Yu. Gromov-Witten classes, quantum cohomology, and enumerative geometry. Comm. Math. Phys. 164 (1994), no. 3, 525–562.

\bibitem[Kuz08]{Kuz08}
 Kuznetsov, Alexander. Derived categories of quadric fibrations and intersections of quadrics. Adv. Math. 218 (2008), no. 5, 1340–1369.

\bibitem[Lew99]{Lew99} 
Lewis, James D. 
A survey of the Hodge conjecture. Second edition. Appendix B by B. Brent Gordon. CRM Monograph Series, 10. American Mathematical Society, Providence, RI, 1999.

\bibitem[LT98]{LT98}
Li, Jun; Tian, Gang.
Virtual moduli cycles and Gromov-Witten invariants of algebraic varieties.
J. Amer. Math. Soc. 11 (1998), no. 1, 119–174.

\bibitem[Man99]{Man99}
Manin, Yuri I. Frobenius manifolds, quantum cohomology, and moduli spaces. American Mathematical Society Colloquium Publications, 47. American Mathematical Society, Providence, RI, 1999.

\bibitem[Rei72]{Rei72}
Reid, Miles A. The complete intersection of two or more quadrics. PhD diss., University of Cambridge, 1972.


\bibitem[Ruan96]{Ruan96}
Ruan, Yongbin.
Topological sigma model and Donaldson-type invariants in Gromov theory.
Duke Math. J. 83 (1996), no. 2, 461–500.

\bibitem[Tel12]{Tel12}
Teleman, Constantin.
The structure of 2D semi-simple field theories.
Invent. Math. 188 (2012), no. 3, 525–588.

\bibitem[Zin08]{Zin08}
Zinger, Aleksey.
Standard versus reduced genus-one Gromov-Witten invariants.
Geom. Topol. 12 (2008), no. 2, 1203–1241.


\bibitem[Zin14]{Zin14}
Zinger, Aleksey. The genus 0 Gromov-Witten invariants of projective complete intersections. Geom. Topol. 18 (2014), no. 2, 1035–1114.

\end{thebibliography}
\end{document}